%% file: smc_bias_main.tex
\documentclass[final,onefignum,onetabnum,10pt,a4paper]{siamart171218}

\input{smc_bias_shared}

\begin{document}

\date{}

\maketitle

\begin{abstract}
In many applications, a state-space model depends on a parameter which needs to be inferred from
data in an online manner. In the maximum likelihood approach, this can be achieved using stochastic gradient search,
where the underlying gradient estimation is based on the optimal filter and the optimal filter derivative.
However, the optimal filter and its derivative are not analytically tractable
for a non-linear state-space model and need to be approximated numerically.
In \cite{poyiadjis&doucet&singh}, a particle approximation to this derivative
has been proposed, while the corresponding central limit theorem and $L_{p}$ error bounds
have been established in \cite{delmoral&doucet&singh2}.
We derive here bounds on the bias of this particle approximation.
Under mixing conditions, these bounds are uniform in time and inversely proportional to the number of particles.
\end{abstract}



\section{Introduction}\label{section0}

State-space models, also known as continuous-state hidden Markov models,
are a class of stochastic processes used to model complex time-series data
and stochastic dynamical systems.
A state-space model can be described as a latent discrete-time Markov process
observed only through noisy measurements of its states.
In this context, one of the most important problems is the optimal estimation
of the present state given the noisy observations of the present and past states.
This problem is known as optimal filtering.
For non-linear state-space models, the optimal filter does not typically admit a closed-form expression
and needs to be approximated numerically.
Numerous numerical methods for optimal filtering have been proposed and studied
in the literature --- see e.g., \cite{crisan&rozovskii} for a recent overview.
Among them, particle methods (also known as sequential Monte Carlo methods)
have gained significant attention.
Their convergence properties have been thoroughly investigated in a number of papers and books
--- see, e.g.,
\cite{cappe&moulines&ryden}, \cite{crisan&rozovskii}, \cite{delmoral1},
\cite{delmoral2}, \cite{douc&moulines&stoffer},
\cite{doucet&defreitas&gordon}.

In many scenarios of practical interest, a state-space model depends on a parameter
whose value needs to be estimated given a set of observations.
When the number of observations is very large,
it is desirable, for the sake of computational efficiency, to perform parameter estimation online.
In the maximum likelihood approach, this can be achieved using stochastic gradient search,
where the corresponding gradient estimator is a non-linear functional of the optimal filter and its derivative
--- see, e.g., \cite{kantas&doucet&singh&maciejowski&chopin},
\cite{legland&mevel2}, \cite{olsson&westerborn}, \cite{poyiadjis&doucet&singh}.
Since the optimal filter derivative is analytically intractable for non-linear
state-space models, it needs to be approximated numerically.
To the best of our knowledge, only the particle approximations to the optimal filter derivative
proposed in \cite{olsson&westerborn} and \cite{poyiadjis&doucet&singh} are numerically stable.
In the particle estimator proposed in \cite{olsson&westerborn},
the average iteration complexity is linear in the number of particles,
while the iteration running times are random.
In \cite{olsson&westerborn}, concentration inequalities and a central limit theorem have been shown
for this scheme.
The particle estimator proposed in \cite{poyiadjis&doucet&singh} has quadratic iteration complexity,
deterministic iteration running times and lower variance than the estimator in \cite{olsson&westerborn}.
In \cite{delmoral&doucet&singh2},
$L_{p}$ error bounds and a central limit theorem have been established for the scheme proposed
in \cite{poyiadjis&doucet&singh}.

In this paper, we analyze the bias of the particle approximation to the optimal filter derivative
proposed in \cite{poyiadjis&doucet&singh}. Using the stability properties of the optimal filter and its derivative,
we derive bounds on this bias in terms of the number of particles.
These bounds cover several classes of state-space models met in practice.
Moreover, under mixing conditions, these bounds are uniform in time and
inversely proportional to the number of particles.
To the best of our knowledge, the results presented here are the first results on
the bias of the particle approximation to the optimal filter derivative proposed in \cite{poyiadjis&doucet&singh}.
They are also one of the first and most important stepping stones to analyze the asymptotic properties of online maximum likelihood
estimation in non-linear state-space models --- see \cite{tadic&doucet2}.

The rest of this paper is organized as follows.
In Section \ref{section1}, we define the optimal filter derivative and its particle approximation.
In the same section, we present the main results of the paper.
These results are proved in Sections \ref{section2} -- \ref{section4}.

\section{Main Results}\label{section1}

\subsection{State-Space Models, Optimal Filter and Optimal Filter Derivative}\label{ssection1.1}

To define state-space models and state the optimal filtering problem, we use the following notation.
For a set ${\cal Z}$ in a finite dimensional space,
${\cal B}({\cal Z} )$ denotes the collection of Borel subsets of ${\cal Z}$.
$d_{x}\geq 1$ and $d_{y}\geq 1$ are integers,
while ${\cal X}\in{\cal B}(\mathbb{R}^{d_{x}} )$ and
${\cal Y}\in{\cal B}(\mathbb{R}^{d_{y}} )$.
Let $(\Omega,{\cal F},P)$ be a probability space.
A state-space model can be described as
an ${\cal X}\times{\cal Y}$-valued stochastic process
$\left\{ (X_{n},Y_{n}) \right\}_{n\geq 0}$ defined on $(\Omega,{\cal F},P)$,
where the process $\{X_{n} \}_{n\geq 0}$ is unobservable
and any information on $\{X_{n} \}_{n\geq 0}$ is only available through
the observation process $\{Y_{n} \}_{n\geq 0}$.
In this context, random variables $X_{n}$ and $Y_{n}$ are (respectively) called the state
and observation at discrete-time $n$,
while sets ${\cal X}$ and ${\cal Y}$ are (respectively) referred to as the state and observation spaces.
One of the most important problems related to state-space models is the estimation of the states
$X_{n}$ and $X_{n+1}$ given observations $Y_{0:n}:=(Y_0,...,Y_n)$.
This problem is known as filtering.

In the Bayesian approach,
the estimation of states $X_{n}$ and $X_{n+1}$ given $Y_{0:n}$
is based on the optimal filtering distributions
$P(X_{n}\in dx_{n}|Y_{0:n} )$ and $P(X_{n+1}\in dx_{n+1}|Y_{0:n})$.
In practice, the filtering distributions are usually evaluated using
some approximate models.
In this paper, we assume that the model
$\left\{ (X_{n},Y_{n}) \right\}_{n\geq 0}$
can be accurately approximated by a parametric family of non-linear state-space models.
To specify such a family, we rely on the following notation:
$d\geq 1$ is an integer, while $\Theta\subseteq\mathbb{R}^{d}$ is an open set.
$\mu(dx)$ and $\nu(dy)$ are positive measures on ${\cal X}$ and ${\cal Y}$ (respectively).
$p_{\theta}(x'|x)$ and $q_{\theta}(y|x)$ are Borel-measurable functions which map
$\theta\in\Theta$, $x,x'\in{\cal X}$, $y\in{\cal Y}$ to $[0,\infty )$
and are probability densities in $x'$, $y$ with respect to $\mu(dx)$, 
$\nu(dy)$. 
$\xi_{\theta}(dx)$ is a parameterized probability measure on ${\cal X}$,
i.e., $\xi_{\theta}(B)$ maps $\theta\in\Theta$, $B\in{\cal B}({\cal X} )$
to $[0,1]$ and is a probability measure in $B$ and
Borel-measurable in $\theta$.
With this notation, we can define a parametric family of state-space models
as an ${\cal X}\times{\cal Y}$-valued stochastic process
$\left\{ \big( X_{n}^{\theta}, Y_{n}^{\theta} \big) \right\}_{n\geq 0}$
which is defined on $(\Omega,{\cal F},P)$,
parameterized by $\theta\in\Theta$ 
and satisfies
\begin{align*}
	&
	P\left( \big( X_{0}^{\theta}, Y_{0}^{\theta} \big) \in B \right)
	=
	\int\int I_{B}(x,y) q_{\theta}(y|x) \nu(dy)\xi_{\theta}(dx),
	\\
	&
	P\left(\left. \big( X_{n+1}^{\theta}, Y_{n+1}^{\theta} \big) \in B
	\right| X_{0:n}^{\theta}, Y_{0:n}^{\theta}\right)
	=
	\int\int I_{B}(x,y) q_{\theta}(y|x) p_{\theta}\big(x|X_{n}^{\theta} \big)
	\nu(dy)\mu(dx)
\end{align*}
almost surely for each $\theta\in\Theta$, $B\in{\cal B}({\cal X}\times{\cal Y})$,
$n\geq 0$.\footnote
{To evaluate the values of $\theta$ for which
$\big\{ (X_{n}^{\theta}, Y_{n}^{\theta} ) \big\}_{n\geq 0}$
provides the best approximation to
$\left\{ (X_{n}, Y_{n} ) \right\}_{n\geq 0}$,
we usually rely on the maximum likelihood principle.
For further details on maximum likelihood estimation in state-space and hidden Markov models,
see, e.g., \cite{cappe&moulines&ryden}, \cite{douc&moulines&stoffer}.  }

Throughout the paper,
we assume that $p_{\theta}(x'|x)$ and $q_{\theta}(y|x)$ are differentiable in $\theta$
for each $\theta\in\Theta$, $x,x'\in{\cal X}$, $y\in{\cal Y}$.
To show how the filtering distribution and its derivative are computed using
the approximate model
$\big\{ \big( X_{n}^{\theta}, Y_{n}^{\theta} \big) \big\}_{n\geq 0}$,
we use the following notation.
If $k\geq 1$ is an integer, ${\cal Z}$ is a finite dimensional space and
${\zeta}(dz)$ is a $k$-dimensional signed vector measure on ${\cal Z}$,
then $\langle\zeta\rangle$ denotes the quantity
$\langle\zeta\rangle = \zeta({\cal Z} )$.
$w_{\theta}(x)$ a Borel-measurable function mapping $\theta\in\Theta$, $x\in{\cal X}$
to $\mathbb{R}^{d}$.
$r_{\theta,\boldsymbol y}^{n}(x'|x)$ and $t_{\theta,\boldsymbol y}^{n}(x'|x)$
are the functions defined by
\begin{align}\label{1.1}
	r_{\theta,\boldsymbol y}^{n}(x'|x)
	=
	p_{\theta}(x'|x) q_{\theta}(y_{n-1}|x),
	\;\;\;\;\;
	t_{\theta,\boldsymbol y}^{n}(x'|x)
	=
	\nabla_{\theta}\log\left( r_{\theta,\boldsymbol y}^{n}(x'|x) \right)
\end{align}
for $\theta\in\Theta$, $x,x'\in{\cal X}$, $n\geq 1$
and a sequence $\boldsymbol y = \{y_{n}\}_{n\geq 0}$ in ${\cal Y}$.
$r_{\theta,\boldsymbol y}^{m:n}(x_{m:n} )$ and $t_{\theta,\boldsymbol y}^{m:n}(x_{m:n} )$
are the functions defined by
\begin{align}
	\label{1.3}
	r_{\theta,\boldsymbol y}^{m:m}(x_{m:m} )
	=
	1,
	&
	\;\;\;\;\;
	r_{\theta,\boldsymbol y}^{m:n}(x_{m:n} )
	=
	\prod_{k=m+1}^{n} r_{\theta,\boldsymbol y}^{k}(x_{k}|x_{k-1} ),
	\\
	t_{\theta,\boldsymbol y}^{m:m}(x_{m:m} )
	=
	w_{\theta}(x_{m} ),
	&\;\;\;\;\;
	t_{\theta,\boldsymbol y}^{m:n}(x_{m:n} )
	=
	w_{\theta}(x_{m} )
	+
	\sum_{k=m+1}^{n} t_{\theta,\boldsymbol y}^{k}(x_{k}|x_{k-1} )
	\nonumber
\end{align}
for $x_{m},\dots,x_{n}\in{\cal X}$, $n>m\geq 0$.
$\mathbb{R}_{\theta,\boldsymbol y}^{n}(dx_{0:n} )$ and
$\mathbb{T}_{\theta,\boldsymbol y}^{n}(dx_{0:n} )$ are the measures defined by
\begin{align*}
	&
	\mathbb{R}_{\theta,\boldsymbol y}^{n}(A)
	=
	\int_{A}
	r_{\theta,\boldsymbol y}^{0:n}(x_{0:n} )
	(\xi_{\theta}\times\mu^{n} )(dx_{0:n} ),
	\\
	&
	\mathbb{T}_{\theta,\boldsymbol y}^{n}(A)
	=
	\int_{A}
	t_{\theta,\boldsymbol y}^{0:n}(x_{0:n} ) r_{\theta,\boldsymbol y}^{0:n}(x_{0:n} )
	(\xi_{\theta}\times\mu^{n} )(dx_{0:n} )
\end{align*}
for $A\in{\cal B}({\cal X}^{n+1} )$, $n\geq 1$,
where $\mu^{n}(d_{1:n})=\mu(dx_{1})\cdots\mu(dx_{n})$
and $(\xi_{\theta}\times\mu^{n} )(dx_{0:n} )=\xi_{\theta}(dx_{0})\mu^{n}(dx_{1:n})$.
$\mathbb{P}_{\theta,\boldsymbol y}^{n}(dx_{0:n} )$ and
$P_{\theta,\boldsymbol y}^{n}(dx)$ are the measures defined for $B\in{\cal B}({\cal X} )$ by
\begin{align}\label{1.1'}
	\mathbb{P}_{\theta,\boldsymbol y}^{n}(A)
	=
	\frac{\mathbb{R}_{\theta,\boldsymbol y}^{n}(A) }
	{\big\langle\mathbb{R}_{\theta,\boldsymbol y}^{n}\big\rangle },
	&\;\;\;\;\;
	P_{\theta,\boldsymbol y}^{n}(B)
	=
	\mathbb{P}_{\theta,\boldsymbol y}^{n}({\cal X}^{n}\times B ).
\end{align}
$\mathbb{Q}_{\theta,\boldsymbol y}^{n}(dx_{0:n} )$ and
$Q_{\theta,\boldsymbol y}^{n}(dx)$ are the measures defined by
\begin{align}\label{1.3'}
	\mathbb{Q}_{\theta,\boldsymbol y}^{n}(A)
	=
	\frac{\mathbb{T}_{\theta,\boldsymbol y}^{n}(A) }
	{\big\langle\mathbb{R}_{\theta,\boldsymbol y}^{n}\big\rangle }
	-
	\mathbb{P}_{\theta,\boldsymbol y}^{n}(A)
	\frac{\big\langle\mathbb{T}_{\theta,\boldsymbol y}^{n}\big\rangle }
	{\big\langle\mathbb{R}_{\theta,\boldsymbol y}^{n}\big\rangle },
	&\;\;\;\;\;
	Q_{\theta,\boldsymbol y}^{n}(B)
	=
	\mathbb{Q}_{\theta,\boldsymbol y}^{n}({\cal X}^{n}\times B ).
\end{align}

All results presented in this paper are based on Assumptions \ref{a1} -- \ref{a3}
(see Subsection \ref{ssection1.3}, below).
Using elementary calculus, it can be shown that the functions and measures defined above
are well-defined under these assumptions. 
It can also be verified 
\begin{align}\label{1.301}
	P_{\theta,\boldsymbol y}^{n}(B)
	=
	P\left(X_{n}^{\theta}\in B| Y_{0:n-1}^{\theta}=y_{0:n-1} \right).
\end{align}
Hence, $P_{\theta,\boldsymbol y}^{n}(dx)$ is the optimal filter (i.e., one-step predictor) for the model
\linebreak
$\left\{ (X_{n}^{\theta}, Y_{n}^{\theta} ) \right\}_{n\geq 0}$.

Let $\lambda(dx)$ be a finite positive measure on ${\cal X}$.
Moreover, let $l_{\theta}(x)$ be a Borel-measurable function
which maps $\theta\in\Theta$, $x\in{\cal X}$ to $(0,\infty )$
and satisfies the following conditions:
(i) $l_{\theta}(x)$ is differentiable in $\theta$
for each $\theta\in\Theta$, $x\in{\cal X}$, and
(ii) $\|\nabla_{\theta}l_{\theta}(x) \|$ is uniformly bounded in
$(\theta,x)$ on $\Theta\times{\cal X}$.
Suppose that Assumptions \ref{a1} -- \ref{a3} hold and
that $\xi_{\theta}(dx)$, $w_{\theta}(x)$ are of the form
\begin{align}\label{1.305}
	\xi_{\theta}(B) = \int_{B} l_{\theta}(x')\lambda(dx'),
	\;\;\;\;\;
	w_{\theta}(x) = \nabla_{\theta}\log\left(l_{\theta}(x) \right).
\end{align}
Then, it is straightforward to verify 
\begin{align}\label{1.303}
	Q_{\theta,\boldsymbol y}^{n}(B)
	=
	\nabla_{\theta} P\left(X_{n}^{\theta}\in B| Y_{0:n-1}^{\theta}=y_{0:n-1} \right). 
\end{align}
Thus, $Q_{\theta,\boldsymbol y}^{n}(B)$ is the optimal filter derivative
--- see \cite[Section 2]{delmoral&doucet&singh2} for further details.

\subsection{Particle Approximation to Optimal Filter Derivative}\label{ssection1.2}

Unless the model
$\big\{ \big( X_{n}^{\theta}, Y_{n}^{\theta} \big) \big\}_{n\geq 0}$
is linear Gaussian (or the state-space ${\cal X}$ has finitely many elements),
$P_{\theta,\boldsymbol y}^{n}(dx)$
and its gradient $Q_{\theta,\boldsymbol y}^{n}(dx)$ do not admit
closed-form expressions and need to be approximated numerically.

For a given $\theta \in \Theta$, the particle method proposed in \cite{poyiadjis&doucet&singh} approximates
$P_{\theta,\boldsymbol y}^{n}(dx)$ and
$Q_{\theta,\boldsymbol y}^{n}(dx)$
respectively by the empirical distributions
\begin{align}\label{1.31}
	\hat{\xi}_{n}^{\theta}(dx)
	=
	\frac{1}{N} \sum_{i=1}^{N} \delta_{\hat{X}_{n,i}^{\theta} }(dx),
	\;\;\;
	\hat{\zeta}_{n}^{\theta}(dx)
	=
	\frac{1}{N} \sum_{i=1}^{N}
	\left( W_{n,i}^{\theta} - \frac{1}{N} \sum_{j=1}^{N} W_{n,j}^{\theta} \right)
	\delta_{\hat{X}_{n,i}^{\theta} }(dx).
\end{align}
Here, $N\geq 2$ is a fixed integer and $\left\{ W_{n,i}^{\theta}: n\geq 0, 1\leq i\leq N \right\}$ are random vectors
generated through the recursion
\begin{align}\label{1.33}
	W_{n+1,i}^{\theta}
	=&
	\frac{\sum_{j=1}^{N} \left(
	p_{\theta}\big(\hat{X}_{n+1,j}^{\theta} |\hat{X}_{n,j}^{\theta} \big)
	\nabla_{\theta} q_{\theta}\big(Y_{n}|\hat{X}_{n,j}^{\theta} \big)
	+
	\nabla_{\theta} p_{\theta}\big(\hat{X}_{n+1,j}^{\theta} |\hat{X}_{n,j}^{\theta} \big)
	q_{\theta}\big(Y_{n}|\hat{X}_{n,j}^{\theta} \big)
	\right) }
	{\sum_{j=1}^{N} p_{\theta}\big(\hat{X}_{n+1,j}^{\theta} |\hat{X}_{n,j}^{\theta} \big)
	q_{\theta}\big(Y_{n}|\hat{X}_{n,j}^{\theta} \big) }
	\\ \nonumber
	&+
	\frac{\sum_{j=1}^{N} p_{\theta}\big(\hat{X}_{n+1,j}^{\theta} |\hat{X}_{n,j}^{\theta} \big)
	q_{\theta}\big(Y_{n}|\hat{X}_{n,j}^{\theta} \big) W_{n,j}^{\theta} }
	{\sum_{j=1}^{N} p_{\theta}\big(\hat{X}_{n+1,j}^{\theta} |\hat{X}_{n,j}^{\theta} \big)
	q_{\theta}\big(Y_{n}|\hat{X}_{n,j}^{\theta} \big) },
\end{align}
where $\big\{ \hat{X}_{n,i}^{\theta}: n\geq 0, 1\leq i\leq N \big\}$ are random samples called particles
generated through 
\begin{align}\label{1.35}
	\hat{X}_{n+1,i}^{\theta}
	\sim
	\frac{\sum_{j=1}^{N} p_{\theta}\big(x|\hat{X}_{n,j}^{\theta} \big)
	q_{\theta}\big(Y_{n}|\hat{X}_{n,j}^{\theta} \big) \mu(dx) }
	{\sum_{j=1}^{N} q_{\theta}\big(Y_{n}|\hat{X}_{n,j}^{\theta} \big) }.
\end{align}
In recursion (\ref{1.33}), $\big\{ W_{0,i}^{\theta}: 1\leq i\leq N \big\}$
are selected as $W_{0,i}^{\theta} = w_{\theta}\big( \hat{X}_{0,i}^{\theta} \big)$.
In recursion (\ref{1.35}),
$\big\{ \hat{X}_{n+1,i}^{\theta}: 1\leq i\leq N \big\}$ are sampled independently
from one another. In the same recursion, $\big\{ \hat{X}_{0,i}^{\theta}: 1\leq i\leq N \big\}$
are sampled from $\xi_{\theta}(dx)$ independently one from another and independently from $Y_{0}$.

\begin{vremark}
Let $w_{\theta,\boldsymbol y}^{n}(x)$ be the function defined by
\begin{align*}
	w_{\theta,\boldsymbol y}^{n}(x)
	=
	E\left(
	t_{\theta,\boldsymbol y}^{0:n}(X_{0:n}^{\theta} )
	|X_{n}^{\theta}=x,Y_{0:n-1}^{\theta}=y_{0:n-1}
	\right)
\end{align*}
for $\theta\in\Theta$, $x\in{\cal X}$, $n\geq 1$ and a sequence
$\boldsymbol y = \{y_{n} \}_{n\geq 0}$ in ${\cal Y}$.
Then, we have $W_{n,i}^{\theta} = w_{\theta,\boldsymbol Y}^{n}(\hat{X}_{n,i}^{\theta} )$
for each $1\leq i\leq N$, $n\geq 1$,
where $\boldsymbol Y = \{Y_{n} \}_{n\geq 0}$.
The equivalence between this representation of random vector $W_{n,i}^{\theta}$
and recursion (\ref{1.33}) is shown and discussed
in \cite[Section 3]{delmoral&doucet&singh2}.
Recursion (\ref{1.33}) is derived in \cite[Section 2.2]{poyiadjis&doucet&singh}.
\end{vremark}

\subsection{Bias of Particle Approximation to Optimal Filter Derivative}\label{ssection1.3}

We analyze here the bias of the particle approximations (\ref{1.31}).
The analysis relies on the following notation.
For $z\in\mathbb{R}^{k}$, $k\geq 1$,
$\|z\|$ denotes the $l_{\infty}$ norm of $z$.
If $\psi_{\theta}(x)$ is a Borel-measurable function mapping
$\theta\in\Theta$, $x\in{\cal X}$ to $\mathbb{R}^{k}$,
then $\|\psi_{\theta}\|$ denotes the $L_{\infty}$ norm of $\psi_{\theta}(x)$
in $x$, i.e., $\|\psi_{\theta} \| = \sup_{x\in{\cal X} } \|\psi_{\theta}(x) \|$.
If $\varphi:{\cal X}\rightarrow\mathbb{R}$ is a Borel-measurable function
and $\zeta(dx)$ is a $k$-dimensional signed vector measure on ${\cal X}$,
then $\zeta(\varphi)$ denotes the integral $\zeta(\varphi) = \int \varphi(x)\zeta(dx)$.

The analysis carried out in this paper relies on the following assumptions.

\begin{assumption}\label{a1}
There exists a real number $\varepsilon\in(0,1)$ such that
\begin{align*}
	\varepsilon \leq p_{\theta}(x'|x) \leq \frac{1}{\varepsilon},
	\;\;\;\;\;
	\varepsilon \leq q_{\theta}(y|x) \leq \frac{1}{\varepsilon}
\end{align*}
for all $\theta\in\Theta$, $x,x'\in{\cal X}$, $y\in{\cal Y}$.
\end{assumption}

\begin{assumption}\label{a2}
There exists a real number $K\in[1,\infty)$ such that
\begin{align*}
	\max\{\|\nabla_{\theta} p_{\theta}(x'|x) \|,
	\|\nabla_{\theta} q_{\theta}(y|x) \| \}
	\leq K
\end{align*}
for all $\theta\in\Theta$, $x,x'\in{\cal X}$, $y\in{\cal Y}$.
\end{assumption}

\begin{assumption}\label{a3}
$\|w_{\theta} \|=\sup_{x\in{\cal X} } \|w_{\theta}(x) \| < \infty$
for all $\theta\in\Theta$.
\end{assumption}

Assumption \ref{a1} is a standard strong mixing condition
and is a crucial ingredient of many results on
optimal filtering and statistical inference in state-space and hidden Markov models
--- see e.g., \cite{cappe&moulines&ryden}, \cite{delmoral&guionnet}, \cite{delmoral&doucet&singh2},
\cite{douc&moulines&ryden}, \cite{legland&mevel}, \cite{legland&oudjane},
\cite{oudjane&rubenthaler}
\cite{tadic&doucet1}.
Together with Assumption \ref{a2}, it ensures that
the optimal filter and its derivative forget initial conditions exponentially fast
--- see Proposition \ref{proposition3.2}, Section \ref{section3}.
This assumption, together with Assumptions \ref{a2} and \ref{a3}, also ensures
the stability of particle approximations $(\ref{1.31})$
--- see Proposition \ref{proposition4.1}, Section \ref{section4}.
Assumption \ref{a1} is restrictive as it implicitly requires the state and observation spaces ${\cal X}$ and ${\cal Y}$
to be bounded.

Let $\boldsymbol Y$ denote stochastic process $\{Y_{n} \}_{n\geq 0}$,
i.e., $\boldsymbol Y = \{Y_{n} \}_{n\geq 0}$.
The main results of our paper are stated in the next theorem.

\begin{theorem}\label{theorem1.1}
Let $\theta$ be any element of $\Theta$, 
while $\boldsymbol y = \{y_{n} \}_{n\geq 0}$ is any sequence in ${\cal Y}$. 
Moreover, let $\varphi:{\cal X}\rightarrow[-1,1]$ be any Borel-measurable function, 
while $n$ is any positive integer. 

(i) Suppose that Assumption \ref{a1} holds.
Then, there exists a real number $L\in[1,\infty)$
(independent of $N$, $\theta$, $\boldsymbol y$, $\varphi(x)$, $n$ 
and depending only on $\varepsilon$) such that
\begin{align}\label{t1.1.1*}
	\left|
	E\left(\left.
	\hat{\xi}_{n}^{\theta}(\varphi)
	-
	P_{\theta,\boldsymbol Y}^{n}(\varphi)
	\right|
	\boldsymbol Y = \boldsymbol y
	\right)
	\right|
	\leq
	\frac{L}{N}. 
\end{align}

(ii) Suppose that Assumptions \ref{a1} -- \ref{a3} hold.
Then, there exist real numbers $\rho\in(0,1)$, $M\in[1,\infty)$
(independent of $N$, $\theta$, $\boldsymbol y$, $\varphi(x)$, $n$ 
and depending only on $\varepsilon$, $d$, $K$) such that
\begin{align}\label{t1.1.3*}
	\left\|
	E\left(\left.
	\hat{\zeta}_{n}^{\theta}(\varphi)
	-
	Q_{\theta,\boldsymbol Y}^{n}
	(\varphi)
	\right|
	\boldsymbol Y = \boldsymbol y
	\right)
	\right\|
	\leq
	\frac{M(1 + \rho^{n} \|w_{\theta}\| ) }{N}. 
\end{align}
\end{theorem}

The proof of Theorem \ref{theorem1.1} is provided in Section \ref{section5}
--- see Proposition \ref{proposition5.1}.

The empirical measures $\hat{\xi}_{n}^{\theta}(dx)$ and $\hat{\zeta}_{n}^{\theta}(dx)$
are estimators of the optimal predictor
$P_{\theta, \boldsymbol Y}^{n}(dx)$
and its gradient
$Q_{\theta, \boldsymbol Y}^{n}(dx)$.
Hence, the conditional expectations in (\ref{t1.1.1*}), (\ref{t1.1.3*})
can be viewed as the bias of particle approximations (\ref{1.31}) for which 
Theorem \ref{theorem1.1} provides bounds. These bounds are inversely proportional to $N$ and uniform in discrete-time $n$
as $\rho^{n}\leq 1$. They depend on $\left\{ \big( X_{n}^{\theta}, Y_{n}^{\theta} \big) \right\}_{n\geq 0}$
through constants $\rho$, $L$, $M$ and the initial conditions in recursion (\ref{1.33})
(through $\|w_{\theta}\|$).

Due to their practical and theoretical importance,
particle methods have extensively been studied
in a number of papers and books
--- see e.g., \cite{andrieu&doucet&singh&tadic},
\cite{cappe&moulines&ryden}, \cite{crisan&rozovskii},
\cite{delmoral1}, \cite{delmoral2},
\cite{douc&moulines&stoffer},
\cite{doucet&defreitas&gordon},
\cite{douc&moulines&stoffer} -- \cite{kantas&doucet&singh&maciejowski&chopin}.
Within a broader analysis of the propagation of chaos in Feynman-Kac models,
the bias of particle approximations to the optimal filter has been addressed in
\cite{delmoral1} -- \cite{delmoral&jacob&lee&murray&peters},
\cite{delmoral&jasra}.
Under conditions similar or identical to Assumption \ref{a1},
the results of
\cite{delmoral1} -- \cite{delmoral&jacob&lee&murray&peters},
\cite{delmoral&jasra}
lead to Part (i) of Theorem \ref{theorem1.1}.\footnote
{Although Part (i) of Theorem \ref{theorem1.1} is
a particular case in the analysis carried out in
\cite{delmoral1}, \cite{delmoral2}, \cite{delmoral&doucet&peters},
\cite{delmoral&doucet&singh1}, \cite{delmoral&jacob&lee&murray&peters},
\cite{delmoral&jasra},
we include it in the main results for the following reasons:
(i) $\xi_{n}^{\theta}(dx)$ is an integral part
of the particle approximation (\ref{1.31}) -- (\ref{1.35}),
(ii) the bound (\ref{t1.1.1*}) is an essential prerequisite for
Part (ii) of Theorem \ref{theorem1.1},
and (iii) the proof of Part (i) of Theorem \ref{theorem1.1} presented here seems
more direct than the analysis carried out in
\cite{delmoral1}, \cite{delmoral2}, \cite{delmoral&doucet&peters},
\cite{delmoral&doucet&singh1}, \cite{delmoral&jacob&lee&murray&peters},
\cite{delmoral&jasra}. }
As opposed to particle approximations to the optimal filter,
the optimal filter derivative and its particle approximations have attracted
much less attention.
Part (ii) of Theorem \ref{theorem1.1} fills this gap in the literature
on optimal filtering and particle methods.
To the best of our knowledge, Part (ii) of Theorem \ref{theorem1.1} is the first result on
the bias of the particle approximation (\ref{1.31}) -- (\ref{1.35}).
In \cite{tadic&doucet2}, we use this result, together with the results of \cite{tadic&doucet4},
to analyze the asymptotic behavior of
recursive maximum likelihood estimation in non-linear state-space models.

\section{Results Related to Empirical Measures}\label{section2}

In this section, we present an auxiliary result on
the ratio of integrals approximated using empirical measures.
This result is a crucial ingredient in the proof of Lemma \ref{lemma5.2}
which itself is a corner-stone of the main results --- see Proposition \ref{proposition5.1}.
This result has already appeared in
\cite[Lemma A.1]{delmoral&doucet&singh2} and \cite[Lemmas B.3, B.4]{delmoral&jasra}. 
For completeness, a proof is provided in the supplementary material (Section SM1).

We use the following additional notation. 
${\cal Z}$ is a finite dimensional space, 
while $\xi(dz)$ is a probability measure on ${\cal Z}$.
$\{Z_{k} \}_{k\geq 1}$ are independent ${\cal Z}$-valued random variables
which are defined on a probability space $(\Omega, {\cal F}, P)$ and distributed according to $\xi(dz)$
(i.e., $P(Z_{k}\in B ) = \xi(B)$
for each $B\in{\cal B}({\cal Z} )$).
$\xi_{k}(dz)$ is the empirical measure defined for $k\geq 1$ by
\begin{align*}
	\xi_{k}(B)
	=
	\frac{1}{k} \sum_{i=1}^{k} \delta_{Z_{i} }(B).
\end{align*}
\begin{proposition}\label{proposition2.1}
Let $f:{\cal Z}\rightarrow \mathbb{R}$ and $g:{\cal Z}\rightarrow (0,\infty )$
be Borel-measurable functions such that
\begin{align*}
	\sup_{z\in{\cal Z} } |f(z)| < \infty,
	\;\;\;\;\;
	\sup_{z\in{\cal Z} } g(z) < \infty,
	\;\;\;\;\;
	\inf_{z\in{\cal Z} } g(z) > 0.
\end{align*}
Then, we have
\begin{align}\label{p2.1.3*}
	\left|
	E\left(
	\frac{\xi_{k}(f) }{\xi_{k}(g) } \right)
	-
	\frac{\xi(f) }{\xi(g) }
	\right|
	\leq
	\frac{2\alpha\beta^{2} }{k},
	\;\;\;\;\;
	\left(
	E\left(
	\left|
	\frac{\xi_{k}(f) }{\xi_{k}(g) }
	-
	\frac{\xi(f) }{\xi(g) }
	\right|^{2}
	\right)
	\right)^{1/2}
	\leq
	\frac{2\alpha\beta}{\sqrt{k} }
\end{align}
for any $k\geq 1$,
where $\alpha$, $\beta$ are defined by
\begin{align}\label{p2.1.1*}
	\alpha
	=
	\sup_{z',z''\in{\cal Z} }
	\left|
	\frac{f(z')}{g(z')} - \frac{f(z'')}{g(z'')}
	\right|,
	\;\;\;\;\;
	\beta
	=
	\sup_{z',z''\in{\cal Z} }
	\frac{g(z')}{g(z'')}.
\end{align}
\end{proposition}

\section{Results Related to Stability of Optimal Filter and Its Derivative}\label{section3}

In this section, we present results on the stability properties of the optimal predictor
$P_{\theta,\boldsymbol y}^{n}(dx)$ and its gradient
$Q_{\theta,\boldsymbol y}^{m:n}(dx)$.
These results are prerequisites for the proof
of the main results --- see Lemmas \ref{lemma5.1}, \ref{lemma5.2} and
Proposition \ref{proposition5.1}.

The following additional notation is used here.
${\cal P}({\cal X} )$ is the collection of probability measures on ${\cal X}$,
while the set of Borel-measurable functions mapping ${\cal X}$ to $\mathbb{R}$
is denoted by ${\cal F}({\cal X})$.
${\cal M}_{p}({\cal X} )$ is the set of positive measures on ${\cal X}$,
while the collection of signed measures on ${\cal X}$ is denoted by ${\cal M}_{s}({\cal X} )$.
If $k\geq 1$ is an integer,
then ${\cal M}_{s}^{k}({\cal X} )$ is the set of $k$-dimensional signed vector measures on ${\cal X}$
and ${\cal F}^{k}({\cal X})$ is the collection of Borel-measurable functions
mapping ${\cal X}$ to $\mathbb{R}^{k}$.
For $\xi\in{\cal M}_{s}({\cal X})$,
$|\xi|(dx)$ and $\|\xi\|$ denote (respectively) the total variation
and the total variation norm of $\xi(dx)$.
For $\zeta\in{\cal M}_{s}^{k}({\cal X} )$,
$|\zeta|(dx)$ and $\|\zeta\|$ denote (respectively) the total variation
and the total variation norm of $\zeta(dx)$ induced by $l_{1}$ vector norm.\footnote
{If $\zeta\in{\cal M}_{s}^{k}({\cal X})$,
then $|\zeta|(dx) = \sum_{i=1}^{k} |e_{i}^{T}\zeta|(dx)$
and $\|\zeta\| = \sum_{i=1}^{k} \|e_{i}^{T}\zeta\|$,
where $e_{i}$ is the $i$-th standard unit vector in $\mathbb{R}^{k}$. }
If $H(\zeta)$ is a function mapping $\zeta\in{\cal M}_{s}^{k}$ to ${\cal M}_{s}^{l}$,
then $H(\zeta)(B)$ stands for the measure of $B\in{\cal B}({\cal X})$
with respect to $H(\zeta)$ .
Moreover, if $\xi\in{\cal M}_{s}({\cal X})$, $\zeta\in{\cal M}_{s}^{k}({\cal X})$
and $R(x,dx')$, $S(x,dx')$ are integral operators from
${\cal F}({\cal X})$ to ${\cal F}({\cal X})$, ${\cal F}^{k}({\cal X})$ (respectively),
then $(\zeta R)(dx)$ and $(\xi S)(dx)$ denote the elements of
${\cal M}_{s}^{k}({\cal X})$ defined for $B\in{\cal B}({\cal X})$ by
\begin{align*}
	(\zeta R)(B) = \int_{B} R(x,B)\zeta(dx),
	\;\;\;\;
	(\xi S)(B) = \int_{B} S(x,B)\xi(dx).
\end{align*}
$s_{\theta,\boldsymbol y}^{m:n}(x_{m:n} )$ is the function defined by
\begin{align}\label{1.5}
	s_{\theta,\boldsymbol y}^{m:m}(x_{m:m} )=0,
	\;\;\;\;\;
	s_{\theta,\boldsymbol y}^{m:n}(x_{m:n} )
	=
	\nabla_{\theta}r_{\theta,\boldsymbol y}^{m:n}(x_{m:n} )
\end{align}
for $\theta\in\Theta$, $x_{m},\dots,x_{n}\in{\cal X}$,
$n>m\geq 0$ and a sequence $\boldsymbol y = \{y_{n} \}_{n\geq 0}$ in ${\cal Y}$
($r_{\theta,\boldsymbol y}^{m:n}(x_{m:n} )$ is specified in (\ref{1.3})).
$R_{\theta,\boldsymbol y}^{m:n}(x,dx')$ and
$S_{\theta,\boldsymbol y}^{m:n}(x,dx')$
are the integral operators from ${\cal F}({\cal X})$ to
${\cal F}({\cal X})$, ${\cal F}^{d}({\cal X})$ (respectively)
defined by
\begin{align}
	\label{1.7}
	R_{\theta,\boldsymbol y}^{m:m}(x,B)
	=
	\delta_{x}(B),
	\;\;\;\;\;
	&
	\begin{aligned}[t]
	R_{\theta,\boldsymbol y}^{m:n}(x,B)
	=
	&\int_{{\cal X}^{n-m}\times B} r_{\theta,\boldsymbol y}^{m:n}(x_{m:n} )
	(\delta_{x}\times\mu^{n-m})(dx_{m:n}),
	\end{aligned}
	\\
	\label{1.9}
	S_{\theta,\boldsymbol y}^{m:m}(x,B)
	=
	0,
	\;\;\;\;\;
	&
	\begin{aligned}[t]
	S_{\theta,\boldsymbol y}^{m:n}(x,B)
	=
	&\int_{{\cal X}^{n-m}\times B} s_{\theta,\boldsymbol y}^{m:n}(x_{m:n} )
	(\delta_{x}\times\mu^{n-m})(dx_{m:n}),
	\end{aligned}
\end{align}
where $(\delta_{x}\times\mu^{n-m})(dx_{m:n})=\delta_{x}(dx_{m})\mu(dx_{m+1})\cdots\mu(dx_{n})$.
$F_{\theta,\boldsymbol y}^{m:n}(\xi)$,
$G_{\theta,\boldsymbol y}^{m:n}(\xi,\zeta)$ and
$H_{\theta,\boldsymbol y}^{m:n}(\xi,\zeta)$
are the functions mapping $\xi\in{\cal P}({\cal X})$, $\zeta\in{\cal M}_{s}^{d}({\cal X})$
to ${\cal P}({\cal X})$, ${\cal M}_{s}^{d}({\cal X})$, ${\cal M}_{s}^{d}({\cal X})$ (respectively)
defined by
\begin{align}
	\label{1.21}
	F_{\theta,\boldsymbol y}^{m:m}(\xi)(B)
	=
	\xi(B),
	&\;\;\;
	F_{\theta,\boldsymbol y}^{m:n}(\xi)(B)
	=
	\frac{(\xi R_{\theta,\boldsymbol y}^{m:n} )(B) }
	{\big\langle \xi R_{\theta,\boldsymbol y}^{m:n} \big\rangle },
	\\
	\label{1.23}
	H_{\theta,\boldsymbol y}^{m:m}(\xi,\zeta)(B)
	\!=\!
	\zeta(B),
	&\;\;\;
	H_{\theta,\boldsymbol y}^{m:n}(\xi,\zeta)(B)
	\!=\!
	\frac{(\zeta R_{\theta,\boldsymbol y}^{m:n} )(B)
	+ (\xi S_{\theta,\boldsymbol y}^{m:n} )(B) }
	{\big\langle \xi R_{\theta,\boldsymbol y}^{m:n} \big\rangle },
	\\
	\label{1.25}
	G_{\theta,\boldsymbol y}^{m:m}(\xi,\zeta)(B)
	\!=\!
	\zeta(B),
	&\;\;\;
	G_{\theta,\boldsymbol y}^{m:n}(\xi,\zeta)(B)
	\!=\!
	H_{\theta,\boldsymbol y}^{m:n}(\xi,\zeta)(B)
	\!-\!
	F_{\theta,\boldsymbol y}^{m:n}(\xi)(B)
	\big\langle H_{\theta,\boldsymbol y}^{m:n}(\xi,\zeta) \big\rangle.
\end{align}

\begin{vremark}
For $\theta\in\Theta$, let $\zeta_{\theta}(dx)$ be the element of ${\cal M}_{s}^{d}({\cal X} )$
specified in (\ref{5.1}) (below).
Suppose that Assumptions \ref{a1} -- \ref{a3} hold.
Then, it easy to show
\begin{align}\label{3.301}
	P_{\theta,\boldsymbol y}^{n}(B)
	=
	F_{\theta,\boldsymbol y}^{0:n}(\xi_{\theta} )(B),
	\;\;\;\;\;
	Q_{\theta,\boldsymbol y}^{n}(B)
	=
	G_{\theta,\boldsymbol y}^{0:n}(\xi_{\theta},\zeta_{\theta} )(B)
\end{align}
for $n\geq 1$.
Given (\ref{1.301}), (\ref{1.303}), (\ref{3.301}),
$F_{\theta,\boldsymbol y}^{m:n}(\xi)$ and
$G_{\theta,\boldsymbol y}^{m:n}(\xi,\zeta)$
can be considered as a generalization of the optimal filter and its gradient.
In this context, $\xi$, $\zeta$ can be viewed as initial conditions
in the recursion generating the optimal filter and its gradient
(for further details, see \cite{legland&mevel}, \cite{tadic&doucet1}).
\end{vremark}

\begin{proposition}\label{proposition3.2}
Let $\theta$ be any element of $\Theta$, 
while $\boldsymbol y = \{y_{n} \}_{n\geq 0}$ is any sequence in ${\cal Y}$. 
Moreover, let $\xi,\xi'$ be any elements of ${\cal P}({\cal X})$,
while $\zeta,\zeta'$ are any elements of ${\cal M}_{s}^{d}({\cal X})$. 
Further to this, let $n$, $m$ be any integers satisfying $n\geq m\geq 0$. 

(i) Suppose that Assumption \ref{a1} holds.
Then, there exist real numbers $\rho_{1}\in(0,1)$, $C_{1}\in[1,\infty)$
(independent of $\theta$, $\boldsymbol y$, $\xi$, $\xi'$, $n$, $m$
and depending only on $\varepsilon$) such that
\begin{align}
	&\label{p3.2.3*}
	\left\| F_{\theta,\boldsymbol y}^{m:n}(\xi)
	-
	F_{\theta,\boldsymbol y}^{m:n}(\xi') \right\|
	\leq
	C_{1}\rho_{1}^{n-m},
	\;\;\;\;\;
	\frac{\big\langle \xi R_{\theta,\boldsymbol y}^{m:n} \big\rangle }
	{\big\langle \xi' R_{\theta,\boldsymbol y}^{m:n} \big\rangle }
	\leq
	C_{1}. 
\end{align}

(ii) Suppose that Assumptions \ref{a1} and \ref{a2} hold.
Then, there exist real numbers $\rho_{2}\in(0,1)$, $C_{2}\in[1,\infty)$
(independent of $\theta$, $\boldsymbol y$, $\xi$, $\xi'$, $\zeta$, $\zeta'$, $n$, $m$
and depending only on $\varepsilon$, $d$, $K$) such that
\begin{align}
	&\label{p3.2.5*}
	\begin{aligned}[t]
	\left\| G_{\theta,\boldsymbol y}^{m:n}(\xi,\zeta)
	-
	G_{\theta,\boldsymbol y}^{m:n}(\xi',\zeta') \right\|
	\leq &
	C_{2}\rho_{2}^{n-m} (1 + \|\zeta\| + \|\zeta'\| ),
	\end{aligned}
	\\
	&
	\label{p3.2.7*}
	\left\| H_{\theta,\boldsymbol y}^{m:n}(\xi,\zeta)\right\|
	\leq
	C_{2} (n-m+\|\zeta\| ). 
\end{align}
\end{proposition}

Proposition \ref{proposition3.2} is a relatively straightforward extension
of the results of \cite{legland&mevel}, \cite{legland&oudjane}, \cite{tadic&doucet1} to the optimal predictor and its gradient.
A detailed proof of the proposition is provided in the supplementary material
(Section SM2).

\section{Results Related to Stability of Particle Approximations}\label{section4}

In this section, we consider the particle approximation $\hat{\zeta}_{n}^{\theta}(dx)$
and its stability.
Using results on the (Dobrushin) ergodicity coefficient,
we show that the sequence $\big\{ \big\|\hat{\zeta}_{n}^{\theta} \big\| \big\}_{n\geq 0}$
is bounded uniformly in $\theta$.
The results presented here are prerequisites for the proof of
the main results --- see Lemma \ref{lemma5.2} and Proposition \ref{proposition5.1}.

To state the results of this section, additional notation needs to be introduced.
${\cal P}^{N}$ is the set of $N$-dimensional probability vectors.
${\cal P}^{N\times N}$ is the set of $N\times N$ (column) stochastic matrices
(i.e., $A\in{\cal P}^{N\times N}$ if and only if the columns of $A$ are elements of
${\cal P}^{N}$).
$e$ is the element of $\mathbb{R}^{N}$ whose all elements are one.
For $1\leq i\leq N$,
$e_{i}$ is the $i$-th standard unit vector in $\mathbb{R}^{N}$
(i.e., $e_{i}$ is the element of ${\cal P}^{N}$ whose $i$-th element is one).
For $z\in\mathbb{R}^{N}$,
$\|z\|_{1}$ and $\|z\|$ are (respectively) the $l_{1}$ and $l_{\infty}$ norm of $z$.
For $B\in\mathbb{R}^{d\times N}$,
$\|B\|$ is the $l_{\infty}$ norm of $B$
(i.e., $\|B\|$ is the maximum absolute value of the entries of $B$).
For $A\in{\cal P}^{N\times N}$,
$\tau(A)$ is the ergodicity coefficient of $A$, i.e.,
\begin{align}\label{4.7}
	\tau(A)
	=
	\frac{1}{2}
	\max_{1\leq j',j''\leq N} \sum_{i=1}^{N} |A_{i,j'} - A_{i,j''} |
	=
	1
	-
	\min_{1\leq j',j''\leq N}
	\sum_{i=1}^{N} \min\{A_{i,j'},A_{i,j''} \},
\end{align}
where $A_{i,j}$ is the $(i,j)$ entry of $A$
(see \cite[Section 15.2.1]{bremaud} for more details on the ergodicity coefficient and its equivalent forms).
$A_{n}^{\theta}$ and $B_{n}^{\theta}$ are (respectively)
the $N\times N$ and $d\times N$ random matrices defined by
\begin{align}\label{4.5}
	A_{n,i,j}^{\theta}
	=
	\frac{r_{\theta,\boldsymbol Y}^{n}\big(\hat{X}_{n,j}^{\theta} |\hat{X}_{n-1,i}^{\theta} \big) }
	{\sum_{k=1}^{N} r_{\theta,\boldsymbol Y}^{n}
	\big(\hat{X}_{n,j}^{\theta} |\hat{X}_{n-1,k}^{\theta} \big) },
	\;\;\;\;\;
	B_{n,j}^{\theta}
	=
	\frac{\sum_{k=1}^{N} \nabla_{\theta} r_{\theta,\boldsymbol Y}^{n}\big(\hat{X}_{n,j}^{\theta} |\hat{X}_{n-1,k}^{\theta} \big) }
	{\sum_{k=1}^{N} r_{\theta,\boldsymbol Y}^{n}
	\big(\hat{X}_{n,j}^{\theta} |\hat{X}_{n-1,k}^{\theta} \big) },
\end{align}
where $A_{n,i,j}^{\theta}$ is the $(i,j)$ entry of $A_{n}^{\theta}$
and $B_{n,j}^{\theta}$ is the $j$-th column of $B_{n}^{\theta}$.
$r_{\theta,\boldsymbol y}^{n}(x'|x)$ and $\boldsymbol Y$
are specified in (\ref{1.1}) and Subsection \ref{ssection1.3} (respectively).
$V_{n}^{\theta}$, $W_{n}^{\theta}$ and $V_{n,i}^{\theta}$
are the $d\times N$ random matrices and the $d$-dimensional random vector
defined by
\begin{align}\label{4.3}
	V_{n,i}^{\theta}
	=
	W_{n,i}^{\theta}
	-
	\frac{1}{N} \sum_{j=1}^{N} W_{n,j}^{\theta},
	\;\;\;
	V_{n}^{\theta}
	=
	\big( V_{n,1}^{\theta}, \dots, V_{n,N}^{\theta} \big),
	\;\;\;
	W_{n}^{\theta}
	=
	\big( W_{n,1}^{\theta}, \dots, W_{n,N}^{\theta} \big)
\end{align}
for $n\geq 0$.
Notice here that $V_{n,i}^{\theta}$ and $W_{n,i}^{\theta}$ are the $i$-th columns of
$V_{n}^{\theta}$ and $W_{n}^{\theta}$.
Then, it is easy to show
$A_{n}^{\theta}\in{\cal P}^{N\times N}$ and
\begin{align}\label{4.1}
	V_{n}^{\theta}
	=
	W_{n}^{\theta}
	\left( I - \frac{ee^{T} }{N} \right),
	\;\;\;\;\;
	W_{n+1}^{\theta} = W_{n}^{\theta} A_{n+1}^{\theta} + B_{n+1}^{\theta},
\end{align}
where  $I$ is the $N\times N$ unit matrix.

\begin{vremark}
Throughout this and subsequent sections, the following convention is applied.
Diacritic $\tilde{}$ is used to denote a locally defined quantity,
i.e., a quantity whose definition holds only within the proof where
the quantity appears.
\end{vremark}

\begin{proposition}\label{proposition4.1}
Let $\theta$ be any element of $\Theta$, 
while $n$ is any non-negative integer. 
Suppose that Assumptions \ref{a1} -- \ref{a3} hold.
Then, there exist real numbers $\rho_{3}\in(0,1)$, $C_{3}\in[1,\infty)$
(independent of $N$, $\theta$, $n$ and depending only on $\varepsilon$, $d$, $K$) such that
\begin{align}\label{l4.1.1*}
	\big\| \hat{\zeta}_{n}^{\theta} \big\|
	\leq
	C_{3} \left(1 + \rho_{3}^{n} \|w_{\theta} \| \right). 
\end{align}
\end{proposition}

\begin{proof}
Throughout the proof, the following notation is used.
$\rho_{3}$, $C_{3}$ are the real numbers defined by
$\rho_{3}=1-\varepsilon^{4}$, $C_{3}=8Kd\varepsilon^{-7}$
($\varepsilon$, $K$ are specified in Assumptions \ref{a1} and \ref{a2}).
$\tilde{A}_{k,l}^{\theta}$ is the matrix defined by
\begin{align}
	\tilde{A}_{k,k}^{\theta} = I,
	\;\;\;\;\;
	\tilde{A}_{k,l}^{\theta} = A_{k+1}^{\theta} \cdots A_{l}^{\theta}
\end{align}\label{l4.1.301}
for $l > k\geq 0$.

Iterating the second part of (\ref{4.1}), we get
\begin{align}\label{l4.1.1}
	W_{n}^{\theta}
	=
	W_{0}^{\theta} \tilde{A}_{0,n}^{\theta}
	+
	\sum_{k=1}^{n} B_{k}^{\theta} \tilde{A}_{k,n}^{\theta}
\end{align}
for $n\geq 1$.
Since $\tilde{A}_{0,n}^{\theta}\in{\cal P}^{N\times N}$,
we also have $e^{T}\tilde{A}_{0,n}^{\theta}=e^{T}$.
Consequently, the first part of (\ref{4.1}) implies
\begin{align*}
	V_{0}^{\theta} \tilde{A}_{0,n}^{\theta} \left(I - \frac{ee^{T} }{N} \right)
	=&
	W_{0}^{\theta} \tilde{A}_{0,n}^{\theta} \left(I - \frac{ee^{T} }{N} \right)
	-
	\frac{W_{0}^{\theta} e}{N}
	e^{T} \tilde{A}_{0,n}^{\theta} \left(I - \frac{ee^{T} }{N} \right)
	\\
	=&
	W_{0}^{\theta} \tilde{A}_{0,n}^{\theta} \left(I - \frac{ee^{T} }{N} \right).
\end{align*}
Combining this with the first part of (\ref{4.1}) and (\ref{l4.1.1}), we get
\begin{align}\label{l4.1.3}
	V_{n}^{\theta}
	=&
	W_{0}^{\theta} \tilde{A}_{0,n}^{\theta} \left(I - \frac{ee^{T} }{N} \right)
	+
	\sum_{k=1}^{n} B_{k}^{\theta} \tilde{A}_{k,n}^{\theta}
	\left(I - \frac{ee^{T} }{N} \right)
	\\ \nonumber
	=&
	V_{0}^{\theta} \tilde{A}_{0,n}^{\theta} \left(I - \frac{ee^{T} }{N} \right)
	+
	\sum_{k=1}^{n} B_{k}^{\theta} \tilde{A}_{k,n}^{\theta}
	\left(I - \frac{ee^{T} }{N} \right).
\end{align}

Owing to Assumptions \ref{a1}, \ref{a2}, we have
$\varepsilon^{2} \leq r_{\theta,\boldsymbol Y}^{k}\big(\hat{X}_{k,j}^{\theta} |\hat{X}_{k-1,i}^{\theta} \big)
\leq 1/\varepsilon^{2}$
and
\begin{align}\label{l4.1.501}
	\|\nabla_{\theta} r_{\theta,\boldsymbol Y}^{k}\big(\hat{X}_{k,j}^{\theta} |\hat{X}_{k-1,i}^{\theta} \big)
	\|
	\leq &
	q_{\theta}(Y_{k-1}|\hat{X}_{k-1,i} )
	\|\nabla_{\theta} p_{\theta}\big(\hat{X}_{k,j}^{\theta} |\hat{X}_{k-1,i}^{\theta} \big) \|
	\\ \nonumber
	&+
	p_{\theta}\big(\hat{X}_{k,j}^{\theta} |\hat{X}_{k-1,i}^{\theta} \big)
	\|\nabla_{\theta} q_{\theta}(Y_{k-1}|\hat{X}_{k-1,i} ) \|
	\leq
	\frac{2K}{\varepsilon}
\end{align}
for $1\leq i,j\leq N$, $k\geq 1$.
Therefore, we get
\begin{align*}
	&
	N\varepsilon^{2}
	\leq
	\sum_{i=1}^{N}
	r_{\theta,\boldsymbol Y}^{k}\big(\hat{X}_{k,j}^{\theta} |\hat{X}_{k-1,i}^{\theta} \big)
	\leq
	\frac{N}{\varepsilon^{2} },
	\;\;\;\;\;
	\sum_{i=1}^{N}
	\left\|
	\nabla_{\theta}
	r_{\theta,\boldsymbol Y}^{k}\big(\hat{X}_{k,j}^{\theta} |\hat{X}_{k-1,i}^{\theta} \big)
	\right\|
	\leq
	\frac{2KN}{\varepsilon}.
\end{align*}
Consequently, (\ref{4.5}) implies
\begin{align}\label{l4.1.5}
	A_{k,i,j}^{\theta}
	\geq
	\frac{\varepsilon^{4} }{N},
	\;\;\;\;\;
	\left\| B_{k,j}^{\theta} \right\|
	\leq
	\frac{2K}{\varepsilon^{3} },
	\;\;\;\;\;
	\left\| B_{k}^{\theta} \right\|
	=
	\max_{1\leq j\leq N} \|B_{k,j}^{\theta} \|
	\leq
	\frac{2K}{\varepsilon^{3} }.
\end{align}
Hence, (\ref{4.7}) yields
$\tau(A_{k}^{\theta} )\leq 1 - \varepsilon^{4} = \rho_{3}$.

Due to the well-known results in Markov chain theory (e.g. \cite[Theorems 15.2.4, 15.2.5]{bremaud}), we have
\begin{align}
	\tau(A'A'')\leq\tau(A')\tau(A''),
	\;\;\;\;\;
	\|A(z'-z'')\|_{1}\leq\tau(A)\|z'-z''\|_{1}
\end{align}
for any $A,A',A''\in{\cal P}^{N\times N}$, $z',z''\in{\cal P}^{N}$.
Then, using (\ref{l4.1.301}), we get
$\tau(\tilde{A}_{k,k}^{\theta} )=1$ and
\begin{align*}
	\tau(\tilde{A}_{k,l}^{\theta} )
	\leq
	\tau(A_{k+1}^{\theta} ) \cdots \tau(A_{l}^{\theta} )
	\leq
	\rho_{3}^{l-k}
\end{align*}
for $l>k\geq 0$.
Since $e_{i},\frac{e}{N}\in{\cal P}^{N}$, we deduce
\begin{align}\label{l4.1.23}
	\left\|
	\tilde{A}_{k,l}^{\theta} \left( e_{i} - \frac{e}{N} \right)
	\right\|_{1}
	\leq
	\tau(\tilde{A}_{k,l}^{\theta} )
	\left\|
	e_{i} - \frac{e}{N}
	\right\|_{1}
	\leq
	2\rho_{3}^{l-k}
\end{align}
for $1\leq i\leq N$, $l\geq k\geq 0$.
Consequently, (\ref{l4.1.5}) yields
\begin{align*}
	&
	\left\|
	V_{0}^{\theta} \tilde{A}_{0,l}^{\theta} \left( e_{i} - \frac{e}{N} \right)
	\right\|
	\leq
	\left\|V_{0}^{\theta}  \right\|
	\left\|
	\tilde{A}_{0,l}^{\theta} \left( e_{i} - \frac{e}{N} \right)
	\right\|_{1}
	\leq
	2\rho_{3}^{l} \left\|V_{0}^{\theta}  \right\|,
	\\
	&
	\left\|
	B_{k}^{\theta} \tilde{A}_{k,l}^{\theta} \left( e_{i} - \frac{e}{N} \right)
	\right\|
	\leq
	\left\|B_{k}^{\theta}  \right\|
	\left\|
	\tilde{A}_{k,l}^{\theta} \left( e_{i} - \frac{e}{N} \right)
	\right\|_{1}
	\leq
	\frac{4K\rho_{3}^{l-k} }{\varepsilon^{3} }.
\end{align*}
As vectors
$B_{k}^{\theta} \tilde{A}_{k,l}^{\theta} \left(e_{i} - \frac{e}{N} \right)$,
$V_{0}^{\theta} \tilde{A}_{0,l}^{\theta} \left(e_{i} - \frac{e}{N} \right)$
are the $i$-th columns of matrices
\linebreak
$B_{k}^{\theta} \tilde{A}_{k,l}^{\theta} \left(I - \frac{ee^{T} }{N} \right)$,
$V_{0}^{\theta} \tilde{A}_{0,l}^{\theta} \left(I - \frac{ee^{T} }{N} \right)$
(respectively), we conclude
\begin{align*}
	&
	\left\|
	V_{0}^{\theta} \tilde{A}_{0,l}^{\theta} \left( I - \frac{ee^{T} }{N} \right)
	\right\|
	=
	\max_{1\leq i\leq N}
	\left\|
	V_{0}^{\theta} \tilde{A}_{0,l}^{\theta} \left( e_{i} - \frac{e}{N} \right)
	\right\|
	\leq
	2\rho_{3}^{l} \left\|V_{0}^{\theta}  \right\|,
	\\
	&
	\left\|
	B_{k}^{\theta} \tilde{A}_{k,l}^{\theta} \left( I - \frac{ee^{T} }{N} \right)
	\right\|
	=
	\max_{1\leq i\leq N}
	\left\|
	B_{k}^{\theta} \tilde{A}_{k,l}^{\theta} \left( e_{i} - \frac{e}{N} \right)
	\right\|
	\leq
	\frac{4K\rho_{3}^{l-k} }{\varepsilon^{3} }.
\end{align*}
Hence, (\ref{l4.1.3}) implies 
\begin{align}\label{l4.1.25}
	\left\| V_{n}^{\theta} \right\|
	\leq &
	\left\|
	V_{0}^{\theta} \tilde{A}_{0,n}^{\theta} \left( I - \frac{ee^{T} }{N} \right)
	\right\|
	+
	\sum_{k=1}^{n}
	\left\|
	B_{k}^{\theta} \tilde{A}_{k,n}^{\theta} \left( I - \frac{ee^{T} }{N} \right)
	\right\|
	\\ \nonumber
	\leq &
	2\rho_{3}^{n} \left\| V_{0}^{\theta} \right\|
	+
	\frac{4K}{\varepsilon^{3} } \sum_{k=1}^{n} \rho_{3}^{n-k}
	\leq
	\frac{4K}{\varepsilon^{7} }
	\left(
	1
	+
	\rho_{3}^{n} \left\| V_{0}^{\theta} \right\|
	\right)
\end{align}
for $n\geq 1$. 
Since $W_{0,i}^{\theta} = w_{\theta}\big(\hat{X}_{0,i}^{\theta} \big)$,
Assumption \ref{a3} and (\ref{4.3}) yield 
\begin{align*}
	\left\|V_{0,i}^{\theta} \right\|
	\leq
	\left\|W_{0,i}^{\theta} \right\|
	+
	\frac{1}{N} \sum_{j=1}^{N} \left\|W_{0,j}^{\theta} \right\|
	\leq
	2\|w_{\theta} \|
\end{align*}
for $1\leq i\leq N$.
Thus, we have $\|V_{0}^{\theta}\|\leq\|w_{\theta}\|$.
Consequently, (\ref{l4.1.25}) implies 
\begin{align*}
	\left\|V_{n,i}^{\theta} \right\|
	\leq
	\left\|V_{n}^{\theta} \right\|
	\leq
	\frac{8K}{\varepsilon^{7} }
	\left(1 + \rho_{3}^{n} \|w_{\theta} \| \right)
\end{align*}
for $n\geq 0$. 
As $\hat{\zeta}_{n}^{\theta}(dx)=
\frac{1}{N}\sum_{j=1}^{N}V_{n,j}^{\theta}\delta_{\hat{X}_{n,j}^{\theta} }(dx)$
(due to (\ref{1.31}), (\ref{4.3})),
we get 
\begin{align*}
	\big\| \hat{\zeta}_{n}^{\theta}(B) \big\|
	\leq
	\frac{1}{N} \sum_{i=1}^{N} \big\| V_{n,i}^{\theta} \big\|
	\leq
	\frac{8K}{\varepsilon^{7} }
	\left(1 + \rho_{3}^{n} \|w_{\theta} \| \right)
\end{align*}
for $B\in{\cal B}({\cal X} )$. 
Hence, we have
\begin{align*}
	\big\| \hat{\zeta}_{n} \big\|
	\leq
	\frac{8dK}{\varepsilon^{7} }
	\left(1 + \rho_{3}^{n} \|w_{\theta} \| \right)
	=
	C_{3} \left(1 + \rho_{3}^{n} \|w_{\theta} \| \right).
\end{align*}
\end{proof}

\section{Proof of Main Results}\label{section5}

In this section, Proposition \ref{proposition5.1} is proved, while
Theorem \ref{theorem1.1} directly follows from it.
Lemma \ref{lemma5.2} and decompositions (\ref{l5.3.25*}), (\ref{p5.1.1}), (\ref{p5.1.29})
can be considered as the corner-stones in the proof of Proposition \ref{proposition5.1}
--- see inequalities
(\ref{p5.1.3}) -- (\ref{p5.1.25}), (\ref{p5.1.31}),
(\ref{p5.1.33}).
Proposition \ref{proposition2.1},
conditional distributions (\ref{l5.2.7}), (\ref{l5.2.9}) and
identities (\ref{l5.2.3}), (\ref{l5.2.53}), (\ref{l5.2.49})
are the main ingredients in the proof of Lemma \ref{lemma5.2}
--- see inequalities (\ref{l5.2.25}), (\ref{l5.2.27}), (\ref{l5.2.59}), (\ref{l5.2.61}),
(\ref{l5.2.67}).
Propositions \ref{proposition3.2}, \ref{proposition4.1} and Lemma \ref{lemma5.3} are important ingredients
of the proof of Lemma \ref{lemma5.2}, too
--- see inequalities (\ref{l5.2.21}), (\ref{l5.2.23}), (\ref{l5.2.57}),
(\ref{l5.2.201}), (\ref{l5.2.63}).
Proposition \ref{proposition3.2} plays an important role in the proof of Lemma \ref{lemma5.2}, either
--- see inequalities (\ref{l5.1.5}) -- (\ref{l5.1.9}).

Throughout this section, the following notation is used: $u_{\theta}(x)$, $\bar{w}_{\theta}$ and $\zeta_{\theta}(dx)$ are the functions
and the element of ${\cal M}_{s}^{d}({\cal X} )$ defined by
\begin{align}\label{5.1}
	\bar{w}_{\theta} = \int w_{\theta}(x') \xi_{\theta}(dx'),
	\;\;\;\;\;
	u_{\theta}(x) = w_{\theta}(x) - \bar{w}_{\theta},
	\;\;\;\;\;
	\zeta_{\theta}(B) = \int_{B} u_{\theta}(x')\xi_{\theta}(dx')
\end{align}
for $\theta\in\Theta$, $x\in{\cal X}$, $B\in{\cal B}({\cal X} )$.
$\hat{\xi}_{-1}^{\theta}(dx)$ and $\hat{\zeta}_{-1}^{\theta}(dx)$
are the elements of
${\cal P}({\cal X})$ and ${\cal M}_{s}^{d}({\cal X})$ (respectively)
defined by
\begin{align}\label{5.5}
	\hat{\xi}_{-1}^{\theta}(B) = \xi_{\theta}(B),
	\;\;\;\;\;
	\hat{\zeta}_{-1}^{\theta}(B) = \zeta_{\theta}(B).
\end{align}
$\hat{v}_{n}^{\theta}(x)$
is the (random) function defined by
\begin{align}\label{5.3}
	\hat{v}_{0}^{\theta}(x)
	\!=\!
	u_{\theta}(x),
	\;\;\;
	\hat{v}_{n}^{\theta}(x)
	\!=\!
	\frac{\int r_{\theta,\boldsymbol Y}^{n}(x|x') \hat{\zeta}_{n-1}^{\theta}(dx')
	+ \int \nabla_{\theta} r_{\theta,\boldsymbol Y}^{n}(x|x') \hat{\xi}_{n-1}^{\theta}(dx') }
	{\int r_{\theta,\boldsymbol Y}^{n}(x|x') \hat{\xi}_{n-1}^{\theta}(dx') }
\end{align}
for $n\geq 1$ ($r_{\theta,\boldsymbol Y}^{n}(x|x')$ and $\boldsymbol Y$ are defined in
(\ref{1.1}) and Subsection \ref{ssection1.3}, respectively).
$\hat{\alpha}_{n}^{\theta}(\xi)$ is the (random) function
mapping $\xi\in{\cal P}({\cal X})$ to ${\cal M}_{s}^{d}({\cal X})$ defined by
\begin{align}\label{5.9}
	\hat{\alpha}_{n}^{\theta}(\xi)(B)
	=
	\int_{B} \hat{v}_{n}^{\theta}(x) \xi(dx)
\end{align}
for $\xi\in{\cal P}({\cal X})$, $n\geq 0$.
$\hat{F}_{m:n}^{\theta}(dx)$,
$\hat{G}_{m:n}^{\theta}(dx)$ and
$\hat{H}_{m:n}^{\theta}(dx)$
are the (random) elements of ${\cal P}({\cal X})$, ${\cal M}_{s}^{d}({\cal X})$ and
${\cal M}_{s}^{d}({\cal X})$ (respectively) defined by
\begin{align}
	\label{5.7.1}
	\hat{F}_{-1:n}^{\theta}(B)
	=
	F_{\theta,\boldsymbol Y}^{0:n}(\hat{\xi}_{-1}^{\theta} )(B),
	&\;\;\;\;\;
	\hat{F}_{m:n}^{\theta}(B)
	=
	F_{\theta,\boldsymbol Y}^{m:n}(\hat{\xi}_{m}^{\theta} )(B),
	\\
	\label{5.7.3}
	\hat{G}_{-1:n}^{\theta}(B)
	=
	G_{\theta,\boldsymbol Y}^{0:n}(\hat{\xi}_{-1}^{\theta}, \hat{\zeta}_{-1}^{\theta} )(B),
	&	\;\;\;\;\;
	\hat{G}_{m:n}^{\theta}(B)
	=
	G_{\theta,\boldsymbol Y}^{m:n}(\hat{\xi}_{m}^{\theta}, \hat{\zeta}_{m}^{\theta} )(B),
	\\
	\label{5.7.5}
	\hat{H}_{-1:n}^{\theta}(B)
	=
	H_{\theta,\boldsymbol Y}^{0:n}(\hat{\xi}_{-1}^{\theta}, \hat{\zeta}_{-1}^{\theta} )(B),
	&	\;\;\;\;\;
	\hat{H}_{m:n}^{\theta}(B)
	=
	H_{\theta,\boldsymbol Y}^{m:n}(\hat{\xi}_{m}^{\theta}, \hat{\zeta}_{m}^{\theta} )(B)
\end{align}
for $n\geq m\geq 0$
($F_{\theta,\boldsymbol y}^{m:n}(\xi)$,
$G_{\theta,\boldsymbol y}^{m:n}(\xi,\zeta)$, $H_{\theta,\boldsymbol y}^{m:n}(\xi,\zeta)$
are specified in (\ref{1.21}) -- (\ref{1.25})).
$\hat{R}_{m:n}^{\theta}(x,dx')$ and
$\hat{S}_{m:n}^{\theta}(x,dx')$
are the (random) integral operators from ${\cal F}({\cal X})$ to
${\cal F}({\cal X})$, ${\cal F}^{d}({\cal X})$ (respectively)
defined by
\begin{align}
	&\label{5.21'}
	\hat{R}_{m:n}^{\theta}(x,B)=R_{\theta,\boldsymbol Y}^{m:n}(x,B),
	\;\;\;\;\;
	\hat{S}_{m:n}^{\theta}(x,B)=S_{\theta,\boldsymbol Y}^{m:n}(x,B).
\end{align}
$\hat{\Psi}_{\theta,\boldsymbol Y}^{m:n}(x,dx')$ and
$\hat{\Phi}_{\theta,\boldsymbol Y}^{m:n}(x,dx')$
are the (random) integral operators from ${\cal F}({\cal X})$ to ${\cal F}^{d}({\cal X})$
defined by
\begin{align}
	&\label{5.21}
	\hat{\Psi}_{m:n}^{\theta}(x,B)
	=
	\hat{R}_{m:n}^{\theta}(x,B) \hat{v}_{m}^{\theta}(x)
	+
	\hat{S}_{m:n}^{\theta}(x,B),
	\\
	&\label{5.23}
	\hat{\Phi}_{m:n}^{\theta}(x,B)
	=
	\hat{\Psi}_{m:n}^{\theta}(x,B)
	-
	\hat{F}_{m-1:n}^{\theta}(B)
	\hat{\Psi}_{m:n}^{\theta}(x,{\cal X}).
\end{align}
$\hat{C}_{m:n}^{\theta}(dx)$,
$\hat{B}_{m:n}^{\theta}(dx)$ and
$\hat{A}_{m:n}^{\theta}(dx)$
are the (random) elements of ${\cal M}_{s}^{d}({\cal X})$ defined by
\begin{align}
	&\label{5.25}
	\hat{C}_{m:n}^{\theta}(B)
	=
	\frac{(\hat{\xi}_{m}^{\theta}\hat{\Psi}_{m:n}^{\theta} )(B) }
	{\big\langle\hat{\xi}_{m}^{\theta}\hat{R}_{m:n}^{\theta}\big\rangle },
	\\
	&\label{5.27}
	\hat{B}_{m:n}^{\theta}(B)
	=
	-
	\big(
	\hat{F}_{m:n}^{\theta}(B)
	-
	\hat{F}_{m-1:n}^{\theta}(B)
	\big)
	\big\langle\hat{C}_{m:n}^{\theta}\big\rangle,
	\\
	&\label{5.29}
	\hat{A}_{m:n}^{\theta}(B)
	=
	\hat{C}_{m:n}^{\theta}(B)
	-
	\hat{F}_{m-1:n}^{\theta}(B)
	\big\langle\hat{C}_{m:n}^{\theta}\big\rangle.
\end{align}

\begin{lemma}\label{lemma5.3}
Let $\theta$, $B$, $\xi$ be any elements of 
$\Theta$, ${\cal B}({\cal X})$, ${\cal P}({\cal X})$ (respectively). 
Moreover, let $n$, $m$ be any integers satisfying $n\geq m\geq 0$. 

(i) Suppose that Assumption \ref{a1} holds. Then, we have
\begin{align}\label{l5.3.1*}
	\hat{F}_{m-1:n}^{\theta}(B)
	=
	\frac{\big( \hat{F}_{m-1:m}^{\theta} \hat{R}_{m:n}^{\theta} \big)(B) }
	{\big\langle \hat{F}_{m-1:m}^{\theta} \hat{R}_{m:n}^{\theta} \big\rangle }. 
\end{align}

(ii) Suppose that Assumptions \ref{a1} -- \ref{a3} hold. Then, we have
\begin{align}
	&\label{l5.3.25*}
	\hat{G}_{m:n}^{\theta}(B)
	=
	\hat{A}_{m:n}^{\theta}(B)
	+
	\hat{B}_{m:n}^{\theta}(B),
	\\
	&\label{l5.3.21*}
	\frac{\big( \xi\hat{\Psi}_{m:n}^{\theta} \big)(B) }
	{\big\langle \xi\hat{R}_{m:n}^{\theta} \big\rangle }
	=
	H_{\theta,\boldsymbol Y}^{m:n}
	\big(\xi,\hat{\alpha}_{m}^{\theta}(\xi) \big)(B),
	\\
	&\label{l5.3.23*}
	\begin{aligned}[t]
	\frac{\big( \xi\hat{\Phi}_{m:n}^{\theta} \big)(B) }
	{\big\langle \xi\hat{R}_{m:n}^{\theta} \big\rangle }
	=&
	\big(
	F_{\theta,\boldsymbol Y}^{m:n}(\xi)(B)
	-
	\hat{F}_{m-1:n}^{\theta}(B)
	\big)
	\big\langle
	H_{\theta,\boldsymbol Y}^{m:n}\big(\xi, \hat{\alpha}_{m}^{\theta}(\xi) \big)
	\big\rangle
	\\
	&+
	G_{\theta,\boldsymbol Y}^{m:n}
	\big(\xi, \hat{\alpha}_{m}^{\theta}(\xi) \big)(B). 
	\end{aligned}
\end{align}
We also have
\begin{align}
	\label{l5.3.3*}
	\hat{H}_{m-1:n}^{\theta}(B)
	\!=\!
	\frac{\big( \hat{F}_{m-1:m}^{\theta} \hat{\Psi}_{m:n}^{\theta} \big)(B) }
	{\big\langle \hat{F}_{m-1:m}^{\theta} \hat{R}_{m:n}^{\theta} \big\rangle },
	\;\:
	\hat{G}_{m-1:n}^{\theta}(B)
	\!=\!
	\frac{\big( \hat{F}_{m-1:m}^{\theta} \hat{\Phi}_{m:n}^{\theta} \big)(B) }
	{\big\langle \hat{F}_{m-1:m}^{\theta} \hat{R}_{m:n}^{\theta} \big\rangle }. 
\end{align}
\end{lemma}

Lemma \ref{lemma5.3} summarizes relatively straightforward relationships between
the measures defined in (\ref{5.7.1}) -- (\ref{5.7.5}) and (\ref{5.25}) -- (\ref{5.29}).
A detailed proof of the lemma is provided in the supplementary material
(Section SM3).

\begin{lemma}\label{lemma5.1}
Let $\theta$ be any element of $\Theta$, 
while $\xi$, $\xi'$ are any elements of ${\cal P}({\cal X})$.
Moreover, let $n$, $m$ be any integers safisfying $n\geq m\geq 0$.  
Suppose that Assumptions \ref{a1} and \ref{a2} hold.
Then, there exist real numbers $\rho_{4}\in(0,1)$, $C_{4}\in[1,\infty)$
(independent of $\theta$, $\xi$, $\xi'$, $n$, $m$
and depending only on $\varepsilon$, $d$, $K$) such that
\begin{align}
	&\label{l5.1.1*}
	\max\left\{
	\big\|
	\hat{C}_{m:n}^{\theta}\big\|,
	\big\| \hat{H}_{m-1:n}^{\theta} \big\|
	\right\}
	\leq
	C_{4} \left(1 + n - m + \rho_{4}^{m} \|w_{\theta} \| \right),
	\\
	&\label{l5.1.3*}
	\left\|
	\frac{\xi\hat{\Psi}_{m:n}^{\theta} }
	{\big\langle \xi\hat{R}_{m:n}^{\theta} \big\rangle }
	-
	\frac{\xi'\hat{\Psi}_{m:n}^{\theta} }
	{\big\langle \xi'\hat{R}_{m:n}^{\theta} \big\rangle }
	\right\|
	\leq
	C_{4} \left(1 + n - m + \rho_{4}^{m} \|w_{\theta} \| \right),
	\\
	&\label{l5.1.5*}
	\left\|
	\frac{\xi\hat{\Phi}_{m:n}^{\theta} }
	{\big\langle \xi\hat{R}_{m:n}^{\theta} \big\rangle }
	-
	\frac{\xi'\hat{\Phi}_{m:n}^{\theta} }
	{\big\langle \xi'\hat{R}_{m:n}^{\theta} \big\rangle }
	\right\|
	\leq
	C_{4} \rho_{4}^{n-m} \left(1 + n - m + \rho_{4}^{m} \|w_{\theta} \| \right).
\end{align}
\end{lemma}

\begin{proof}
Throughout the proof, the following notation is used.
$x$, $x'$ are any elements of ${\cal X}$, while $B$ is any element of ${\cal B}({\cal X} )$. 
$\rho_{4}$ is the real number defined by
$\rho_{4}=\max\big\{\rho_{1}, \rho_{2}, \rho_{3} \big\}$,
while $\tilde{C}_{1}$, $\tilde{C}_{2}$, $\tilde{C}_{3}$, $\tilde{C}_{4}$, $C_{4}$
are the real numbers defined as
$\tilde{C}_{1}=2dC_{3}$,
$\tilde{C}_{2}=4\tilde{C}_{1}K\varepsilon^{-4}\rho_{4}^{-1}$,
$\tilde{C}_{3}=d\tilde{C}_{2}$, $\tilde{C}_{4}=3C_{2}\tilde{C}_{1}\tilde{C}_{3}$,
$C_{4}=3C_{1}\tilde{C}_{4}$
($\varepsilon$, $\rho_{1}$, $\rho_{2}$, $K$, $C_{1}$, $C_{2}$ are specified in Assumptions \ref{a1}, \ref{a2}
and Proposition \ref{proposition3.2}).
$n$, $m$ are any integers satisfying $n\geq m\geq 0$.

Relying on Assumption \ref{a3} and (\ref{5.1}), (\ref{5.3}), we conclude
\begin{align}\label{l5.1.301}
	\|\hat{v}_{0}^{\theta}(x) \|
	=
	\|u_{\theta}(x) \|
	\leq
	\|w_{\theta}(x) \|
	+
	\int \|w_{\theta}(x') \|\xi_{\theta}(dx')
	\leq
	2\|w_{\theta} \|.
\end{align}
Consequently, (\ref{5.1}), (\ref{5.5}) imply
\begin{align*}
	\|\hat{\zeta}_{-1}^{\theta}(B) \|
	=
	\|\zeta_{\theta}(B) \|
	\leq
	\int_{B} \|u_{\theta}(x) \|\xi_{\theta}(dx)
	\leq
	2\|w_{\theta}\|.
\end{align*}
Hence, we have
$\|\hat{\zeta}_{-1}^{\theta}\|=\|\zeta_{\theta}\|\leq 2d\|w_{\theta}\|$.
Combining this with Proposition \ref{proposition4.1}, we get
\begin{align}\label{l5.1.303}
	\|\hat{\zeta}_{k}^{\theta} \|
	\leq
	2dC_{3}(1 + \rho_{3}^{k}\|w_{\theta}\| )
	\leq
	\tilde{C}_{1}(1 + \rho_{4}^{k}\|w_{\theta}\| )
\end{align}
for $k\geq -1$.

Using Assumptions \ref{a1}, \ref{a2} and the same arguments as in Proposition \ref{proposition4.1}
(see (\ref{l4.1.501})), we deduce
\begin{align*}
	\varepsilon^{2}
	\leq
	r_{\theta,\boldsymbol Y}^{k}(x'|x)
	\leq
	\frac{1}{\varepsilon^{2} },
	\;\;\;\;\;
	\big\|
	\nabla_{\theta} r_{\theta,\boldsymbol Y}^{k}(x'|x)
	\big\|
	\leq
	\frac{2K}{\varepsilon}
\end{align*}
for $k\geq 1$.
Then, (\ref{5.3}), (\ref{l5.1.303}) yield
\begin{align*}
	\left\|\hat{v}_{k}^{\theta}(x) \right\|
	\leq &
	\frac{
	\int r_{\theta,\boldsymbol Y}^{k}(x|x') \: |\hat{\zeta}_{k-1}^{\theta} |(dx')
	+
	\int \|\nabla_{\theta} r_{\theta,\boldsymbol Y}^{k}(x|x') \| \: \hat{\xi}_{k-1}^{\theta}(dx') }
	{\int r_{\theta,\boldsymbol Y}^{k}(x|x') \: \hat{\xi}_{k-1}^{\theta}(dx') }
	\\
	\leq &
	\frac{2K}{\varepsilon^{3} }
	+
	\frac{\|\hat{\zeta}_{k-1}^{\theta}\|}{\varepsilon^{4} }
	\leq
	\frac{4\tilde{C}_{1}K (1+\rho_{3}^{k-1}\|w_{\theta} \| ) }{\varepsilon^{4} }
	\leq
	\tilde{C}_{2} (1+\rho_{4}^{k}\|w_{\theta} \|  ).
\end{align*}
Combining this with (\ref{5.9}), (\ref{l5.1.301}), we get
\begin{align}\label{l5.1.1}
	\left\|\hat{\alpha}_{k}^{\theta}(\xi)(B) \right\|
	\leq
	\int_{B} \left\|\hat{v}_{k}^{\theta}(x) \right\| \xi(dx)
	\leq
	\tilde{C}_{2} (1+\rho_{4}^{k}\|w_{\theta} \|  )
\end{align}
for $k\geq 0$.
Thus, we have
\begin{align}\label{l5.1.3}
	\left\|\hat{\alpha}_{k}^{\theta}(\xi) \right\|
	\leq
	d\tilde{C}_{2}(1 + \rho_{4}^{k} \|w_{\theta} \| )
	=
	\tilde{C}_{3}(1 + \rho_{4}^{k} \|w_{\theta} \| ).
\end{align}
Consequently, Proposition \ref{proposition3.2} implies
\begin{align}\label{l5.1.5}
	\left\| H_{\theta,\boldsymbol Y}^{m:n}(\xi,\hat{\alpha}_{m}^{\theta}(\xi) ) \right\|
	\leq &
	C_{2}
	\left(
	n-m
	+
	\left\|\hat{\alpha}_{m}^{\theta}(\xi) \right\|
	\right)
	\\ \nonumber
	\leq &
	C_{2}\tilde{C}_{3}
	\left(
	1 + n - m + \rho_{4}^{m}\|w_{\theta}\|
	\right)
	\\ \nonumber
	\leq&
	\tilde{C}_{4}
	\left(
	1 + n - m + \rho_{4}^{m}\|w_{\theta}\|
	\right).
\end{align}
Then, relying on Lemma \ref{lemma5.3} and (\ref{5.25}), we deduce
\begin{align}
	\label{l5.1.5''}
	\left\| \hat{C}_{m:n}^{\theta} \right\|
	=
	\left\| H_{\theta,\boldsymbol Y}^{m:n}\big(\hat{\xi}_{m}^{\theta},
	\hat{\alpha}_{m}^{\theta}(\hat{\xi}_{m}^{\theta}) \big)\right\|
	\leq
	\tilde{C}_{4}
	\left(
	1 + n - m + \rho_{4}^{m}\|w_{\theta}\|
	\right).
\end{align}
Moreover, if $m\geq 1$, Proposition \ref{proposition3.2}
and (\ref{5.7.5}), (\ref{l5.1.303}), (\ref{l5.1.5}) yield
\begin{align}\label{l5.1.5'}
	\left\| \hat{H}_{m-1:n}^{\theta} \right\|
	=
	\left\| H_{\theta,\boldsymbol Y}^{m-1:n}(\hat{\xi}_{m-1}^{\theta}, \hat{\zeta}_{m-1}^{\theta} ) \right\|
	\leq &
	C_{2} \big(n - m + \|\hat{\zeta}_{m-1}^{\theta} \| \big)
	\\ \nonumber
	\leq &
	C_{2}\tilde{C}_{1}(1 + n - m + \rho_{4}^{m-1} \|w_{\theta} \| )
	\\ \nonumber
	\leq &
	\tilde{C}_{4}(1 + n - m + \rho_{4}^{m} \|w_{\theta} \| ).
\end{align}
The same arguments also imply
\begin{align}\label{l5.1.305}
	\left\| \hat{H}_{-1:n}^{\theta} \right\|
	=
	\left\| H_{\theta,\boldsymbol Y}^{0:n}(\hat{\xi}_{-1}^{\theta}, \hat{\zeta}_{-1}^{\theta} ) \right\|
	\leq
	C_{2} \big(n + \|\hat{\zeta}_{-1}^{\theta} \| \big)
	\leq &
	C_{2}\tilde{C}_{1}(1 + n + \rho_{4}^{-1} \|w_{\theta} \| )
	\\ \nonumber
	\leq &
	\tilde{C}_{4}(1 + n + \rho_{4}^{-1} \|w_{\theta} \| ).
\end{align}
Using (\ref{l5.1.5}) -- (\ref{l5.1.305}), we conclude that (\ref{l5.1.1*}) holds.

Owing to Proposition \ref{proposition3.2}, Lemma \ref{lemma5.3}
and (\ref{5.7.1}), we have
\begin{align}\label{l5.1.7}
	\left\|
	F_{\theta,\boldsymbol Y}^{m:n}(\xi)
	-
	\hat{F}_{m-1:n}^{\theta}
	\right\|
	=
	\left\|
	F_{\theta,\boldsymbol Y}^{m:n}(\xi)
	-
	F_{\theta,\boldsymbol Y}^{m:n}
	\big( \hat{F}_{m-1:m}^{\theta} \big)
	\right\|
	\leq
	C_{1} \rho_{1}^{n-m}.
\end{align}
Due to the same proposition and (\ref{l5.1.3}), we also have
\begin{align}\label{l5.1.9}
	\left\|
	G_{\theta,\boldsymbol Y}^{m:n}(\xi,\hat{\alpha}_{m}^{\theta}(\xi) )
	-
	G_{\theta,\boldsymbol Y}^{m:n}(\xi',\hat{\alpha}_{m}^{\theta}(\xi') )
	\right\|
	\leq &
	C_{2}\rho_{2}^{n-m}
	\!\left(
	1
	+\!
	\left\|\hat{\alpha}_{m}^{\theta}(\xi) \right\|
	+\!
	\left\|\hat{\alpha}_{m}^{\theta}(\xi') \right\|
	\right)
	\\ \nonumber
	\leq &
	3C_{2}\tilde{C}_{3}\rho_{2}^{n-m}
	\left( 1 + \rho_{4}^{m} \|w_{\theta}\| \right)
	\\ \nonumber
	\leq &
	\tilde{C}_{4}\rho_{2}^{n-m}
	\left( 1 + \rho_{4}^{m} \|w_{\theta}\| \right).
\end{align}
Combining Lemma \ref{lemma5.3} and (\ref{l5.1.5}), (\ref{l5.1.7}),  (\ref{l5.1.9}), we get
\begin{align}\label{l5.1.21}
	\left\|
	\frac{\xi\hat{\Phi}_{m:n}^{\theta} }
	{\big\langle \xi\hat{R}_{m:n}^{\theta} \big\rangle }
	-
	\frac{\xi'\hat{\Phi}_{m:n}^{\theta} }
	{\big\langle \xi'\hat{R}_{m:n}^{\theta} \big\rangle }
	\right\|
	\leq &
	\left\|
	G_{\theta,\boldsymbol Y}^{m:n}(\xi,\hat{\alpha}_{m}^{\theta}(\xi) )
	-
	G_{\theta,\boldsymbol Y}^{m:n}(\xi',\hat{\alpha}_{m}^{\theta}(\xi') )
	\right\|
	\\ \nonumber
	&+
	\left\|
	F_{\theta,\boldsymbol Y}^{m:n}(\xi)
	-
	\hat{F}_{m-1:n}^{\theta}
	\right\|
	\left\| H_{\theta,\boldsymbol Y}^{m:n}(\xi,\hat{\alpha}_{m}^{\theta}(\xi) ) \right\|
	\\ \nonumber
	&+
	\left\|
	F_{\theta,\boldsymbol Y}^{m:n}(\xi')
	-
	\hat{F}_{m-1:n}^{\theta}
	\right\|
	\left\| H_{\theta,\boldsymbol Y}^{m:n}(\xi',\hat{\alpha}_{m}^{\theta}(\xi') ) \right\|
	\\ \nonumber
	\leq &
	\tilde{C}_{4}\rho_{2}^{n-m}
	\!\left( 1 + \!\rho_{4}^{m} \|w_{\theta}\| \right)
	\\
	&+
	2C_{1}\tilde{C}_{4}\rho_{1}^{n-m}
	\!\left(
	1 + n - m + \!\rho_{4}^{m} \|w_{\theta}\|
	\right)
	\\ \nonumber
	\leq &
	3C_{1}\tilde{C}_{4}\rho_{4}^{n-m}
	\left( 1 + n - m + \!\rho_{4}^{m} \|w_{\theta}\| \right).
\end{align}
Similarly, relying on Lemma \ref{lemma5.3} and (\ref{l5.1.5}), we get
\begin{align}\label{l5.1.21'}
	\left\|
	\frac{\xi\hat{\Psi}_{m:n}^{\theta} }
	{\big\langle \xi\hat{R}_{m:n}^{\theta} \big\rangle }
	-
	\frac{\xi'\hat{\Psi}_{m:n}^{\theta} }
	{\big\langle \xi'\hat{R}_{m:n}^{\theta} \big\rangle }
	\right\|
	\leq &
	\left\| H_{\theta,\boldsymbol Y}^{m:n}(\xi,\hat{\alpha}_{m}^{\theta}(\xi) ) \right\|
	+
	\left\| H_{\theta,\boldsymbol Y}^{m:n}(\xi',\hat{\alpha}_{m}^{\theta}(\xi') ) \right\|
	\\ \nonumber
	\leq &
	2\tilde{C}_{4}
	\!\left(
	1 + n - m + \!\rho_{4}^{m} \|w_{\theta}\|
	\right).
\end{align}
Using (\ref{l5.1.21}), (\ref{l5.1.21'}),
we deduce that (\ref{l5.1.3*}) holds.
\end{proof}

\begin{lemma}\label{lemma5.2}
Let $\theta$ be any element of $\Theta$, 
while $\varphi:{\cal X}\rightarrow[-1,1]$ is any Borel-measurable function. 
Moreover, let $n$, $m$ be any integers satisfying $n\geq m\geq 0$.  

(i) Suppose that Assumption \ref{a1} holds.
Then, there exists a real number $C_{5}\in[1,\infty)$
(independent of $N$, $\theta$, $\varphi(x)$, $n$, $m$ 
and depending only on $\varepsilon$) such that
\begin{align}
	&\label{l5.2.21*}
	\left|
	E\left(\left.
	\hat{F}_{m:n}^{\theta}(\varphi)
	-
	\hat{F}_{m-1:n}^{\theta}(\varphi)
	\right|
	\boldsymbol Y
	\right)
	\right|
	\leq
	\frac{C_{5}\rho_{1}^{n-m} }{N},
	\\
	&\label{l5.2.23*}
	\left(
	E\left(\left.
	\left|
	\hat{F}_{m:n}^{\theta}(\varphi)
	-
	\hat{F}_{m-1:n}^{\theta}(\varphi)
	\right|^{2}
	\right|
	\boldsymbol Y
	\right)
	\right)^{1/2}
	\leq
	\frac{C_{5}\rho_{1}^{n-m} }{\sqrt{N} }
\end{align}
almost surely
($\rho_{1}$ is specified in Proposition \ref{proposition3.2}).

(ii) Suppose that Assumptions \ref{a1} -- \ref{a3} hold.
Then, there exist real numbers $\rho_{5}\in(0,1)$, $C_{6}\in[1,\infty)$
(independent of $N$, $\theta$, $\varphi(x)$, $n$, $m$ 
and depending only on $\varepsilon$, $d$, $K$) such that
\begin{align}
	&\label{l5.2.1*}
	\left\|
	E\left(\left.
	\hat{A}_{m:n}^{\theta}(\varphi)
	-
	\hat{G}_{m-1:n}^{\theta}(\varphi)
	\right|
	\boldsymbol Y
	\right)
	\right\|
	\leq
	\frac{C_{6}(\rho_{5}^{n-m} + \rho_{5}^{n} \|w_{\theta}\| ) }{N},
	\\
	&\label{l5.2.3*}
	\left(
	E\left(\left.
	\left\|
	\hat{A}_{m:n}^{\theta}(\varphi)
	-
	\hat{G}_{m-1:n}^{\theta}(\varphi)
	\right\|^{2}
	\right|
	\boldsymbol Y
	\right)
	\right)^{1/2}
	\leq
	\frac{C_{6}(\rho_{5}^{n-m} + \rho_{5}^{n} \|w_{\theta}\| ) }{\sqrt{N} },
	\\
	&\label{l5.2.5*}
	\left\|
	E\left(\left.
	\hat{B}_{m:n}^{\theta}(\varphi)
	\right|
	\boldsymbol Y
	\right)
	\right\|
	\leq
	\frac{C_{6}(\rho_{5}^{n-m} + \rho_{5}^{n} \|w_{\theta}\| ) }{N},
	\\
	&\label{l5.2.7*}
	\left(
	E\left(\left.
	\left\|
	\hat{B}_{m:n}^{\theta}(\varphi)
	\right\|^{2}
	\right|
	\boldsymbol Y
	\right)
	\right)^{1/2}
	\leq
	\frac{C_{6}(\rho_{5}^{n-m} + \rho_{5}^{n} \|w_{\theta}\| ) }{\sqrt{N} }
\end{align}
almost surely.
\end{lemma}

\begin{proof}
(i) Throughout this part of the proof, the following notation is used. 
$\xi$, $\xi'$ are any elements of ${\cal P}({\cal X} )$.
$\boldsymbol 1(x)$ is the function which maps $x\in{\cal X}$ to one.

Relying on (\ref{1.1}), (\ref{1.3}), (\ref{1.35}),
(\ref{1.7}), (\ref{1.21}), (\ref{5.7.1}), (\ref{5.21'}), we conclude
\begin{align}\label{l5.2.7}
	P\!\left(\!\left.
	\hat{X}_{k,1}^{\theta}\in B_{1}, \dots, \hat{X}_{k,N}^{\theta}\in B_{N}
	\right| \!\boldsymbol Y, \hat{\xi}_{k-1}^{\theta}
	\!\right)
	\!=&
	\!\prod_{i=1}^{N}
	\!\left(\!
	\frac{\int\!\left(\!\int_{B_{i} }(x') r_{\theta,\boldsymbol Y}^{k}(x'|x)
	\mu(dx') \right)
	\hat{\xi}_{k-1}^{\theta}(dx) }
	{\int\!\left(\!\int r_{\theta,\boldsymbol Y}^{k}(x'|x)
	\mu(dx') \right)
	\hat{\xi}_{k-1}^{\theta}(dx) }
	\right)
	\\ \nonumber
	\!=&
	\prod_{i=1}^{N}
	\frac{\big(\hat{\xi}_{k-1}\hat{R}_{k-1:k}^{\theta}\big)(B_{i} )}
	{\big\langle\hat{\xi}_{k-1}\hat{R}_{k-1:k}^{\theta}\big\rangle }
	=
	\prod_{i=1}^{N}
	\hat{F}_{n-1:n}^{\theta} (B_{i} )
\end{align}
almost surely for any $B_{1},\dots,B_{N}\in{\cal B}({\cal X})$, $k\geq 1$.
Similarly, using (\ref{5.5}), (\ref{5.7.1}), we deduce
\begin{align}\label{l5.2.9}
	P\left(\left.
	\hat{X}_{0,1}^{\theta}\in B_{1}, \dots, \hat{X}_{0,N}^{\theta}\in B_{N}
	\right| \boldsymbol Y, \hat{\xi}_{-1}^{\theta}
	\right)
	=
	\prod_{i=1}^{N}
	\hat{\xi}_{-1}^{\theta}(B_{i} )
	=
	\prod_{i=1}^{N}
	\hat{F}_{-1:0}^{\theta} (B_{i} )
\end{align}
almost surely for the same $B_{1},\dots,B_{N}$.
Moreover, Lemma \ref{lemma5.3} and (\ref{1.21}), (\ref{5.7.1}), (\ref{5.21'}) imply
\begin{align}
	&\label{l5.2.3}
	\hat{F}_{m-1:n}^{\theta}(\varphi)
	=
	\frac{\big( \hat{F}_{m-1:m}^{\theta}\hat{R}_{m:n}^{\theta} \big)(\varphi) }
	{\big( \hat{F}_{m-1:m}^{\theta}\hat{R}_{m:n}^{\theta} \big)(\boldsymbol 1) },
	\;\;\;\;\;
	\hat{F}_{m:n}^{\theta}(\varphi)
	=
	\frac{\big( \hat{\xi}_{m}^{\theta}\hat{R}_{m:n}^{\theta} \big)(\varphi) }
	{\big( \hat{\xi}_{m}^{\theta}\hat{R}_{m:n}^{\theta} \big)(\boldsymbol 1) }.
\end{align}

Let $C_{5}=2C_{1}^{3}$
($C_{1}$ is specified in Proposition \ref{proposition3.2}).
Owing to Proposition \ref{proposition3.2}
and (\ref{1.21}), (\ref{5.7.1}), (\ref{5.21'}), we have
\begin{align}\label{l5.2.21}
	\left|
	\frac{\big( \xi\hat{R}_{m:n}^{\theta} \big)(\varphi) }
	{\big( \xi\hat{R}_{m:n}^{\theta} \big)(\boldsymbol 1) }
	\!-\!
	\frac{\big( \xi'\hat{R}_{m:n}^{\theta} \big)(\varphi) }
	{\big( \xi'\hat{R}_{m:n}^{\theta} \big)(\boldsymbol 1) }
	\right|
	\!=\!
	\left|
	F_{\theta,\boldsymbol Y}^{m:n}(\xi)(\varphi)
	\!-\!
	F_{\theta,\boldsymbol Y}^{m:n}(\xi')(\varphi)
	\right|
	\!\leq\!
	C_{1} \rho_{1}^{n\!-\!m}.
\end{align}
Due to the same arguments, we also have
\begin{align}\label{l5.2.23}
	\frac{\big( \xi\hat{R}_{m:n}^{\theta} \big)(\boldsymbol 1) }
	{\big( \xi'\hat{R}_{m:n}^{\theta} \big)(\boldsymbol 1) }
	\leq
	C_{1}.
\end{align}
Using Proposition \ref{proposition2.1}
and (\ref{l5.2.7}) -- (\ref{l5.2.23}), we conclude
\begin{align}\label{l5.2.25}
	&
	\left|
	E\left(\left.
	\hat{F}_{m:n}^{\theta}(\varphi)
	-
	\hat{F}_{m-1:n}^{\theta}(\varphi)
	\right|
	\boldsymbol Y, \hat{\xi}_{m-1}^{\theta}
	\right)
	\right|
	\\ \nonumber
	&=
	\left|
	E\left(\left.
	\frac{\big( \hat{\xi}_{m}^{\theta}\hat{R}_{m:n}^{\theta} \big)(\varphi) }
	{\big( \hat{\xi}_{m}^{\theta}\hat{R}_{m:n}^{\theta} \big)(\boldsymbol 1) }
	-
	\frac{\big( \hat{F}_{m-1:m}^{\theta}\hat{R}_{m:n}^{\theta} \big)(\varphi) }
	{\big( \hat{F}_{m-1:m}^{\theta}\hat{R}_{m:n}^{\theta} \big)(\boldsymbol 1) }
	\right|
	\boldsymbol Y, \hat{\xi}_{m-1}^{\theta}
	\right)
	\right|
	\\ \nonumber
	&
	\leq
	\frac{2C_{1}^{3}\rho_{1}^{n-m} }{N}
	=
	\frac{C_{5}\rho_{1}^{n-m} }{N}
\end{align}
almost surely.\footnote
{To get (\ref{l5.2.25}), (\ref{l5.2.27}), the following should be done:
In Proposition \ref{proposition2.1}, set $z=x$, $k=N$ and replace
$f(z)$, $g(z)$, $\xi_{k}(dz)$, $\xi(dz)$ with
$\big(\hat{R}_{m:n}^{\theta}\varphi\big)(x)$, $\big(\hat{R}_{m:n}^{\theta}\boldsymbol 1\big)(x)$,
$\hat{\xi}_{m}^{\theta}(dx)$, $\hat{F}_{m-1:m}^{\theta}(dx)$. }
Relying on the same arguments, we deduce
\begin{align}\label{l5.2.27}
	&
	E\left(\left.
	\left|
	\hat{F}_{m:n}^{\theta}(\varphi)
	-
	\hat{F}_{m-1:n}^{\theta}(\varphi)
	\right|^{2}
	\right|
	\boldsymbol Y, \hat{\xi}_{m-1}^{\theta}
	\!\right)
	\\ \nonumber
	&=
	E\left(\left.
	\left|
	\frac{\big( \hat{\xi}_{m}^{\theta}\hat{R}_{m:n}^{\theta} \big)(\varphi) }
	{\big( \hat{\xi}_{m}^{\theta}\hat{R}_{m:n}^{\theta} \big)(\boldsymbol 1) }
	-
	\frac{\big( \hat{F}_{m-1:m}^{\theta}\hat{R}_{m:n}^{\theta} \big)(\varphi) }
	{\big( \hat{F}_{m-1:m}^{\theta}\hat{R}_{m:n}^{\theta} \big)(\boldsymbol 1) }
	\right|^{2}
	\right|
	\boldsymbol Y, \hat{\xi}_{m-1}^{\theta}
	\!\right)
	\\ \nonumber
	&
	\leq
	\left( \frac{2C_{1}^{2}\rho_{1}^{n-m} }{\sqrt{N} } \right)^{2}
	\leq
	\left( \frac{C_{5}\rho_{1}^{n-m} }{\sqrt{N} } \right)^{2}
\end{align}
almost surely.
Combining (\ref{l5.2.25}), (\ref{l5.2.27}) with the tower property of conditional expectations,
we conclude that (\ref{l5.2.21*}), (\ref{l5.2.23*}) hold almost surely.

(ii) Let $\xi$, $\xi'$, $\boldsymbol 1(x)$ have the same meaning as in (i).
Using (\ref{5.23}), (\ref{5.25}), (\ref{5.29}), it is straightforward to verify
\begin{align}
	&\label{l5.2.53}
	\hat{A}_{m:n}^{\theta}(\varphi)
	=
	\frac{\big( \hat{\xi}_{m}^{\theta}\hat{\Phi}_{m:n}^{\theta} \big)(\varphi) }
	{\big( \hat{\xi}_{m}^{\theta}\hat{R}_{m:n}^{\theta} \big)(\boldsymbol 1) },
	\;\;\;\;\;
	\hat{C}_{m:n}^{\theta}(\varphi)
	=
	\frac{\big( \hat{\xi}_{m}^{\theta}\hat{\Psi}_{m:n}^{\theta} \big)(\varphi) }
	{\big( \hat{\xi}_{m}^{\theta}\hat{R}_{m:n}^{\theta} \big)(\boldsymbol 1) }.
\end{align}
Similarly, Lemma \ref{lemma5.3} yields
\begin{align}\label{l5.2.49}
	&
	\hat{G}_{m-1:n}^{\theta}(\varphi)
	\!=\!
	\frac{\big( \hat{F}_{m-1:m}^{\theta}\hat{\Phi}_{m:n}^{\theta} \big)(\varphi) }
	{\big( \hat{F}_{m-1:m}^{\theta}\hat{R}_{m:n}^{\theta} \big)(\boldsymbol 1) },
	\;\;\;
	\hat{H}_{m-1:n}^{\theta}
	(\varphi)
	\!=\!
	\frac{\big( \hat{F}_{m-1:m}^{\theta}\hat{\Psi}_{m:n}^{\theta} \big)(\varphi) }
	{\big( \hat{F}_{m-1:m}^{\theta}\hat{R}_{m:n}^{\theta} \big)(\boldsymbol 1) }.
\end{align}

Let $\rho_{5}=\max\{\sqrt{\rho_{1} },\sqrt{\rho_{4} } \}$,
$\tilde{C}_{1}=\max_{n\geq 1} n\rho_{5}^{n}$,
$\tilde{C}_{2}=2C_{4}\tilde{C}_{1}$,
$\tilde{C}_{3}=2C_{1}^{2}\tilde{C}_{2}$
($\rho_{1}$, $\rho_{4}$, $C_{1}$, $C_{4}$ are specified in
Proposition \ref{proposition3.2} and  Lemma \ref{lemma5.1}).
Since $\rho_{4}\leq\rho_{5}$, $\rho_{4}^{n-m}(n-m)\leq\tilde{C}_{1}\rho_{5}^{n-m}$,
Lemma \ref{lemma5.1} 
implies
\begin{align}\label{l5.2.57}
	\left\|
	\frac{\big( \xi\hat{\Phi}_{m:n}^{\theta} \big)(\varphi) }
	{\big( \xi\hat{R}_{m:n}^{\theta} \big)(\boldsymbol 1) }
	\!-\!
	\frac{\big( \xi'\hat{\Phi}_{m:n}^{\theta} \big)(\varphi) }
	{\big( \xi'\hat{R}_{m:n}^{\theta} \big)(\boldsymbol 1) }
	\right\|
	&\!\leq\!
	C_{4} \rho_{4}^{n-m}
	\left(1 + n - m + \rho_{4}^{m} \|w_{\theta} \| \right)
	\\ \nonumber
	&\!\leq\!
	2C_{4}\tilde{C}_{1}
	\left(\rho_{5}^{n-m} + \rho_{5}^{m} \|w_{\theta} \| \right)
	\\ \nonumber
	&\!\leq\!
	\tilde{C}_{2} \left(\rho_{5}^{n-m} + \rho_{5}^{n} \|w_{\theta} \| \right).
\end{align}
Using Proposition \ref{proposition2.1}
and (\ref{l5.2.7}), (\ref{l5.2.9}), (\ref{l5.2.23}),
(\ref{l5.2.53}) -- (\ref{l5.2.57}), we conclude
\begin{align}\label{l5.2.59}
	&
	\left\|
	E\left(\left.
	\hat{A}_{m:n}^{\theta}(\varphi)
	-
	\hat{G}_{m-1:n}^{\theta}(\varphi)
	\right| \boldsymbol Y, \hat{\xi}_{m-1}^{\theta}
	\right)
	\right\|
	\\ \nonumber
	&=
	\left\|
	E\left(\left.
	\frac{\big( \hat{\xi}_{m}^{\theta}\hat{\Phi}_{m:n}^{\theta} \big)(\varphi) }
	{\big( \hat{\xi}_{m}^{\theta}\hat{R}_{m:n}^{\theta} \big)(\boldsymbol 1) }
	-
	\frac{\big( \hat{F}_{m-1:m}^{\theta}\hat{\Phi}_{m:n}^{\theta} \big)(\varphi) }
	{\big( \hat{F}_{m-1:m}^{\theta}\hat{R}_{m:n}^{\theta} \big)(\boldsymbol 1) }
	\right| \boldsymbol Y, \hat{\xi}_{m-1}^{\theta}
	\right)
	\right\|
	\\ \nonumber
	&\leq
	\frac{2C_{1}^{2}\tilde{C}_{2} \left(\rho_{5}^{n-m} + \rho_{5}^{n} \|w_{\theta} \| \right) }
	{N}
	=
	\frac{\tilde{C}_{3} \left(\rho_{5}^{n-m} + \rho_{5}^{n} \|w_{\theta} \| \right)}
	{N}
\end{align}
almost surely.\footnote
{To get (\ref{l5.2.59}), (\ref{l5.2.61}), the following should be done:
In Proposition \ref{proposition2.1}, set $z=x$, $k=N$ and replace
$f(z)$, $g(z)$, $\xi_{k}(dz)$, $\xi(dz)$ with
$\big(\hat{\Phi}_{m:n}^{\theta}\varphi\big)(x)$,
$\big(\hat{R}_{m:n}^{\theta}\boldsymbol 1\big)(x)$,
$\hat{\xi}_{m}^{\theta}(dx)$, $\hat{F}_{m-1:m}^{\theta}(dx)$. }
Relying on the same arguments, we deduce
\begin{align}\label{l5.2.61}
	&
	E\left(\left.
	\left\|
	\hat{A}_{m:n}^{\theta}(\varphi)
	-
	\hat{G}_{m-1:n}^{\theta}(\varphi)
	\right\|^{2}
	\right| \boldsymbol Y, \hat{\xi}_{m-1}^{\theta}
	\right)
	\\ \nonumber
	&=
	E\left(\left.
	\left\|
	\frac{\big( \hat{\xi}_{m}^{\theta}\hat{\Phi}_{m:n}^{\theta} \big)(\varphi) }
	{\big( \hat{\xi}_{m}^{\theta}\hat{R}_{m:n}^{\theta} \big)(\boldsymbol 1) }
	-
	\frac{\big( \hat{F}_{m-1:m}^{\theta}\hat{\Phi}_{m:n}^{\theta} \big)(\varphi) }
	{\big( \hat{F}_{m-1:m}^{\theta}\hat{R}_{m:n}^{\theta} \big)(\boldsymbol 1) }
	\right\|^{2}
	\right| \boldsymbol Y, \hat{\xi}_{m-1}^{\theta}
	\right)
	\\ \nonumber
	&\leq
	\left(
	\frac{2C_{1}\tilde{C}_{2} \left(\rho_{5}^{n-m} + \rho_{5}^{n} \|w_{\theta} \| \right)}
	{\sqrt{N} }
	\right)^{2}
	\leq
	\left(
	\frac{\tilde{C}_{3}\left(\rho_{5}^{n-m} + \rho_{5}^{n} \|w_{\theta} \| \right)}
	{\sqrt{N} }
	\right)^{2}
\end{align}
almost surely.

Let $\tilde{C}_{4}=C_{5}\tilde{C}_{2}$
($C_{5}$ is defined in (i)).
Since $\tilde{C}_{1}\rho_{5}^{m-n}\geq n-m$,
Lemma \ref{lemma5.1} implies
\begin{align}\label{l5.2.201}
	\left\|
	\hat{H}_{m-1:n}^{\theta}(\varphi)
	\right\|
	\leq
	C_{4} \left(1 + n - m + \rho_{4}^{m} \|w_{\theta} \| \right)
	\leq &
	2C_{4} \tilde{C}_{1} \left(\rho_{5}^{m-n} + \rho_{5}^{m} \|w_{\theta} \| \right)
	\\ \nonumber
	\leq &
	\tilde{C}_{2} \left(\rho_{5}^{m-n} + \rho_{5}^{m} \|w_{\theta} \| \right).
\end{align}
As $\hat{H}_{m-1:n}^{\theta}(\boldsymbol 1)$
is measurable with respect to $\boldsymbol Y$, $\hat{\xi}_{m-1}^{\theta}$
and
$\rho_{1}^{n-m}\leq \rho_{5}^{2(n-m)}$,
(\ref{l5.2.25}), (\ref{l5.2.201}) yield
\begin{align}\label{l5.2.73}
	&
	\left\|
	E\left(\left(\left.
	\hat{F}_{m:n}^{\theta}(\varphi)
	-
	\hat{F}_{m-1:n}^{\theta}(\varphi)
	\right)
	\hat{H}_{m-1:n}^{\theta}(\boldsymbol 1)
	\right| \boldsymbol Y, \hat{\xi}_{m-1}^{\theta}
	\right)
	\right\|
	\\ \nonumber
	&\leq
	\left|
	E\left(\left.
	\hat{F}_{m:n}^{\theta}(\varphi)
	-
	\hat{F}_{m-1:n}^{\theta}(\varphi)
	\right| \boldsymbol Y, \hat{\xi}_{m-1}^{\theta}
	\right)
	\right|
	\left\|
	\hat{H}_{m-1:n}^{\theta}(\boldsymbol 1)
	\right\|
	\\ \nonumber
	&
	\leq
	\frac{C_{5}\tilde{C}_{2} \rho_{1}^{n-m}
	(\rho_{5}^{m-n} + \rho_{5}^{m}\|w_{\theta}\| ) }
	{N}
	\leq
	\frac{\tilde{C}_{4}
	(\rho_{5}^{n-m} + \rho_{5}^{n}\|w_{\theta}\| ) }
	{N}
\end{align}
almost surely.

Let $\tilde{C}_{5}=C_{5}\tilde{C}_{3}$
($C_{5}$ is defined in (i)).
Since $\tilde{C}_{1}\rho_{5}^{m-n}\geq n-m$,
Lemma \ref{lemma5.1} implies
\begin{align}\label{l5.2.63}
	\left\|
	\hat{C}_{m:n}^{\theta}(\varphi)
	\right\|
	\leq
	C_{4}
	\left(1 + n - m + \rho_{4}^{m} \|w_{\theta} \| \right)
	\leq &
	2C_{4}\tilde{C}_{1}
	\left(\rho_{5}^{m-n} + \rho_{5}^{m} \|w_{\theta} \| \right)
	\\ \nonumber
	\leq &
	\tilde{C}_{2} \left(\rho_{5}^{m-n} + \rho_{5}^{m} \|w_{\theta} \| \right).
\end{align}
Moreover, the same lemma yields
\begin{align}\label{l5.2.65}
	\left\|
	\frac{\big( \xi\hat{\Psi}_{m:n}^{\theta} \big)(\varphi) }
	{\big( \xi\hat{R}_{m:n}^{\theta} \big)(\boldsymbol 1) }
	\!-\!
	\frac{\big( \xi'\hat{\Psi}_{m:n}^{\theta} \big)(\varphi) }
	{\big( \xi'\hat{R}_{m:n}^{\theta} \big)(\boldsymbol 1) }
	\right\|
	&\!\leq\!
	C_{4}
	\left(1 + n - m + \rho_{4}^{m} \|w_{\theta} \| \right)
	\\ \nonumber
	&\!\leq\!
	\tilde{C}_{2} \left(\rho_{5}^{m-n} + \rho_{5}^{m} \|w_{\theta} \| \right).
\end{align}
Using Proposition \ref{proposition2.1}
and (\ref{l5.2.7}), (\ref{l5.2.9}), (\ref{l5.2.23}),
(\ref{l5.2.49}), (\ref{l5.2.53}), (\ref{l5.2.65}),
we conclude
\begin{align}\label{l5.2.67}
	&
	E\left(\left.
	\left\|
	\hat{C}_{m:n}^{\theta}(\varphi)
	-
	\hat{H}_{m-1:n}^{\theta}(\varphi)
	\right\|^{2}
	\right| \boldsymbol Y, \hat{\xi}_{m-1}^{\theta}
	\right)
	\\ \nonumber
	&=
	E\left(\left.
	\left\|
	\frac{\big( \hat{\xi}_{m}^{\theta}\hat{\Psi}_{m:n}^{\theta} \big)(\varphi) }
	{\big( \hat{\xi}_{m}^{\theta}\hat{R}_{m:n}^{\theta} \big)(\boldsymbol 1) }
	-
	\frac{\big( \hat{F}_{m-1:m}^{\theta}\hat{\Psi}_{m:n}^{\theta} \big)(\varphi) }
	{\big( \hat{F}_{m-1:m}^{\theta}\hat{R}_{m:n}^{\theta} \big)(\boldsymbol 1) }
	\right\|^{2}
	\right| \boldsymbol Y, \hat{\xi}_{m-1}^{\theta}
	\right)
	\\ \nonumber
	&\leq
	\left(
	\frac{2C_{1}\tilde{C}_{2} \left(\rho_{5}^{m-n} + \rho_{5}^{m} \|w_{\theta} \| \right)}
	{\sqrt{N} }
	\right)^{2}
	\leq
	\left(
	\frac{\tilde{C}_{3} \left(\rho_{5}^{m-n} + \rho_{5}^{m} \|w_{\theta} \| \right)}
	{\sqrt{N} }
	\right)^{2}
\end{align}
almost surely.\footnote
{To get (\ref{l5.2.67}), the following should be done:
In Proposition \ref{proposition2.1}, set $z=x$, $k=N$ and replace
$f(z)$, $g(z)$, $\xi_{k}(dz)$, $\xi(dz)$ with
$\big(\hat{\Psi}_{m:n}^{\theta}\varphi\big)(x)$,
$\big(\hat{R}_{m:n}^{\theta}\boldsymbol 1\big)(x)$,
$\hat{\xi}_{m}^{\theta}(dx)$, $\hat{F}_{m-1:m}^{\theta}(dx)$. }
As $\rho_{1}^{n-m}\leq \rho_{5}^{2(n-m)}$,
H\"{o}lder inequality and (\ref{l5.2.27}), (\ref{l5.2.67}) imply
\begin{align}\label{l5.2.71}
	&
	\left\|
	E\left(\left.
	\left(
	\hat{F}_{m:n}^{\theta}(\varphi)
	-
	\hat{F}_{m-1:n}^{\theta}(\varphi)
	\right)
	\left(
	\hat{C}_{m:n}^{\theta}(\boldsymbol 1)
	-
	\hat{H}_{m:n}^{\theta}(\boldsymbol 1)
	\right)
	\right| \boldsymbol Y, \hat{\xi}_{m-1}^{\theta}
	\right)
	\right\|
	\\ \nonumber
	&
	\begin{aligned}[b]
	\leq &
	\left(
	E\left(\left.
	\left|
	\hat{F}_{m:n}^{\theta}(\varphi)
	-
	\hat{F}_{m-1:n}^{\theta}(\varphi)
	\right|^{2}
	\right| \boldsymbol Y, \hat{\xi}_{m-1}^{\theta}
	\right)
	\right)^{1/2}
	\\
	&\cdot
	\left(
	E\left(\left.
	\left\|
	\hat{C}_{m:n}^{\theta}(\boldsymbol 1)
	-
	\hat{H}_{m:n}^{\theta}(\boldsymbol 1)
	\right\|^{2}
	\right| \boldsymbol Y, \hat{\xi}_{m-1}^{\theta}
	\right)
	\right)^{1/2}
	\end{aligned}
	\\ \nonumber
	&
	\leq
	\frac{C_{5}\tilde{C}_{3}\rho_{1}^{n-m}
	(\rho_{5}^{m-n} + \rho_{5}^{m}\|w_{\theta}\| ) }
	{N}
	\leq
	\frac{\tilde{C}_{5}
	(\rho_{5}^{n-m} + \rho_{5}^{n}\|w_{\theta}\| ) }
	{N}
\end{align}
almost surely.
Similarly, (\ref{l5.2.27}), (\ref{l5.2.63}) yield
\begin{align}\label{l5.2.75}
	&
	E\left(\left.
	\left\|
	\left(
	\hat{F}_{m:n}^{\theta}(\varphi)
	-
	\hat{F}_{m-1:n}^{\theta}(\varphi)
	\right)
	\hat{C}_{m:n}^{\theta}(\boldsymbol 1)
	\right\|^{2}
	\right| \boldsymbol Y, \hat{\xi}_{m-1}^{\theta}
	\right)
	\\ \nonumber
	&\leq
	\tilde{C}_{2}^{2} (\rho_{5}^{m-n} + \rho_{5}^{m} \|w_{\theta}\| )^{2}
	E\left(\left.
	\left|
	\hat{F}_{m:n}^{\theta}(\varphi)
	-
	\hat{F}_{m-1:n}^{\theta}(\varphi)
	\right|^{2}
	\right| \boldsymbol Y, \hat{\xi}_{m-1}^{\theta}
	\right)
	\\ \nonumber
	&
	\leq
	\left(
	\frac{C_{5}\tilde{C}_{2} \rho_{1}^{n-m}
	(\rho_{5}^{m-n} + \rho_{5}^{m}\|w_{\theta}\| ) }
	{\sqrt{N} }
	\right)^{2}
	\leq
	\left(
	\frac{\tilde{C}_{4}
	(\rho_{5}^{n-m} + \rho_{5}^{n}\|w_{\theta}\| ) }
	{\sqrt{N} }
	\right)^{2}
\end{align}
almost surely.

Let $\tilde{C}_{6}=\tilde{C}_{4}+\tilde{C}_{5}$.
Due to (\ref{5.27}), we have
\begin{align}\label{l5.2.77}
	\hat{B}_{m:n}^{\theta}(\varphi)
	=&
	-\left(
	\hat{F}_{m:n}^{\theta}(\varphi)
	-
	\hat{F}_{m-1:n}^{\theta}(\varphi)
	\right)
	\hat{C}_{m:n}^{\theta}(\boldsymbol 1)
	\\ \nonumber
	=&
	-\left(
	\hat{F}_{m:n}^{\theta}(\varphi)
	-
	\hat{F}_{m-1:n}^{\theta}(\varphi)
	\right)
	\left(
	\hat{C}_{m:n}^{\theta}(\boldsymbol 1)
	-
	\hat{H}_{m:n}^{\theta}(\boldsymbol 1)
	\right)
	\\ \nonumber
	&
	-\left(
	\hat{F}_{m:n}^{\theta}(\varphi)
	-
	\hat{F}_{m-1:n}^{\theta}(\varphi)
	\right)
	\hat{H}_{m:n}^{\theta}(\boldsymbol 1).
\end{align}
Consequently,  (\ref{l5.2.73}), (\ref{l5.2.71}) and the second part of (\ref{l5.2.77}) imply
\begin{align}\label{l5.2.501}
	\left\|
	E\left(\left.
	\hat{B}_{m:n}^{\theta}(\varphi)
	\right| \boldsymbol Y, \hat{\xi}_{m-1}^{\theta}
	\right)
	\right\|
	\leq
	\frac{(\tilde{C}_{4}+\tilde{C}_{5} )
	(\rho_{5}^{n-m} + \rho_{5}^{n}\|w_{\theta}\| ) }
	{N}
	=
	\frac{\tilde{C}_{6}	(\rho_{5}^{n-m} + \rho_{5}^{n}\|w_{\theta}\| ) }
	{N}
\end{align}
almost surely.
Moreover, the first part of (\ref{l5.2.77}) and (\ref{l5.2.75}) yield
\begin{align}\label{l5.2.503}
	E\left(\left.
	\!\left\|
	\hat{B}_{m:n}^{\theta}(\varphi)
	\right\|^{2}
	\right| \boldsymbol Y, \hat{\xi}_{m-1}^{\theta}
	\!\right)
	\leq
	\!\left(
	\frac{\tilde{C}_{4}
	(\rho_{5}^{n-m} + \rho_{5}^{n}\|w_{\theta}\| ) }
	{\sqrt{N} }
	\right)^{\!\!2}
	\!\leq
	\!\left(
	\frac{\tilde{C}_{6}
	(\rho_{5}^{n-m} + \rho_{5}^{n}\|w_{\theta}\| ) }
	{\sqrt{N} }
	\right)^{\!\!2}
\end{align}
almost surely.

Let $C_{6}=\tilde{C}_{3}+\tilde{C}_{6}$.
Then, combining (\ref{l5.2.59}), (\ref{l5.2.61}), (\ref{l5.2.501}), (\ref{l5.2.503})
with the tower property of conditional expectations,
we conclude that (\ref{l5.2.1*}) -- (\ref{l5.2.7*}) hold almost surely.
\end{proof}

\begin{proposition}\label{proposition5.1}
Let $\theta$ be any element of $\Theta$, 
while $\boldsymbol y = \{y_{n} \}_{n\geq 0}$ is any sequence in ${\cal Y}$. 
Moreover, let $\varphi:{\cal X}\rightarrow[-1,1]$ be any Borel-measurable function, 
while $n$ is any non-negative integer.  

(i) Suppose that Assumption \ref{a1} holds.
Then, there exists a real number $L\in[1,\infty)$
(independent of $N$, $\theta$, $\boldsymbol y$, $\varphi(x)$, $n$ 
and depending only on $\varepsilon$) such that
\begin{align}
	&\label{p5.1.1*}
	\left|
	E\left(\left.
	\hat{\xi}_{n}^{\theta}(\varphi)
	-
	F_{\theta,\boldsymbol Y}^{0:n}(\xi_{\theta} )(\varphi)
	\right|
	\boldsymbol Y =\boldsymbol y
	\right)
	\right|
	\leq
	\frac{L}{N},
	\\
	&\label{p5.1.3*}
	\left(
	E\left(\left.
	\left|
	\hat{\xi}_{n}^{\theta}(\varphi)
	-
	F_{\theta,\boldsymbol Y}^{0:n}(\xi_{\theta} )(\varphi)
	\right|^{2}
	\right|
	\boldsymbol Y = \boldsymbol y
	\right)
	\right)^{1/2}
	\leq
	\frac{L}{\sqrt{N} }.
\end{align}

(ii) Suppose that Assumptions \ref{a1} -- \ref{a3} hold.
Then, there exist real numbers $\rho\in(0,1)$, $M\in[1,\infty)$
(independent of $N$, $\theta$, $\boldsymbol y$, $\varphi(x)$, $n$ 
and depending only on $\varepsilon$, $d$, $K$) such that
\begin{align}
	&\label{p5.1.5*}
	\left\|
	E\left(\left.
	\hat{\zeta}_{n}^{\theta}(\varphi)
	-
	G_{\theta,\boldsymbol Y}^{0:n}(\xi_{\theta},\zeta_{\theta} )(\varphi) 
	\right|
	\boldsymbol Y = \boldsymbol y
	\right)
	\right\|
	\leq
	\frac{M(1 + \rho^{n} \|w_{\theta}\| ) }{N},
	\\
	&\label{p5.1.7*}
	\left(
	E\left(\left.
	\left\|
	\hat{\zeta}_{n}^{\theta}(\varphi)
	-
	G_{\theta,\boldsymbol Y}^{0:n}(\xi_{\theta},\zeta_{\theta} )(\varphi)
	\right\|^{2}
	\right|
	\boldsymbol Y = \boldsymbol y
	\right)
	\right)^{1/2}
	\leq
	\frac{M(1 + \rho^{n} \|w_{\theta}\| ) }{\sqrt{N} }. 
\end{align}
\end{proposition}

\begin{vremark}
Relying on (\ref{1.1'}), (\ref{1.3'}), (\ref{1.21}), (\ref{1.25}), 
it is easy to show 
\begin{align*}
	P_{\theta,\boldsymbol Y}^{n}(B)
	=F_{\theta,\boldsymbol Y}^{0:n}\big(\xi_{\theta} \big)(B), 
	\;\;\;\;\; 
	Q_{\theta,\boldsymbol Y}^{n}(B)
	=G_{\theta,\boldsymbol Y}^{0:n}\big(\xi_{\theta},\zeta_{\theta} \big)(B)
\end{align*}
for $B\in{\cal B}({\cal X})$, $n\geq 1$. 
Hence, Proposition \ref{proposition5.1} can be considered as an extended version of Theorem \ref{theorem1.1}. 
Moreover, the bounds in (\ref{p5.1.3*}), (\ref{p5.1.7*}) can be viewed as 
by-products of Theorem \ref{theorem1.1}.
Under the same conditions as in Proposition \ref{proposition5.1},
bounds similar to (\ref{p5.1.3*}), (\ref{p5.1.7*}) have been derived in \cite{delmoral&doucet&singh2}.
\end{vremark}

\begin{proof}
(i) Let $L=C_{5}(1-\rho_{1} )^{-1}$
($\rho_{1}$ is specified in Proposition \ref{proposition3.2}). 
Using (\ref{5.5}), (\ref{5.7.1}),
it is straightforward to verify
\begin{align*}
	\hat{\xi}_{n}^{\theta}(\varphi)
	=\hat{F}_{n:n}^{\theta}(\varphi),
	\;\;\;\;\;
	F_{\theta,\boldsymbol Y}^{0:n}\big(\xi_{\theta} \big)(\varphi)
	=F_{\theta,\boldsymbol Y}^{0:n}\big(\hat{\xi}_{-1}^{\theta} \big)(\varphi)
	=\hat{F}_{-1:n}^{\theta}(\varphi).
\end{align*}
Therefore, we get
\begin{align}\label{p5.1.1}
	\hat{\xi}_{n}^{\theta}(\varphi)
	-
	F_{\theta,\boldsymbol Y}^{0:n}\big(\xi_{\theta} \big)(\varphi)
	=
	\hat{F}_{n:n}^{\theta}(\varphi)
	-
	\hat{F}_{-1:n}^{\theta}(\varphi)
	=
	\sum_{m=0}^{n}
	\left(
	\hat{F}_{m:n}^{\theta}(\varphi)
	-
	\hat{F}_{m-1:n}^{\theta}(\varphi)
	\right).
\end{align}
Then, Lemma \ref{lemma5.2} implies
\begin{align}\label{p5.1.3}
	\left|
	E\left(\left.
	\hat{\xi}_{n}^{\theta}(\varphi)
	-
	F_{\theta,\boldsymbol Y}^{0:n}\big(\xi_{\theta} \big)(\varphi)
	\right|
	\boldsymbol Y
	\right)
	\right|
	\leq &
	\sum_{m=0}^{n}
	\left|
	E\left(\left.
	\hat{F}_{m:n}^{\theta}(\varphi)
	-
	\hat{F}_{m-1:n}^{\theta}(\varphi)
	\right|
	\boldsymbol Y
	\right)
	\right|
	\\ \nonumber
	\leq &
	\frac{C_{5} }{N} \sum_{m=0}^{n} \rho_{1}^{n-m}
	\leq
	\frac{L}{N}
\end{align}
almost surely.
Moreover, Minkowski inequality, Lemma \ref{lemma5.2} and (\ref{p5.1.1}) yield
\begin{align}\label{p5.1.5}
	\left(
	E\left(\left.
	\left|
	\hat{\xi}_{n}^{\theta}(\varphi)
	\!-\!
	F_{\theta,\boldsymbol Y}^{0:n}\big(\xi_{\theta} \big)(\varphi)
	\right|^{2}
	\right|
	\boldsymbol Y
	\right)
	\right)^{\!1/2}
	\leq &
	\sum_{m=0}^{n}
	\left(
	E\left(\left.
	\left|
	\hat{F}_{m:n}^{\theta}(\varphi)
	\!-\!
	\hat{F}_{m-1:n}^{\theta}(\varphi)
	\right|^{2}
	\right|
	\boldsymbol Y
	\right)
	\right)^{\!1/2}
	\\ \nonumber
	\leq &
	\frac{C_{5}}{\sqrt{N} } \sum_{m=0}^{n} \rho_{1}^{n-m}
	\leq
	\frac{L}{\sqrt{N} }
\end{align}
almost surely.
Using (\ref{p5.1.3}), (\ref{p5.1.5}),
we conclude that (\ref{p5.1.1*}), (\ref{p5.1.3*}) hold.

(ii) Let $\rho\!=\!\sqrt{\rho_{5} }$,
$\tilde{C}=\max_{n\geq 1} n\rho^{n}$, while
$M=4C_{6}\tilde{C}(1-\rho )^{-1}$
($\rho_{5}$, $C_{6}$ are specified in
Lemma \ref{lemma5.2}).
Owing to Lemmas \ref{lemma5.3}, \ref{lemma5.2}, we have
\begin{align}\label{p5.1.23}
	\left\|
	E\left(\left.
	\hat{G}_{m:n}^{\theta}(\varphi)
	-
	\hat{G}_{m-1:n}^{\theta}(\varphi)
	\right| \boldsymbol Y
	\right)
	\right\|
	\leq &
	\left\|
	E\left(\left.
	\hat{A}_{m:n}^{\theta}(\varphi)
	-
	\hat{G}_{m-1:n}^{\theta}(\varphi)
	\right| \boldsymbol Y
	\right)
	\right\|
	\\ \nonumber
	&+
	\left\|
	E\left(\left.
	\hat{B}_{m:n}^{\theta}(\varphi)
	\right| \boldsymbol Y
	\right)
	\right\|
	\\ \nonumber
	\leq &
	\frac{2C_{6}(\rho_{5}^{n-m} + \rho_{5}^{n} \|w_{\theta}\| ) }{N}
\end{align}
almost surely.
Similarly, due to Minkowski inequality and Lemmas \ref{lemma5.3}, \ref{lemma5.2},
we have
\begin{align}\label{p5.1.25}
	\left(
	E\!\left(\left.
	\left\|
	\hat{G}_{m:n}^{\theta}(\varphi)
	-
	\hat{G}_{m-1:n}^{\theta}(\varphi)
	\right\|^{2}
	\right| \!\boldsymbol Y
	\!\right)
	\right)^{\!\!1/2}
	\leq &
	\left(
	E\!\left(\left.
	\left\|
	\hat{A}_{m:n}^{\theta}(\varphi)
	-
	\hat{G}_{m-1:n}^{\theta}(\varphi)
	\right\|^{2}
	\right| \!\boldsymbol Y
	\!\right)
	\right)^{\!\!1/2}
	\\ \nonumber
	&+
	\left(
	E\!\left(\left.
	\left\|
	\hat{B}_{m:n}^{\theta}(\varphi)
	\right\|^{2}
	\right| \!\boldsymbol Y
	\!\right)
	\right)^{\!\!1/2}
	\\ \nonumber
	\leq &
	\frac{2C_{6}(\rho_{5}^{n-m} + \rho_{5}^{n} \|w_{\theta}\| ) }{\sqrt{N} }
\end{align}
almost surely.

Using (\ref{5.5}), (\ref{5.7.3}),
it is straightforward to verify
\begin{align*}
	\hat{\zeta}_{n}^{\theta}(\varphi)
	=\hat{G}_{n:n}^{\theta}(\varphi),
	\;\;\;\;\;
	G_{\theta,\boldsymbol Y}^{0:n}\big(\xi_{\theta},\zeta_{\theta} \big)(\varphi)
	=G_{\theta,\boldsymbol Y}^{0:n}\big(\hat{\xi}_{-1}^{\theta}, \hat{\zeta}_{-1}^{\theta} \big)(\varphi)
	=\hat{G}_{-1:n}^{\theta}(\varphi).
\end{align*}
Therefore, we get
\begin{align}\label{p5.1.29}
	\hat{\zeta}_{n}^{\theta}(\varphi)
	\!-
	G_{\theta,\boldsymbol Y}^{0:n}\big(\xi_{\theta},\zeta_{\theta} \big)(\varphi)
	\!=&
	\hat{G}_{n:n}^{\theta}(\varphi)
	\!-
	G_{-1:n}^{\theta}(\varphi)
	\!=
	\sum_{m=0}^{n}
	\left(
	\hat{G}_{m:n}^{\theta}(\varphi)
	\!-
	\hat{G}_{m-1:n}^{\theta}(\varphi)
	\right).
\end{align}
Then, (\ref{p5.1.23}) implies
\begin{align}\label{p5.1.31}
	&
	\left\|
	E\left(\left.
	\hat{\zeta}_{n}^{\theta}(\varphi)
	-
	G_{\theta,\boldsymbol Y}^{0:n}\big(\xi_{\theta},\zeta_{\theta} \big)(\varphi)
	\right| \boldsymbol Y
	\right)
	\right\|
	\leq
	\sum_{m=0}^{n}
	\left\|
	E\left(\left.
	\hat{G}_{m:n}^{\theta}(\varphi)
	-
	\hat{G}_{m-1:n}^{\theta}(\varphi)
	\right| \boldsymbol Y
	\right)
	\right\|
	\\ \nonumber
	&\leq
	\frac{2C_{6} }{N}
	\sum_{m=0}^{n} (\rho_{5}^{n-m} + \rho_{5}^{n} \|w_{\theta}\| )
	\leq
	\frac{2C_{6} }{(1-\rho_{5} ) N}
	+
	\frac{2C_{6} (n+1) \rho_{5}^{n} \|w_{\theta} \| }{N}
	\leq
	\frac{M (1 + \rho^{n} \|w_{\theta} \| ) }{N}
\end{align}
almost surely.
Moreover, Minkowski inequality
and (\ref{p5.1.25}) -- (\ref{p5.1.29}) yield
\begin{align}\label{p5.1.33}
	&
	\left(\!
	E\!\left(\left.
	\left\|
	\hat{\zeta}_{n}^{\theta}(\varphi)
	\!-\!
	G_{\theta,\boldsymbol Y}^{0:n}\big(\xi_{\theta},\zeta_{\theta} \big)(\varphi)
	\right\|^{2}
	\right| \boldsymbol Y\!
	\right)\!
	\right)^{\!\!1/2}
	\!\!\leq\!
	\sum_{m=0}^{n}\!
	\left(\!
	E\!\left(\left.
	\left\|
	\hat{G}_{m:n}^{\theta}(\varphi)
	\!-\!
	\hat{G}_{m-1:n}^{\theta}(\varphi)
	\right\|^{2}
	\right| \boldsymbol Y
	\right)\!
	\right)^{\!\!1/2}
	\\ \nonumber
	&\leq
	\frac{2C_{6} }{\sqrt{N} }
	\sum_{m=0}^{n} (\rho_{5}^{n-m} + \rho_{5}^{n} \|w_{\theta}\| )
	\!\leq\!
	\frac{2C_{6} }{(1-\rho_{5} ) \sqrt{N} }
	+
	\frac{2C_{6} (n+1) \rho_{5}^{n} \|w_{\theta} \| }
	{\sqrt{N} }
	\!\leq\!
	\frac{M(1 + \rho^{n} \|w_{\theta} \| ) }{\sqrt{N} }
	\\ \nonumber
\end{align}
almost surely.
Using (\ref{p5.1.31}), (\ref{p5.1.33}),
we conclude that (\ref{p5.1.5*}), (\ref{p5.1.7*}) hold.
\end{proof}

\end{document}


\maketitle

\section{\nopunct}

In this section, we prove Proposition 3.1. 
The proof is based on the Taylor formula and well-known results 
on sums of independent random variables. 

\begin{proof}[Proof of Proposition 3.1]
Let $\xi'(dz)$, $\xi''(dz)$ be any probability measures on ${\cal Z}$. 
As a direct consequence of (3.2), we get  
\begin{align*}
	g(z')\leq \beta g(z''), 
	\;\;\;\;\; 
	f(z') g(z'') 
	\leq 
	(f(z'')+\alpha g(z'') ) g(z')
\end{align*}
for $z',z''\in{\cal Z}$. 
Then, we have  
\begin{align*}
	\xi'(g) 
	= 
	\int\int g(z') \xi'(dz')\xi''(dz'') 
	\leq 
	\beta \int\int g(z'') \xi'(dz')\xi''(dz'')
	=
	\beta \xi''(g). 
\end{align*}
We also have 
\begin{align*}
	\xi'(f)\xi''(g)
	\!=\!
	\int\int f(z')g(z'') \xi'(dz')\xi''(dz'') 
	\!\leq\! &
	\int\int (f(z'')+\alpha g(z'') ) g(z') \xi'(dz')\xi''(dz'')
	\\
	\!=\! &
	(\xi''(f)+\alpha\xi''(g) ) \xi'(g).
\end{align*}
Consequently, we get 
\begin{align}\label{p2.1.1}
	\frac{\xi'(g)}{\xi''(g)}
	\leq 
	\beta, 
	\;\;\;\;\; 
	\frac{\xi'(f)}{\xi'(g)}
	\leq 
	\frac{\xi''(f)}{\xi''(g)}
	+
	\alpha. 
\end{align}
Reverting the roles of $\xi'(dz)$, $\xi''(dz)$, 
we conclude $\xi''(f)/\xi''(g)\leq\xi'(f)/\xi'(g)+\alpha$. 
Since $\beta\geq 1$, we then deduce  
\begin{align}\label{p2.1.3}
	\left|
	\frac{\xi'(f)}{\xi'(g)} - \frac{\xi''(f)}{\xi''(g)} 
	\right|
	\leq
	\alpha, 
	\;\;\;\;\; 
	\left|
	\frac{\xi'(g)}{\xi''(g)} - 1
	\right|
	\leq 
	\frac{\xi'(g)}{\xi''(g)} + 1
	\leq 
	2\beta. 
\end{align}

Let $k\geq 1$ be any integer. It is straightforward to verify 
\begin{align*}
	\frac{\xi_{k}(f) }{\xi_{k}(g) }
	=&
	\frac{\xi_{k}(f) }{\xi(g) }
	-
	\frac{\xi_{k}(f) }{\xi(g) }
	\left( \frac{\xi_{k}(g) }{\xi(g) } - 1 \right)
	+
	\frac{\xi(f) }{\xi(g) }
	\left( \frac{\xi_{k}(g) }{\xi(g) } - 1 \right)^{2}
	\\
	&+
	\left(\frac{\xi_{k}(f) }{\xi_{k}(g) } - \frac{\xi(f) }{\xi(g) } \right)
	\left( \frac{\xi_{k}(g) }{\xi(g) } - 1 \right)^{2}. 
\end{align*}
Consequently, we get 
\begin{align*}
	E\left(\frac{\xi_{k}(f) }{\xi_{k}(g) } \right)
	-
	\frac{\xi(f) }{\xi(g) }
	=&
	-
	E\left( 
	\frac{\xi_{k}(f) }{\xi(g) }
	\left( \frac{\xi_{k}(g) }{\xi(g) } - 1 \right)
	\right)
	+
	E\left(
	\frac{\xi(f) }{\xi(g) }
	\left( \frac{\xi_{k}(g) }{\xi(g) } - 1 \right)^{2}
	\right)
	\\ \nonumber
	&+
	E\left(
	\left(\frac{\xi_{k}(f) }{\xi_{k}(g) } - \frac{\xi(f) }{\xi(g) } \right)
	\left( \frac{\xi_{k}(g) }{\xi(g) } - 1 \right)^{2}
	\right). 
\end{align*}
Then, owing to (\ref{p2.1.3}), we have 
\begin{align}\label{p2.1.9}
	\left|
	E\left(\frac{\xi_{k}(f) }{\xi_{k}(g) } \right)
	-
	\frac{\xi(f) }{\xi(g) }
	\right|
	\leq & 
	\left|
	E\left( 
	\frac{\xi_{k}(f) }{\xi(g) }
	\left( \frac{\xi_{k}(g) }{\xi(g) } - 1 \right)
	\right)
	-
	E\left(
	\frac{\xi(f) }{\xi(g) }
	\left( \frac{\xi_{k}(g) }{\xi(g) } - 1 \right)^{2}
	\right)
	\right|
	\\ \nonumber
	&+
	E\left(
	\left|\frac{\xi_{k}(f) }{\xi_{k}(g) } - \frac{\xi(f) }{\xi(g) } \right|
	\left( \frac{\xi_{k}(g) }{\xi(g) } - 1 \right)^{2}
	\right)
	\\ \nonumber
	\leq &
	\left|
	E\left( 
	\frac{\xi_{k}(f) }{\xi(g) }
	\left( \frac{\xi_{k}(g) }{\xi(g) } - 1 \right)
	\right)
	-
	E\left(
	\frac{\xi(f) }{\xi(g) }
	\left( \frac{\xi_{k}(g) }{\xi(g) } - 1 \right)^{2}
	\right)
	\right|
	\\ \nonumber
	&+
	\alpha 
	E\left(
	\left( \frac{\xi_{k}(g) }{\xi(g) } - 1 \right)^{2}
	\right). 
\end{align}
As $\{Z_{k} \}_{k\geq 1}$ are distributed according to $\xi(dz)$, we also have 
\begin{align}\label{p2.1.21}
	E\left(\frac{g(Z_{k} ) }{\xi(g) } - 1 \right) = 0. 
\end{align}
Since $\{Z_{k} \}_{k\geq 1}$ are independent, (\ref{p2.1.21}) implies 
\begin{align}\label{p2.1.23}
	E\left( 
	\frac{\xi_{k}(f) }{\xi(g) }
	\left( \frac{\xi_{k}(g) }{\xi(g) } - 1 \right)
	\right) 
	=&
	\frac{1}{k^{2} } \sum_{i=1}^{k} 
	E\left(\frac{f(Z_{i} ) }{\xi(g) } 
	\left(\frac{g(Z_{i} ) }{\xi(g) } - 1 \right)
	\right)
	\\ \nonumber
	&+
	\frac{1}{k^{2} } \sum_{\stackrel{\scriptstyle 1\leq i,j \leq k}{i\neq j} } 
	E\left(\frac{f(Z_{i} ) }{\xi(g) } \right)
	E\left(\frac{g(Z_{j} ) }{\xi(g) } - 1 \right)
	\\ \nonumber
	=&
	\frac{1}{k} 
	\int \frac{f(z)}{\xi(g) } \left(\frac{g(z)}{\xi(g)} - 1 \right)
	\xi(dz)
	\\ \nonumber
	=&
	\frac{1}{k} 
	\left(
	\int \frac{f(z) g(z) }{\xi^{2}(g) } \xi(dz)
	-
	\frac{\xi(f)}{\xi(g) }
	\right). 
\end{align}
For the same reason, (\ref{p2.1.21}) yields  
\begin{align}\label{p2.1.25}
	E\left( 
	\left( \frac{\xi_{k}(g) }{\xi(g) } - 1 \right)^{2}
	\right) 
	=&
	\frac{1}{k^{2} } \sum_{i=1}^{k} 
	E\left( 
	\left(\frac{g(Z_{i} ) }{\xi(g) } - 1 \right)^{2} 
	\right)
	\\ \nonumber
	&+
	\frac{1}{k^{2} } \sum_{\stackrel{\scriptstyle 1\leq i,j \leq k}{i\neq j} } 
	E\left(\frac{g(Z_{i} ) }{\xi(g) } - 1 \right)
	E\left(\frac{g(Z_{j} ) }{\xi(g) } - 1 \right)
	\\ \nonumber
	=&
	\frac{1}{k} 
	\int \left(\frac{g(z)}{\xi(g)} - 1 \right)^{2}
	\xi(dz)
	\\ \nonumber
	=&
	\frac{1}{k} 
	\left(
	\int \frac{g^{2}(z) }{\xi^{2}(g) } \xi(dz)
	-
	1
	\right). 
\end{align}
Using (\ref{p2.1.23}), (\ref{p2.1.25}), we deduce 
\begin{align*}
	&
	E\left( 
	\frac{\xi_{k}(f) }{\xi(g) }
	\left( \frac{\xi_{k}(g) }{\xi(g) } - 1 \right)
	\right)
	-
	E\left(
	\frac{\xi(f) }{\xi(g) }
	\left( \frac{\xi_{k}(g) }{\xi(g) } - 1 \right)^{2}
	\right)
	\\
	&=
	\frac{1}{k}
	\int \frac{g^{2}(z) }{\xi^{2}(g) }
	\left(\frac{f(z) }{g(z) } - \frac{\xi(f) }{\xi(g) } \right)
	\xi(dz). 
\end{align*}
Combining this with (\ref{p2.1.1}), (\ref{p2.1.3}), we get 
\begin{align}\label{p2.1.27}
	&
	\left|
	E\left( 
	\frac{\xi_{k}(f) }{\xi(g) }
	\left( \frac{\xi_{k}(g) }{\xi(g) } - 1 \right)
	\right)
	-
	E\left(
	\frac{\xi(f) }{\xi(g) }
	\left( \frac{\xi_{k}(g) }{\xi(g) } - 1 \right)^{2}
	\right)
	\right|
	\\ \nonumber
	&\leq 
	\frac{1}{k}
	\int \frac{g^{2}(z) }{\xi^{2}(g) }
	\left|\frac{f(z) }{g(z) } - \frac{\xi(f) }{\xi(g) } \right|
	\xi(dz)
	\leq 
	\frac{\alpha\beta^{2} }{k}. 
\end{align}
Similarly, combining (\ref{p2.1.1}), (\ref{p2.1.25}), we get 
\begin{align}\label{p2.1.29}
	E\left( 
	\left( \frac{\xi_{k}(g) }{\xi(g) } - 1 \right)^{2}
	\right) 
	\leq 
	\frac{1}{k} 
	\int \frac{g^{2}(z) }{\xi^{2}(g) } \xi(dz)
	\leq 
	\frac{\beta^{2} }{k}. 
\end{align}
Relying on (\ref{p2.1.9}), (\ref{p2.1.27}), (\ref{p2.1.29}), 
we conclude that the first part of (3.1) is true. 

It is straightforward to verify 
\begin{align*}
	\frac{\xi_{k}(f) }{\xi_{k}(g) }
	=
	\frac{\xi(f) }{\xi(g) }
	+
	\left(
	\frac{\xi_{k}(f) }{\xi(g) } - \frac{\xi(f)\xi_{k}(g) }{\xi^{2}(g) }
	\right)
	-
	\left(\frac{\xi_{k}(f) }{\xi_{k}(g) } - \frac{\xi(f) }{\xi(g) } \right)
	\left( \frac{\xi_{k}(g) }{\xi(g) } - 1 \right). 
\end{align*}
Then, Minkowski inequality and (\ref{p2.1.3}) imply 
\begin{align}\label{p2.1.31}
	\left(
	E\left(
	\left| 
	\frac{\xi_{k}(f) }{\xi_{k}(g) } \!- \frac{\xi(f) }{\xi(g) } 
	\right|^{2}
	\right)
	\right)^{\!\!1/2}
	\!\leq &
	\left(
	E\left(
	\left| \frac{\xi_{k}(f) }{\xi(g) } - \frac{\xi(f)\xi_{k}(g) }{\xi^{2}(g) } \right|^{2}
	\right)
	\right)^{\!\!1/2}
	\\ \nonumber
	&+
	\left(
	E\left(
	\left|\frac{\xi_{k}(f) }{\xi_{k}(g) } - \frac{\xi(f) }{\xi(g) } \right|^{2}
	\left| \frac{\xi_{k}(g) }{\xi(g) } - 1 \right|^{2} 
	\right)
	\right)^{\!\!1/2}
	\\ \nonumber
	\!\leq & 
	\left(
	E\left(
	\left| \frac{\xi_{k}(f) }{\xi(g) } - \frac{\xi(f)\xi_{k}(g) }{\xi^{2}(g) } \right|^{2}
	\right)
	\right)^{\!1/2}
	\\ \nonumber
	&+
	\alpha 
	\left(
	E\left(
	\left| \frac{\xi_{k}(g) }{\xi(g) } - 1 \right|^{2} 
	\right)
	\right)^{\!1/2}. 
\end{align}
As $\{Z_{k} \}_{k\geq 1}$ are distributed according to $\xi(dz)$, we have 
\begin{align}\label{p2.1.33}
	E\left(
	\frac{f(Z_{k} ) }{\xi(g) } - \frac{\xi(f)g(Z_{k} ) }{\xi^{2}(g) }
	\right)
	=
	0. 
\end{align}
Since $\{Z_{k} \}_{k\geq 1}$ are independent, 
(\ref{p2.1.1}), (\ref{p2.1.3}), (\ref{p2.1.33}) yield 
\begin{align}\label{p2.1.35}
	E\left(
	\left| \frac{\xi_{k}(f) }{\xi(g) } - \frac{\xi(f)\xi_{k}(g) }{\xi^{2}(g) } \right|^{2}
	\right)
	=&
	\frac{1}{k^{2} } 
	\sum_{i=1}^{k} 
	E\left(
	\left|\frac{f(Z_{i} ) }{\xi(g) } - \frac{\xi(f)g(Z_{i} ) }{\xi^{2}(g) } \right|^{2}
	\right)
	\\ \nonumber
	&
	\begin{aligned}[t]
	+ 
	\frac{1}{k^{2} } \sum_{\stackrel{\scriptstyle 1\leq i,j \leq k}{i\neq j} } 
	&
	E\left(
	\frac{f(Z_{i} ) }{\xi(g) } - \frac{\xi(f)g(Z_{i} ) }{\xi^{2}(g) } 
	\right)
	\\ 
	&\cdot 
	E\left(
	\frac{f(Z_{j} ) }{\xi(g) } - \frac{\xi(f)g(Z_{j} ) }{\xi^{2}(g) } 
	\right)
	\end{aligned} 
	\\ \nonumber
	= &
	\frac{1}{k} 
	\int \left|\frac{g(z)}{\xi(g)} \right|^{2} 
	\left| \frac{f(z)}{g(z)} - \frac{\xi(f) }{\xi(g) } \right|^{2} \xi(dz) 
	\\ \nonumber
	\leq &
	\frac{\alpha^{2}\beta^{2} }{k}.
\end{align}
Relying on (\ref{p2.1.29}), (\ref{p2.1.31}), (\ref{p2.1.35}), 
we conclude that the second part of (3.1) is true. 
\end{proof}

\section{\nopunct}

In this section, we prove Proposition 4.1. 
The proposition is based on the results of \cite{legland&oudjane} -- \cite{tadic&doucet3}. 
To show Proposition 4.1, 
we use the following additional notation. 
If $R(x,dx')$ and $S(x,dx')$ are integral operators from 
${\cal F}({\cal X})$ to ${\cal F}^{k}({\cal X})$, ${\cal F}^{l}({\cal X})$ (respectively)
and if one of $k$, $l$ is one, 
then $(RS)(x,dx')$ denotes the integral operator defined by 
\begin{align*}
	(RS)(x,B)=\int S(x',B)R(x,dx')
\end{align*}
for $B\in{\cal B}({\cal X})$. 
$\tilde{r}_{\theta,\boldsymbol y}^{n}\!(x'|x)$ is the function defined by 
\begin{align}\label{3.1}
	\tilde{r}_{\theta,\boldsymbol y}^{n}(x'|x)
	=
	q_{\theta}(y_{n}|x') p_{\theta}(x'|x)
\end{align}
for $\theta\in\Theta$, $x,x'\in{\cal X}$, $n\geq 0$ 
and a sequence $\boldsymbol y = \{y_{n} \}_{n\geq 0}$ in ${\cal Y}$. 
$\tilde{r}_{\theta,\boldsymbol y}^{m:n}(x_{m:n} )$ and 
$\tilde{s}_{\theta,\boldsymbol y}^{m:n}(x_{m:n} )$
are the functions defined by 
\begin{align}
	\label{3.3}
	\tilde{r}_{\theta,\boldsymbol y}^{m:m}(x_{m:m} ) = 1, 
	&\;\;\;\;\; 
	\tilde{r}_{\theta,\boldsymbol y}^{m:n}(x_{m:n} )
	=
	\prod_{k=m+1}^{n} 
	\tilde{r}_{\theta,\boldsymbol y}^{k}(x_{k}|x_{k-1} ),
	\\
	\label{3.5}
	\tilde{s}_{\theta,\boldsymbol y}^{m:m}(x_{m:m} ) = 0, 
	&\;\;\;\;\; 
	\tilde{s}_{\theta,\boldsymbol y}^{m:n}(x_{m:n} )
	=
	\nabla_{\theta} \tilde{r}_{\theta,\boldsymbol y}^{m:n}(x_{m:n} )
\end{align}
for $x_{m},\dots,x_{n}\in{\cal X}$, $n>m\geq 0$. 
$\tilde{A}_{\theta}(x,dx')$ and $\tilde{B}_{\theta}(x,dx')$
are the integral operators from ${\cal F}({\cal X})$ to 
${\cal F}({\cal X})$, ${\cal F}^{d}({\cal X})$ (respectively) defined by 
\begin{align}\label{3.51}
	\tilde{A}_{\theta}(x,B)=\int_{B}p_{\theta}(x'|x)\mu(dx'), 
	\;\;\;\;\; 
	\tilde{B}_{\theta}(x,B)=\int_{B}\nabla_{\theta}p_{\theta}(x'|x)\mu(dx'). 
\end{align}
$\tilde{R}_{\theta,\boldsymbol y}^{m:n}(x,dx')$ and 
$\tilde{S}_{\theta,\boldsymbol y}^{m:n}(x,dx')$
are the integral operators from ${\cal F}({\cal X})$ to 
${\cal F}({\cal X})$, ${\cal F}^{d}({\cal X})$ (respectively) defined by 
\begin{align}
	\label{3.7}
	\tilde{R}_{\theta,\boldsymbol y}^{m:m}(x,B)
	=
	\delta_{x}(B), 
	&\;\;\;\;\; 
	\begin{aligned}[t]
	\tilde{R}_{\theta,\boldsymbol y}^{m:n}(x,B)
	=&
	\int_{{\cal X}^{n-m}\times B} 
	\tilde{r}_{\theta,\boldsymbol y}^{m:n}(x_{m:n} ) 
	(\delta_{x}\times\mu^{n-m})(dx_{m:n}),  
	\end{aligned}
	\\
	\label{3.9}
	\tilde{S}_{\theta,\boldsymbol y}^{m:m}(x,B)
	=
	0, 
	&\;\;\;\;\; 
	\begin{aligned}[t]
	\tilde{S}_{\theta,\boldsymbol y}^{m:n}(x,B)
	=&
	\int_{{\cal X}^{n-m}\times B}
	\tilde{s}_{\theta,\boldsymbol y}^{m:n}(x_{m:n} ) 
	(\delta_{x}\times\mu^{n-m})(dx_{m:n}),   
	\end{aligned}
\end{align}
where $(\delta_{x}\times\mu^{n-m})(dx_{m:n})=\delta_{x}(dx_{m})\mu(dx_{m+1})\cdots\mu(dx_{n})$. 
$\tilde{F}_{\theta,\boldsymbol y}^{m:n}(\xi)$, 
$\tilde{G}_{\theta,\boldsymbol y}^{m:n}(\xi,\zeta)$ and 
$\tilde{H}_{\theta,\boldsymbol y}^{m:n}(\xi,\zeta)$
are the functions mapping $\xi\in{\cal P}({\cal X})$, $\zeta\in{\cal M}_{s}^{d}({\cal X})$ to 
${\cal P}({\cal X})$, ${\cal M}_{s}^{d}({\cal X})$, ${\cal M}_{s}^{d}({\cal X})$ (respectively) 
defined by 
\begin{align}\label{3.21}
	\tilde{F}_{\theta,\boldsymbol y}^{m:m}(\xi)(B)
	=
	\xi(B), 
	&\;\;\; 
	\tilde{F}_{\theta,\boldsymbol y}^{m:n}(\xi)(B)
	=
	\frac{\big(\xi\tilde{R}_{\theta,\boldsymbol y}^{m:n} \big)(B)}
	{\big\langle\xi\tilde{R}_{\theta,\boldsymbol y}^{m:n} \big\rangle }, 
	\\
	\label{3.23}
	\tilde{H}_{\theta,\boldsymbol y}^{m:m}(\xi,\zeta)(B)
	\!=\!
	\zeta(B), 
	&\;\;\; 
	\tilde{H}_{\theta,\boldsymbol y}^{m:n}(\xi,\zeta)(B)
	\!=\!
	\frac{\big(\zeta\tilde{R}_{\theta,\boldsymbol y}^{m:n} \big)(B)
	\!+\!
	\big(\xi\tilde{S}_{\theta,\boldsymbol y}^{m:n} \big)(B)}
	{\big\langle\xi\tilde{R}_{\theta,\boldsymbol y}^{m:n} \big\rangle}, 
	\\
	\label{3.25}
	\tilde{G}_{\theta,\boldsymbol y}^{m:m}(\xi,\zeta)(B)
	\!=\!
	\zeta(B), 
	&\;\;\; 
	\tilde{G}_{\theta,\boldsymbol y}^{m:n}(\xi,\zeta)(B)
	\!=\!
	\tilde{H}_{\theta,\boldsymbol y}^{m:n}(\xi,\zeta)(B)
	\!-\!
	\tilde{F}_{\theta,\boldsymbol y}^{m:n}(\xi)(B)
	\big\langle \tilde{H}_{\theta,\boldsymbol y}^{m:n}(\xi,\zeta) \big\rangle.  
\end{align}
$\tilde{\gamma}_{\theta,\boldsymbol y}^{n}(\xi)$ and 
$\tilde{\delta}_{\theta,\boldsymbol y}^{n}(\xi,\zeta)$
are the functions mapping $\xi\in{\cal P}({\cal X})$, $\zeta\in{\cal M}_{s}^{d}({\cal X})$ to
${\cal M}_{s}^{d}({\cal X})$ defined by 
\begin{align}\label{3.27}
	\tilde{\gamma}_{\theta,\boldsymbol y}^{n}(\xi)(B)
	=
	\frac{\int_{B} \nabla_{\theta} q_{\theta}(y_{n}|x) \xi(dx)  }
	{\int q_{\theta}(y_{n}|x) \xi(dx) }, 
	\;\;\;\; 
	\tilde{\delta}_{\theta,\boldsymbol y}^{n}(\xi,\zeta)(B)
	=
	\frac{\int_{B} q_{\theta}(y_{n}|x) \zeta(dx) }
	{\int q_{\theta}(y_{n}|x) \xi(dx) }
\end{align}
for $n\geq 0$. 
$\tilde{\alpha}_{\theta,\boldsymbol y}^{n}(\xi)$ and 
$\tilde{\beta}_{\theta,\boldsymbol y}^{n}(\xi,\zeta)$
are the functions mapping $\xi\in{\cal P}({\cal X})$, $\zeta\in{\cal M}_{s}^{d}({\cal X})$ to 
${\cal P}({\cal X})$, ${\cal M}_{s}^{d}({\cal X})$ (respectively) 
defined by 
\begin{align}\label{3.29}
	\tilde{\alpha}_{\theta,\boldsymbol y}^{n}(\xi)(B) 
	=
	\frac{\int_{B} q_{\theta}(y_{n}|x) \xi(dx) }{\int q_{\theta}(y_{n}|x) \xi(dx) }, 
	\;\;\;\;\; 
	\tilde{\beta}_{\theta,\boldsymbol y}^{n}(\xi,\zeta)(B) 
	=
	\tilde{\gamma}_{\theta,\boldsymbol y}^{n}(\xi)(B)
	+
	\tilde{\delta}_{\theta,\boldsymbol y}^{n}(\xi,\zeta)(B). 
\end{align}

\begin{vremark}
Let Assumptions 2.1 -- 2.3 hold. 
Moreover, let $\xi_{\theta}(dx)$, $\zeta_{\theta}(dx)$ and $w_{\theta}(x)$ be of the form 
(2.6), (6.1). 
Then, it is straightforward to verify 
\begin{align*}
	&
	\tilde{F}_{\theta,\boldsymbol y}^{0:n}(\xi_{\theta} )(B)
	=
	P\left(\left.
	X_{n}^{\theta}\in B
	\right|
	Y_{1:n}^{\theta} = y_{1:n} 
	\right), 
	\;\;\;\;\; 
	\tilde{G}_{\theta,\boldsymbol y}^{0:n}(\xi_{\theta},\zeta_{\theta} )(B)
	=
	\nabla_{\theta} 
	\tilde{F}_{\theta,\boldsymbol y}^{0:n}(\xi_{\theta} )(B). 
\end{align*}
Hence, $\tilde{F}_{\theta,\boldsymbol y}^{m:n}(\xi)$ and 
$\tilde{G}_{\theta,\boldsymbol y}^{m:n}(\xi,\zeta)$ 
can be considered as (a generalization of) 
the optimal filter and its gradient. 
\end{vremark}

\begin{lemma}\label{lemma3.1}
Let $\theta$, $B$, $\xi$, $\zeta$ be any elements of 
$\Theta$, ${\cal B}({\cal X})$, ${\cal P}({\cal X})$, ${\cal M}_{s}^{d}({\cal X})$ (respectively), 
while $\boldsymbol y = \{y_{n} \}_{n\geq 0}$ is any sequence in ${\cal Y}$. 
Moreover, let $n$, $m$ be any integers satisfying $n\geq m\geq 0$. 

(i) Suppose that Assumption 2.1 holds. Then, we 
have 
\begin{align}
	\label{l3.1.3*}
	F_{\theta,\boldsymbol y}^{m:n+1}(\xi)(B)
	=&
	\left(\tilde{F}_{\theta,\boldsymbol y}^{m:n}\big(\tilde{\alpha}_{\theta,\boldsymbol y}^{m}(\xi) \big)
	\tilde{A}_{\theta}
	\right)(B). 
\end{align}

(ii) Suppose that Assumptions 2.1 and 2.2 hold. Then, we 
have 
\begin{align}
	\label{l3.1.7*}
	G_{\theta,\boldsymbol y}^{m:n+1}(\xi,\zeta)(B)
	=&
	\Big(
	\tilde{F}_{\theta,\boldsymbol y}^{m:n}
	\big(\tilde{\alpha}_{\theta,\boldsymbol y}^{m}(\xi) \big)
	\tilde{B}_{\theta} 
	\Big)(B)
	+
	\left( 
	\tilde{G}_{\theta,\boldsymbol y}^{m:n}
	\big(\tilde{\alpha}_{\theta,\boldsymbol y}^{m}(\xi), \tilde{\beta}_{\theta,\boldsymbol y}^{m}(\xi,\zeta) \big)
	\tilde{A}_{\theta} 
	\right)(B).  
\end{align}
\end{lemma}

\begin{proof}
(i) Throughout this part of the proof, $x$ denotes any element of ${\cal X}$. 
Using (2.2), (\ref{3.3}), it is straightforward to verify 
\begin{align}\label{l3.1.1}
	r_{\theta,\boldsymbol y}^{m:n+1}(x_{m:n+1} )
	=
	p_{\theta}(x_{n+1}|x_{n} ) q_{\theta}(y_{m}|x_{m} ) 
	\tilde{r}_{\theta,\boldsymbol y}^{m:n}(x_{m:n} )
\end{align}
for $x_{m},\dots,x_{n+1}\in{\cal X}$. 
Therefore, (4.2), (\ref{3.7}) imply  
\begin{align}\label{l3.1.3}
	R_{\theta,\boldsymbol y}^{m:n+1}(x,B)
	=&
	\begin{aligned}[t]
	&\int\left(
	\int_{B} p_{\theta}(x_{n+1}|x_{n} ) \mu(dx_{n+1} ) \right) 
	q_{\theta}(y_{m}|x_{m} ) 
	\\
	&\cdot 
	\tilde{r}_{\theta,\boldsymbol y}^{m:n}(x_{m:n} )
	(\delta_{x}\times\mu^{n-m} )(dx_{m:n} )
	\end{aligned}
	\\ \nonumber
	=&
	\big(\tilde{R}_{\theta,\boldsymbol y}^{m:n}\tilde{A}_{\theta} \big)(x,B)
	q_{\theta}(y_{m}|x).  
\end{align}
Then, using (\ref{3.29}), we conclude  
\begin{align}\label{l3.1.5}
	\frac{\big( \xi R_{\theta,\boldsymbol y}^{m:n+1} \big)(B) }
	{\int q_{\theta}(y_{m}|x)\xi(dx) }
	=&
	\frac{\int\big( \tilde{R}_{\theta,\boldsymbol y}^{m:n}\tilde{A}_{\theta} \big)(x,B)
	q_{\theta}(y_{m}|x)\xi(dx) }
	{\int q_{\theta}(y_{m}|x)\xi(dx) }
	=
	\big( \tilde{\alpha}_{\theta,\boldsymbol y}^{m}(\xi)
	\tilde{R}_{\theta,\boldsymbol y}^{m:n}\tilde{A}_{\theta} \big)(B). 
\end{align}
Since $\tilde{A}_{\theta}(x,{\cal X})=1$, we have  
\begin{align}\label{l3.1.5'}
	\frac{\big\langle \xi R_{\theta,\boldsymbol y}^{m:n+1} \big\rangle }
	{\int q_{\theta}(y_{m}|x)\xi(dx) }
	=
	\big\langle \tilde{\alpha}_{\theta,\boldsymbol y}^{m}(\xi)
	\tilde{R}_{\theta,\boldsymbol y}^{m:n} \big\rangle. 
\end{align}
Combining this with (4.4), (\ref{3.21}), we get  
\begin{align*}
	F_{\theta,\boldsymbol y}^{m:n+1}(\xi)(B) 
	=&
	\frac{\big( \xi R_{\theta,\boldsymbol y}^{m:n+1} \big)(B) }
	{\big\langle \xi R_{\theta,\boldsymbol y}^{m:n+1} \big\rangle }
	=
	\frac{\big(\tilde{\alpha}_{\theta,\boldsymbol y}^{m}(\xi)
	\tilde{R}_{\theta,\boldsymbol y}^{m:n} \tilde{A}_{\theta} \big)(B) }
	{\big\langle\tilde{\alpha}_{\theta,\boldsymbol y}^{m}(\xi)
	\tilde{R}_{\theta,\boldsymbol y}^{m:n} \big\rangle}
	=
	\left( \tilde{F}_{\theta,\boldsymbol y}^{m:n}\big(\tilde{\alpha}_{\theta,\boldsymbol y}^{m}(\xi) \big)
	\tilde{A}_{\theta}
	\right)(B). 
\end{align*}

(ii) Let $x$ have the same meaning as in (i). 
Differentiating (\ref{l3.1.1}) in $\theta$ and using (4.1), (\ref{3.5}), 
we get 
\begin{align*}
	s_{\theta,\boldsymbol y}^{m:n+1}(x_{m:n+1} )
	=&
	p_{\theta}(x_{n+1}|x_{n} )\nabla_{\theta}q_{\theta}(y_{m}|x_{m} ) 
	\tilde{r}_{\theta,\boldsymbol y}^{m:n}(x_{m:n} )
	\\
	&+
	\nabla_{\theta}p_{\theta}(x_{n+1}|x_{n} )q_{\theta}(y_{m}|x_{m} ) \big)
	\tilde{r}_{\theta,\boldsymbol y}^{m:n}(x_{m:n} )
	\\
	&+
	p_{\theta}(x_{n+1}|x_{n} )q_{\theta}(y_{m}|x_{m} ) 
	\tilde{s}_{\theta,\boldsymbol y}^{m:n}(x_{m:n} )
\end{align*}
for $x_{m},\dots,x_{n+1}\in{\cal X}$. 
Therefore, (4.3), (\ref{3.7}), (\ref{3.9}) imply  
\begin{align*}
	S_{\theta,\boldsymbol y}^{m:n+1}(x,B)
	=&
	\begin{aligned}[t]
	\int&\left(\int_{B}p_{\theta}(x_{n+1}|x_{n})\mu(dx_{n+1} ) \right)
	\nabla_{\theta}q_{\theta}(y_{m}|x_{m} ) 
	\\
	&\cdot 
	\tilde{r}_{\theta,\boldsymbol y}^{m:n}(x_{m:n} )
	(\delta_{x}\times\mu^{n-m})(dx_{m:n} )
	\end{aligned}
	\\
	&+
	\begin{aligned}[t]
	\int&\left(\int_{B}\nabla_{\theta}p_{\theta}(x_{n+1}|x_{n})\mu(dx_{n+1} ) \right)
	q_{\theta}(y_{m}|x_{m} ) 
	\\
	&\cdot 
	\tilde{r}_{\theta,\boldsymbol y}^{m:n}(x_{m:n} )
	(\delta_{x}\times\mu^{n-m})(dx_{m:n} )
	\end{aligned}
	\\
	&+
	\begin{aligned}[t]
	\int&\left(\int_{B}p_{\theta}(x_{n+1}|x_{n})\mu(dx_{n+1} ) \right)
	q_{\theta}(y_{m}|x_{m} ) 
	\\
	&\cdot 
	\tilde{s}_{\theta,\boldsymbol y}^{m:n}(x_{m:n} )
	(\delta_{x}\times\mu^{n-m})(dx_{m:n} )
	\end{aligned}
	\\
	=&
	\big(\tilde{R}_{\theta,\boldsymbol y}^{m:n}\tilde{A}_{\theta} \big)(x,B)
	\nabla_{\theta}q_{\theta}(y_{m}|x)
	+
	\big(\tilde{R}_{\theta,\boldsymbol y}^{m:n}\tilde{B}_{\theta} \big)(x,B)
	q_{\theta}(y_{m}|x)
	\\
	&+
	\big(\tilde{S}_{\theta,\boldsymbol y}^{m:n}\tilde{A}_{\theta} \big)(x,B)
	q_{\theta}(y_{m}|x).   
\end{align*}
Then, using (\ref{3.27}), (\ref{3.29}), 
we conclude 
\begin{align}\label{l3.1.21'}
	\frac{\big(\xi S_{\theta,\boldsymbol y}^{m:n+1} \big)(B) }
	{\int q_{\theta}(y_{m}|x) \xi(dx) }
	=&
	\frac{\int \big(\tilde{R}_{\theta,\boldsymbol y}^{m:n}\tilde{A}_{\theta} \big)(x,B)
	\nabla_{\theta}q_{\theta}(y_{m}|x) \xi(dx)}
	{\int q_{\theta}(y_{m}|x)\xi(dx) }
	\nonumber\\
	&
	+
	\frac{\int \big(\tilde{R}_{\theta,\boldsymbol y}^{m:n}\tilde{B}_{\theta} \big)(x,B)
	q_{\theta}(y_{m}|x) \xi(dx)}
	{\int q_{\theta}(y_{m}|x)\xi(dx) }
	\nonumber\\
	&+
	\frac{\int \big(\tilde{S}_{\theta,\boldsymbol y}^{m:n}\tilde{A}_{\theta} \big)(x,B)
	q_{\theta}(y_{m}|x) \xi(dx)}
	{\int q_{\theta}(y_{m}|x)\xi(dx) }
	\nonumber\\
	=&
	\big(\tilde{\gamma}_{\theta,\boldsymbol y}^{m}(\xi)
	\tilde{R}_{\theta,\boldsymbol y}^{m:n}\tilde{A}_{\theta} \big)(B)
	+
	\big(\tilde{\alpha}_{\theta,\boldsymbol y}^{m}(\xi)
	\tilde{R}_{\theta,\boldsymbol y}^{m:n}\tilde{B}_{\theta} \big)(B)
	\nonumber\\
	&
	+
	\big(\tilde{\alpha}_{\theta,\boldsymbol y}^{m}(\xi)
	\tilde{S}_{\theta,\boldsymbol y}^{m:n}\tilde{A}_{\theta} \big)(B). 
\end{align}
Similarly, relying on (\ref{l3.1.3}), we deduce 
\begin{align*}
	\frac{\big(\zeta R_{\theta,\boldsymbol y}^{m:n+1} \big)(B) }
	{\int q_{\theta}(y_{m}|x) \xi(dx) }
	=&
	\frac{\int \big(\tilde{R}_{\theta,\boldsymbol y}^{m:n}\tilde{A}_{\theta} \big)(x,B)
	q_{\theta}(y_{m}|x) \zeta(dx)}
	{\int q_{\theta}(y_{m}|x)\xi(dx) }
	=
	\big(\tilde{\delta}_{\theta,\boldsymbol y}^{m}(\xi,\zeta)
	\tilde{R}_{\theta,\boldsymbol y}^{m:n}\tilde{A}_{\theta} \big)(B). 
\end{align*}
Consequently, (\ref{3.27}), (\ref{3.29}), (\ref{l3.1.21'}) imply 
\begin{align*}
	\frac{\big(\zeta R_{\theta,\boldsymbol y}^{m:n+1} \big)(B) 
	+ \big(\xi S_{\theta,\boldsymbol y}^{m:n+1} \big)(B) }
	{\int q_{\theta}(y_{m}|x) \xi(dx) }
	=&
	\big(\tilde{\alpha}_{\theta,\boldsymbol y}^{m}(\xi)
	\tilde{R}_{\theta,\boldsymbol y}^{m:n}\tilde{B}_{\theta} \big)(B)
	+
	\big(\tilde{\beta}_{\theta,\boldsymbol y}^{m}(\xi,\zeta)
	\tilde{R}_{\theta,\boldsymbol y}^{m:n}\tilde{A}_{\theta} \big)(B)
	\\
	&+
	\big(\tilde{\alpha}_{\theta,\boldsymbol y}^{m}(\xi)
	\tilde{S}_{\theta,\boldsymbol y}^{m:n}\tilde{A}_{\theta} \big)(B). 
\end{align*}
Combining this with (4.5), (\ref{3.21}), (\ref{3.23}), (\ref{l3.1.5'}), 
we get 
\begin{align}\label{l3.1.21}
	H_{\theta,\boldsymbol y}^{m:n+1}(\xi,\zeta)(B) 
	=&
	\frac{\big(\zeta R_{\theta,\boldsymbol y}^{m:n+1} \big)(B) 
	+
	\big(\xi S_{\theta,\boldsymbol y}^{m:n+1} \big)(B) }
	{\big\langle \xi R_{\theta,\boldsymbol y}^{m:n+1} \big\rangle }
	\\ \nonumber
	=&
	\frac{\big(\tilde{\beta}_{\theta,\boldsymbol y}^{m}(\xi,\zeta)
	\tilde{R}_{\theta,\boldsymbol y}^{m:n}\tilde{A}_{\theta} \big)(B)
	+ 
	\big(\tilde{\alpha}_{\theta,\boldsymbol y}^{m}(\xi)
	\tilde{S}_{\theta,\boldsymbol y}^{m:n}\tilde{A}_{\theta} \big)(B) }
	{\big\langle \tilde{\alpha}_{\theta,\boldsymbol y}^{m}(\xi)
	\tilde{R}_{\theta,\boldsymbol y}^{m:n} \big\rangle }
	\\ \nonumber
	&+
	\frac{\big(\tilde{\alpha}_{\theta,\boldsymbol y}^{m}(\xi)
	\tilde{R}_{\theta,\boldsymbol y}^{m:n}\tilde{B}_{\theta} \big)(B) }
	{\big\langle \tilde{\alpha}_{\theta,\boldsymbol y}^{m}(\xi)
	\tilde{R}_{\theta,\boldsymbol y}^{m:n} \big\rangle }
	\\ \nonumber
	=&
	\left(\tilde{F}_{\theta,\boldsymbol y}^{m:n}\big(\tilde{\alpha}_{\theta,\boldsymbol y}^{m}(\xi) \big)
	\tilde{B}_{\theta} \right)(B)
	+ 
	\left(\tilde{H}_{\theta,\boldsymbol y}^{m:n}\big(\tilde{\alpha}_{\theta,\boldsymbol y}^{m}(\xi), 
	\tilde{\beta}_{\theta,\boldsymbol y}^{m}(\xi,\zeta) \big)
	\tilde{A}_{\theta} \right)(B). 
\end{align}
Since $\tilde{A}_{\theta}(x,{\cal X})=1$ and  
$\tilde{B}_{\theta}(x,{\cal X})=\nabla_{\theta}\tilde{A}_{\theta}(x,{\cal X})=0$,  
we then have  
\begin{align}\label{l3.1.23}
	\big\langle H_{\theta,\boldsymbol y}^{m:n+1}(\xi,\zeta) \big\rangle 
	=
	\big\langle 
	\tilde{H}_{\theta,\boldsymbol y}^{m:n}
	\big(\tilde{\alpha}_{\theta,\boldsymbol y}^{m}(\xi), 
	\tilde{\beta}_{\theta,\boldsymbol y}^{m}(\xi,\zeta) \big)
	\big\rangle. 
\end{align}
Moreover, using (4.6), (\ref{3.25}), (\ref{l3.1.21}), (\ref{l3.1.23}), we conclude  
\begin{align*}
	G_{\theta,\boldsymbol y}^{m:n+1}(\xi,\zeta)(B)
	=&
	H_{\theta,\boldsymbol y}^{m:n+1}(\xi,\zeta)(B)
	-
	F_{\theta,\boldsymbol y}^{m:n+1}(\xi)(B)
	\big\langle H_{\theta,\boldsymbol y}^{m:n+1}(\xi,\zeta) \big\rangle
	\\
	=&
	\left(\tilde{F}_{\theta,\boldsymbol y}^{m:n}\big(\tilde{\alpha}_{\theta,\boldsymbol y}^{m}(\xi) \big)
	\tilde{B}_{\theta} \right)(B)
	+ 
	\left(\tilde{H}_{\theta,\boldsymbol y}^{m:n}\big(\tilde{\alpha}_{\theta,\boldsymbol y}^{m}(\xi), 
	\tilde{\beta}_{\theta,\boldsymbol y}^{m}(\xi,\zeta) \big)
	\tilde{A}_{\theta} \right)(B)
	\\
	&-
	F_{\theta,\boldsymbol y}^{m:n+1}(\xi)(B)
	\big\langle 
	\tilde{H}_{\theta,\boldsymbol y}^{m:n}
	\big(\tilde{\alpha}_{\theta,\boldsymbol y}^{m}(\xi), 
	\tilde{\beta}_{\theta,\boldsymbol y}^{m}(\xi,\zeta) \big)
	\big\rangle
	\\
	=&
	\left(\tilde{F}_{\theta,\boldsymbol y}^{m:n}\big(\tilde{\alpha}_{\theta,\boldsymbol y}^{m}(\xi) \big)
	\tilde{B}_{\theta} \right)(B)
	+ 
	\left(\tilde{G}_{\theta,\boldsymbol y}^{m:n}\big(\tilde{\alpha}_{\theta,\boldsymbol y}^{m}(\xi), 
	\tilde{\beta}_{\theta,\boldsymbol y}^{m}(\xi,\zeta) \big)
	\tilde{A}_{\theta} \right)(B). 
\end{align*}
\end{proof}

\begin{proof}[Proof of Proposition 4.1]
(i) Using \cite[Theorem 4.1]{{legland&oudjane}} (or \cite[Theorem 3.1]{tadic&doucet1}), we conclude that 
there exist real numbers $\rho_{1}\in(0,1)$, $\tilde{C}_{1}\in[1,\infty)$
(independent of $\theta$, $\boldsymbol y$, $\xi$, $\xi'$, $n$, $m$
and depending only on $\varepsilon$) such that 
\begin{align*}
	\left\| \tilde{F}_{\theta,\boldsymbol y}^{m:n}(\xi)
	-
	\tilde{F}_{\theta,\boldsymbol y}^{m:n}(\xi') \right\|
	\leq 
	\tilde{C}_{1}\rho_{1}^{n-m}\|\xi-\xi'\|. 
\end{align*}
As $\|\xi-\xi'\|\leq 1$, we have 
\begin{align}\label{l3.2.501}
	\left\| \tilde{F}_{\theta,\boldsymbol y}^{m:n}(\xi)
	-
	\tilde{F}_{\theta,\boldsymbol y}^{m:n}(\xi') \right\|
	\leq 
	\tilde{C}_{1}\rho_{1}^{n-m}. 
\end{align}
Throughout this part of the proof, the following notation is used. 
$x$, $x'$ are any elements of ${\cal X}$, while  
$B$ is any element of ${\cal B}({\cal X} )$. 
$C_{1}$ is the real number defined by 
$C_{1}=\tilde{C}_{1}\varepsilon^{-4}\rho_{1}^{-1}$
($\varepsilon$ is specified in Assumption 2.1). 

When, $m=n$, the left-hand sides of inequalities (4.8) reduce to 
$\|\xi-\xi'\|$ and $1$ (respectively). 
Hence, it is sufficient to show (4.8) for $n>m\geq 0$. 
In what follows in this part of the proof, 
$n$, $m$ are any integers satisfying $n>m\geq 0$. 

Since $0\leq \tilde{A}_{\theta}(x,B)\leq 1$, we have 
\begin{align*}
	\left|(\xi\tilde{A}_{\theta} )(B) - (\xi'\tilde{A}_{\theta} )(B) \right|
	\leq 
	\int \tilde{A}_{\theta}(x,B)|\xi-\xi'|(dx)
	\leq 
	\|\xi-\xi'\|. 
\end{align*}
Hence, we get 
\begin{align}\label{l3.2.501'}
	\left\|\xi\tilde{A}_{\theta}-\xi'\tilde{A}_{\theta} \right\|
	\leq 
	\|\xi-\xi'\|. 
\end{align}
Moreover, relying on (\ref{l3.2.501}), we conclude  
\begin{align}\label{l3.2.3}
	\left\| 
	\tilde{F}_{\theta,\boldsymbol y}^{m:n-1}
	\big(\tilde{\alpha}_{\theta,\boldsymbol y}^{m}(\xi) \big)
	-
	\tilde{F}_{\theta,\boldsymbol y}^{m:n-1}
	\big(\tilde{\alpha}_{\theta,\boldsymbol y}^{m}(\xi') \big)
	\right\|
	\leq &
	\tilde{C}_{1}\rho_{1}^{n-m-1} 
	\leq 
	C_{1}\rho_{1}^{n-m}. 
\end{align}
Then, using Lemma \ref{lemma3.1}, we deduce  
\begin{align*}
	\left\|
	F_{\theta,\boldsymbol y}^{m:n}(\xi)
	-
	F_{\theta,\boldsymbol y}^{m:n}(\xi')
	\right\|
	=&
	\left\|
	\tilde{F}_{\theta,\boldsymbol y}^{m:n-1}
	\big(\tilde{\alpha}_{\theta,\boldsymbol y}^{m}(\xi) \big) \tilde{A}_{\theta} 
	-
	\tilde{F}_{\theta,\boldsymbol y}^{m:n-1}
	\big(\tilde{\alpha}_{\theta,\boldsymbol y}^{m}(\xi') \big) \tilde{A}_{\theta} 
	\right\|
	\\ \nonumber
	&\leq  
	\left\|
	\tilde{F}_{\theta,\boldsymbol y}^{m:n-1}
	\big(\tilde{\alpha}_{\theta,\boldsymbol y}^{m}(\xi) \big)  
	-
	\tilde{F}_{\theta,\boldsymbol y}^{m:n-1}
	\big(\tilde{\alpha}_{\theta,\boldsymbol y}^{m}(\xi') \big)  
	\right\|
	\leq 
	C_{1}\rho_{1}^{n-m}.  
\end{align*}
Thus, the first part of (4.8) holds.

Owing to Assumption 2.1 and (2.1), we have 
\begin{align}\label{l3.2.7}
	\varepsilon^{2} \leq r_{\theta,\boldsymbol y}^{n}(x'|x) \leq \frac{1}{\varepsilon^{2} }. 
\end{align}
for $n\geq 1$. 
Moreover, due to (2.2), we have  
\begin{align*}
	r_{\theta,\boldsymbol y}^{m:n}(x_{m:n} )
	=
	r_{\theta,\boldsymbol y}^{m+1:n}(x_{m+1:n} ) r_{\theta,\boldsymbol y}^{m+1}(x_{m+1}|x_{m} ). 
\end{align*}
Hence, we get 
\begin{align}\label{l3.2.303}
	\varepsilon^{2} r_{\theta,\boldsymbol y}^{m+1:n}(x_{m+1:n} )
	\leq 
	r_{\theta,\boldsymbol y}^{m:n}(x_{m:n} )
	\leq 
	\frac{1}{\varepsilon^{2} } r_{\theta,\boldsymbol y}^{m+1:n}(x_{m+1:n} ). 
\end{align}
Consequently, (4.2) implies 
\begin{align}
	&\label{l3.2.307} 
	\big\langle \xi R_{\theta,\boldsymbol y}^{m:n} \big\rangle  
	\geq 
	\varepsilon^{2} \langle\xi\rangle  
	\int r_{\theta,\boldsymbol y}^{m+1:n}(x_{m:n} ) 
	\mu^{n-m}(dx_{m+1:n} ), 
	\\ \nonumber
	&
	\big\langle \xi R_{\theta,\boldsymbol y}^{m:n} \big\rangle  
	\leq 
	\frac{\langle\xi\rangle}{\varepsilon^{2} }
	\int r_{\theta,\boldsymbol y}^{m+1:n}(x_{m+1:n} ) 
	\mu^{n-m}(dx_{m+1:n} ). 
\end{align}
As $\langle\xi\rangle=1$, we have 
\begin{align*}
	\frac{\big\langle \xi R_{\theta,\boldsymbol y}^{m:n} \big\rangle }
	{\big\langle \xi' R_{\theta,\boldsymbol y}^{m:n} \big\rangle }
	\leq 
	\frac{1}{\varepsilon^{4}}
	\leq 
	C_{1}. 
\end{align*}

(ii) Relying on \cite[Theorem 3.2]{tadic&doucet1}
(or \cite[Theorem 2.2]{tadic&doucet3}), we deduce that
there exist real numbers $\rho_{2}\in[\rho_{1},1)$, $\tilde{C}_{2}\in[C_{1},\infty)$
(independent of $\theta$, $\boldsymbol y$, $\xi$, $\xi'$, $\zeta$, $\zeta'$, $n$, $m$
and depending only on $\varepsilon$, $d$, $K$) such that 
\begin{align*}
	\left\| \tilde{G}_{\theta,\boldsymbol y}^{m:n}(\xi,\zeta)
	-
	\tilde{G}_{\theta,\boldsymbol y}^{m:n}(\xi',\zeta') \right\|
	\leq &
	\tilde{C}_{2}\rho_{2}^{n-m}\|\xi-\xi'\|(1 + \|\zeta\| + \|\zeta'\| )
	\\ 
	&+
	\tilde{C}_{2}\rho_{2}^{n-m}\|\zeta-\zeta'\|. 
\end{align*}
Since $\|\xi-\xi'\|\leq 1$ and $\|\zeta-\zeta'\|\leq\|\zeta\|+\|\zeta'\|$, 
we have 
\begin{align}\label{l3.2.503}
	\left\| \tilde{G}_{\theta,\boldsymbol y}^{m:n}(\xi,\zeta)
	-
	\tilde{G}_{\theta,\boldsymbol y}^{m:n}(\xi',\zeta') \right\|
	\leq &
	2\tilde{C}_{2}\rho_{2}^{n-m}(1 + \|\zeta\| + \|\zeta'\| ). 
\end{align}
Throughout the rest of the proof, the following notation is used. 
$x$, $x'$, $B$ have the same meaning as in (i).
$\tilde{C}_{3}$, $\tilde{C}_{4}$, $\tilde{C}_{5}$, $C_{2}$ are the real numbers defined by 
$\tilde{C}_{3}=2K\varepsilon^{-2}$, 
$\tilde{C}_{4}=4d\tilde{C}_{2}\tilde{C}_{3}\rho_{2}^{-1}$, 
$\tilde{C}_{5}=C_{1}\tilde{C}_{3}\tilde{C}_{4}$, 
$C_{2}=d\tilde{C}_{3}\tilde{C}_{5}$
($\varepsilon$, $K$ are specified in Assumptions 2.1, 2.2). 

When $m=n$, the left-hand sides of (4.9), (4.10)
reduce to $\|\zeta-\zeta'\|$ and $\|\zeta\|$
(respectively). 
Hence, it is sufficient to show (4.9), (4.10) for $n>m\geq 0$. 
In the rest of the proof, $n$, $m$ are any integers satisfying 
$n>m\geq 0$. 

Owing to Assumptions 2.1, 2.2 and (\ref{3.51}), we have 
\begin{align*}
	\|\tilde{B}_{\theta}(x,B)\|
	\leq 
	\int_{B} \|\nabla_{\theta} p_{\theta}(x'|x) \| \mu(dx') 
	\leq 
	\int \left\|\frac{\nabla_{\theta} p_{\theta}(x'|x) }{p_{\theta}(x'|x) } \right\|
	p_{\theta}(x'|x) \mu(dx')
	\leq 
	\frac{K}{\varepsilon}
	\leq 
	\tilde{C}_{3}.
\end{align*}
Therefore, we get 
\begin{align*}
	\left\|(\xi\tilde{B}_{\theta} )(B) - (\xi'\tilde{B}_{\theta} )(B) \right\|
	\leq 
	\int \|\tilde{B}_{\theta}(x,B)\| |\xi-\xi'|(dx) 
	\leq 
	\tilde{C}_{3}\|\xi-\xi'\|. 
\end{align*}
Thus, we have 
\begin{align}\label{l3.2.301'}
	\left\|\xi\tilde{B}_{\theta} - \xi'\tilde{B}_{\theta} \right\| 
	\leq 
	\tilde{C}_{3}\|\xi-\xi'\|. 
\end{align}
Moreover, due to the same assumptions and (\ref{3.29}), we have 
\begin{align*}
	\left\| 
	\tilde{\beta}_{\theta,\boldsymbol y}^{m}(\xi,\zeta)(B)  
	\right\|
	\leq &
	\frac{\int_{B} q_{\theta}(y_{m}|x) |\zeta|(dx) 
	+ \int_{B} \|\nabla_{\theta} q_{\theta}(y_{m}|x)\| \xi(dx) }
	{\int q_{\theta}(y_{m}|x) \xi(dx) }
	\\
	\leq &
	\frac{K}{\varepsilon} 
	+
	\frac{\|\zeta\| }{\varepsilon^{2} }
	\leq 
	\tilde{C}_{3} (1 + \|\zeta\| ). 
\end{align*}
Hence, we get 
\begin{align*}
	\left\| \tilde{\beta}_{\theta,\boldsymbol y}^{m}(\xi,\zeta) \right\|
	\leq 
	d\tilde{C}_{3}(1 + \|\zeta\| ). 
\end{align*}
Consequently, (\ref{l3.2.503}) implies 
\begin{align}\label{l3.2.23} 
	&
	\left\| 
	\tilde{G}_{\theta,\boldsymbol y}^{m:n-1}
	\big(\tilde{\alpha}_{\theta,\boldsymbol y}^{m}(\xi), 
	\tilde{\beta}_{\theta,\boldsymbol y}^{m}(\xi,\zeta) \big)
	-
	\tilde{G}_{\theta,\boldsymbol y}^{m:n-1}
	\big(\tilde{\alpha}_{\theta,\boldsymbol y}^{m}(\xi'), 
	\tilde{\beta}_{\theta,\boldsymbol y}^{m}(\xi',\zeta') \big)
	\right\|
	\\ \nonumber
	&
	\begin{aligned}[b]
	\leq &
	2\tilde{C}_{2}\rho_{2}^{n-m-1} 
	\left(
	1 
	+ 
	\big\|\tilde{\beta}_{\theta,\boldsymbol y}^{m}(\xi,\zeta) \big\| 
	+
	\big\|\tilde{\beta}_{\theta,\boldsymbol y}^{m}(\xi',\zeta') \big\| 
	\right)
	\end{aligned} 
	\\ \nonumber
	&\leq
	4d\tilde{C}_{2}\tilde{C}_{3}\rho_{2}^{n-m-1} 
	(1 + \|\zeta\| + \|\zeta'\| ) 
	\leq
	\tilde{C}_{4} \rho_{2}^{n-m} 
	(1 + \|\zeta\| + \|\zeta'\| ). 
\end{align}
Then, Lemma \ref{lemma3.1} and (\ref{l3.2.501'}), 
(\ref{l3.2.3}), (\ref{l3.2.301'}),  yield 
\begin{align*}
	&
	\left\| 
	G_{\theta,\boldsymbol y}^{m:n}(\xi,\zeta)
	-
	G_{\theta,\boldsymbol y}^{m:n}(\xi',\zeta')
	\right\|
	\leq   
	\left\|
	\tilde{F}_{\theta,\boldsymbol y}^{m:n-1}
	\big(\tilde{\alpha}_{\theta,\boldsymbol y}^{m}(\xi) \big) \tilde{B}_{\theta} 
	-
	\tilde{F}_{\theta,\boldsymbol y}^{m:n-1}
	\big(\tilde{\alpha}_{\theta,\boldsymbol y}^{m}(\xi') \big) \tilde{B}_{\theta} 
	\right\|
	\\ \nonumber
	&\;\;\;\;\;\;\; +
	\left\|
	\tilde{G}_{\theta,\boldsymbol y}^{m:n-1}
	\big(\tilde{\alpha}_{\theta,\boldsymbol y}^{m}(\xi), 
	\tilde{\beta}_{\theta,\boldsymbol y}^{m}(\xi,\zeta) \big) \tilde{A}_{\theta} 
	-
	\tilde{G}_{\theta,\boldsymbol y}^{m:n-1}
	\big(\tilde{\alpha}_{\theta,\boldsymbol y}^{m}(\xi'), 
	\tilde{\beta}_{\theta,\boldsymbol y}^{m}(\xi',\zeta')\big) \tilde{A}_{\theta} 
	\right\|
	\\ \nonumber
	&\begin{aligned}[t]
	\leq &
	\tilde{C}_{3} 
	\left\|
	\tilde{F}_{\theta,\boldsymbol y}^{m:n-1}
	(\tilde{\alpha}_{\theta,\boldsymbol y}^{m}(\xi) )
	-
	\tilde{F}_{\theta,\boldsymbol y}^{m:n-1}(\tilde{\alpha}_{\theta,\boldsymbol y}^{m}(\xi') )
	\right\|
	\\
	&\;\;\;\;\;\;\; +
	\left\| 
	\tilde{G}_{\theta,\boldsymbol y}^{m:n-1}
	(\tilde{\alpha}_{\theta,\boldsymbol y}^{m}(\xi), 
	\tilde{\beta}_{\theta,\boldsymbol y}^{m}(\xi,\zeta) )
	-
	\tilde{G}_{\theta,\boldsymbol y}^{m:n-1}
	(\tilde{\alpha}_{\theta,\boldsymbol y}^{m}(\xi'), 
	\tilde{\beta}_{\theta,\boldsymbol y}^{m}(\xi',\zeta') )
	\right\|
	\end{aligned}
	\\
	&\leq 
	C_{1}\tilde{C}_{3}\rho_{1}^{n-m} 
	+
	\tilde{C}_{4}\rho_{2}^{n-m}(1 + \|\zeta\| + \|\zeta' \| )
	\leq 
	\tilde{C}_{5}\rho_{2}^{n-m}(1 + \|\zeta\| + \|\zeta' \| ). 
\end{align*}
Thus, we get 
\begin{align*}
	\left\| 
	G_{\theta,\boldsymbol y}^{m:n}(\xi,\zeta)
	-
	G_{\theta,\boldsymbol y}^{m:n}(\xi',\zeta')
	\right\|
	\leq &
	d\tilde{C}_{5}(1 + \|\zeta\| + \|\zeta' \| )
	\leq 
	C_{2}\rho_{2}^{n-m}(1 + \|\zeta\| + \|\zeta' \| ). 
\end{align*}

Owing to (2.2), (4.1), we have 
\begin{align*}
	s_{\theta,\boldsymbol y}^{m:n}(x_{m:n} )
	=
	r_{\theta,\boldsymbol y}^{m:n}(x_{m:n} )
	\left(
	\sum_{k=m+1}^{n} 
	\frac{\nabla_{\theta} r_{\theta,\boldsymbol y}^{k}(x_{k}|x_{k-1}) }
	{r_{\theta,\boldsymbol y}^{k}(x_{k}|x_{k-1}) }
	\right)
\end{align*}
for $x_{m},\dots,x_{n}\in{\cal X}$. 
Moreover, due to Assumptions 2.1, 2.2, we have 
\begin{align}\label{l3.2.25}
	\left\| 
	\frac{\nabla_{\theta} r_{\theta,\boldsymbol y}^{n}(x'|x) }
	{r_{\theta,\boldsymbol y}^{n}(x'|x) }
	\right\|
	\leq 
	\left\| 
	\frac{\nabla_{\theta} p_{\theta}(x'|x) }{p_{\theta}(x'|x) }
	\right\|
	+
	\left\| 
	\frac{\nabla_{\theta} q_{\theta}(y_{n-1}|x) }{q_{\theta}(y_{n-1}|x) }
	\right\|
	\leq 
	\frac{2K}{\varepsilon} 
	\leq 
	\tilde{C}_{3}.  
\end{align}
Therefore, we get 
\begin{align*}
	\left\| s_{\theta,\boldsymbol y}^{m:n}(x_{m:n} ) \right\|
	\leq &
	r_{\theta,\boldsymbol y}^{m:n}(x_{m:n} )
	\left(
	\sum_{k=m+1}^{n} 
	\left\|
	\frac{\nabla_{\theta} r_{\theta,\boldsymbol y}^{k}(x_{k}|x_{k-1}) }
	{r_{\theta,\boldsymbol y}^{k}(x_{k}|x_{k-1}) }
	\right\| 
	\right)
	\\
	\leq &
	\tilde{C}_{3}(n-m)r_{\theta,\boldsymbol y}^{m:n}(x_{m:n} ). 
\end{align*}
Then, relying on (4.3), (\ref{l3.2.303}), we deduce  
\begin{align}\label{l3.2.27}
	\left\|
	\big(\xi S_{\theta,\boldsymbol y}^{m:n} \big)(B) 
	\right\|
	\leq &
	\int 
	\left\| s_{\theta,\boldsymbol y}^{m:n}(x_{m:n} ) \right\|
	(\xi\times\mu^{n-m} )(dx_{m:n} )
	\\ \nonumber
	\leq & 
	\frac{\tilde{C}_{3}\langle\xi\rangle(n\!-\!m) }{\varepsilon^{2} }\!
	\int 
	r_{\theta,\boldsymbol y}^{m+1:n}(x_{m+1:n} ) 
	\mu^{n-m}(dx_{m+1:n} ).  
\end{align}
Similarly, using (4.2), (\ref{l3.2.303}), we conclude  
\begin{align}\label{l3.2.29}
	\left\|
	\big(\zeta R_{\theta,\boldsymbol y}^{m:n} \big)(B) 
	\right\|
	\leq &
	\int r_{\theta,\boldsymbol y}^{m:n}(x_{m:n} ) 
	(|\zeta|\times\mu^{n-m})(dx_{m:n} ) 
	\\ \nonumber
	\leq & 
	\frac{\|\zeta\| }{\varepsilon^{2}}\!
	\int 
	r_{\theta,\boldsymbol y}^{m+1:n}(x_{m+1:n} ) 
	\mu^{n-m}(dx_{m+1:n} ). 
\end{align}
Since $\langle\xi\rangle=1$, (4.5), (\ref{l3.2.307}), (\ref{l3.2.27}), (\ref{l3.2.29}) imply 
\begin{align*}
	\left\| H_{\theta,\boldsymbol y}^{m:n}(\xi,\zeta)(B)\right\|
	\leq 
	\frac{\big\|
	\big(\zeta R_{\theta,\boldsymbol y}^{m:n} \big)(B) 
	\big\|
	+
	\big\|
	\big(\xi S_{\theta,\boldsymbol y}^{m:n} \big)(B) 
	\big\|}
	{\big\langle \xi R_{\theta,\boldsymbol y}^{m:n} \big\rangle }
	\leq 
	\frac{\tilde{C}_{3} }{\varepsilon^{4} } (n-m+\|\zeta\| ). 
\end{align*} 
Hence, we have 
\begin{align*}
	\left\| H_{\theta,\boldsymbol y}^{m:n}(\xi,\zeta)\right\|
	\leq 
	\frac{d\tilde{C}_{3} }{\varepsilon^{4} }(n-m+\|\zeta\| )
	\leq 
	C_{2} (n-m+\|\zeta\| ). 
\end{align*}
\end{proof}

\section{\nopunct}

In this section, we prove Lemma 6.1. 
The proof is mainly based on the same arguments as Lemma \ref{lemma3.1}.

\begin{proof}[Proof of Lemma 6.1]
(i) 
When $m=0$, (6.14) is trivially satisfied. 
Hence, it is sufficient to show (6.14) for $n\geq m>0$. 
In this part of the proof, we assume $n\geq m>0$. 

Owing to (2.2), we have 
\begin{align}\label{l5.3.303}
	r_{\theta,\boldsymbol Y}^{m-1:n}(x_{m-1:n} )
	=
	r_{\theta,\boldsymbol Y}^{m:n}(x_{m:n} ) r_{\theta,\boldsymbol Y}^{m-1:m}(x_{m-1:m} )
\end{align}
for $x_{m-1},\dots,x_{n}\in{\cal X}$. 
Consequently, (4.2), (6.8) implies 
\begin{align*}
	\big(\hat{\xi}_{m-1}^{\theta}\hat{R}_{m-1:n}^{\theta} \big)(B)
	=&
	\begin{aligned}[t]
	\int&\left(\int_{{\cal X}^{n-m-1}\times B}r_{\theta,\boldsymbol Y}^{m:n}(x_{m:n} ) 
	\mu^{n-m}(dx_{m+1:n}) \right)
	\\
	&\cdot 
	r_{\theta,\boldsymbol Y}^{m-1:m}(x_{m-1:m} )
	(\hat{\xi}_{m-1}^{\theta}\times\mu)(dx_{m-1:m} )
	\end{aligned}
	\\
	=&
	\int \hat{R}_{m:n}^{\theta}(B) 
	\big( \hat{\xi}_{m-1}^{\theta}\hat{R}_{m-1:m}^{\theta} \big)(dx). 
\end{align*}
Combining this with (4.4), (6.5), (6.8), we get 
\begin{align}\label{l5.3.1}
	\big( \hat{F}_{m-1:m}^{\theta}\hat{R}_{m:n}^{\theta} \big) (B) 
	=&
	\frac{\int \hat{R}_{m:n}^{\theta}(B) 
	\big( \hat{\xi}_{m-1}^{\theta}\hat{R}_{m-1:m}^{\theta} \big)(dx) }
	{\big\langle \hat{\xi}_{m-1}^{\theta}\hat{R}_{m-1:m}^{\theta} \big\rangle }
	=
	\frac{\big(\hat{\xi}_{m-1}^{\theta}\hat{R}_{m-1:n}^{\theta} \big)(B) }
	{\big\langle \hat{\xi}_{m-1}^{\theta}\hat{R}_{m-1:m}^{\theta} \big\rangle }. 
\end{align}
Thus, we have 
\begin{align}\label{l5.3.1'}
	\big\langle \hat{F}_{m-1:m}^{\theta}\hat{R}_{m:n}^{\theta} \big\rangle
	=
	\frac{\big\langle \hat{\xi}_{m-1}^{\theta}\hat{R}_{m-1:n}^{\theta} \big\rangle }
	{\big\langle \hat{\xi}_{m-1}^{\theta}\hat{R}_{m-1:m}^{\theta} \big\rangle }. 
\end{align}
Then, (4.4), (6.5) yield  
\begin{align*}
	\hat{F}_{m-1:n}^{\theta}(B)
	=
	\frac{\big( \hat{\xi}_{m-1}^{\theta}\hat{R}_{m-1:n}^{\theta} \big)(B) }
	{\big\langle \hat{\xi}_{m-1}^{\theta}\hat{R}_{m-1:n}^{\theta} \big\rangle }
	=
	\frac{\big( \hat{F}_{m-1:m}^{\theta}\hat{R}_{m:n}^{\theta} \big) (B)  }
	{\big\langle \hat{F}_{m-1:m}^{\theta}\hat{R}_{m:n}^{\theta} \big\rangle }. 
\end{align*}
Hence, (6.14) holds.

(ii) Owing to (6.4), (6.9), we have 
\begin{align}\label{l5.3.301}
	\big(\xi\hat{\Psi}_{m:n}^{\theta} \big)(B) 
	=&
	\int \hat{R}_{m:n}^{\theta}(x,B) \hat{v}_{m}^{\theta}(x) \xi(dx)
	+
	\int \hat{S}_{m:n}^{\theta}(x,B) \xi(dx)
	\\ \nonumber
	=&
	\big(\hat{\alpha}_{m}^{\theta}(\xi)\hat{R}_{m:n}^{\theta} \big)(B) 
	+
	\big(\xi\hat{S}_{m:n}^{\theta} \big)(B). 
\end{align}
Combining this with (4.5), we get 
\begin{align}\label{l5.3.501'}
	\frac{\big(\xi\hat{\Psi}_{m:n}^{\theta} \big)(B) }
	{\big\langle\xi\hat{R}_{m:n}^{\theta} \big\rangle }
	=
	\frac{\big(\hat{\alpha}_{m}^{\theta}(\xi)\hat{R}_{m:n}^{\theta} \big)(B) 
	+
	\big(\xi\hat{S}_{m:n}^{\theta} \big)(B) }
	{\big\langle\xi\hat{R}_{m:n}^{\theta} \big\rangle }
	=
	H_{\theta,\boldsymbol Y}^{m:n}(\xi,\hat{\alpha}_{m}^{\theta}(\xi) )(B). 
\end{align}
Consequently, (4.6), (6.10) imply 
\begin{align*}
	\frac{\big(\xi\hat{\Phi}_{m:n}^{\theta} \big)(B) }
	{\big\langle\xi\hat{R}_{m:n}^{\theta} \big\rangle }
	=&
	\frac{\big(\xi\hat{\Psi}_{m:n}^{\theta} \big)(B) 
	- \hat{F}_{m-1:n}^{\theta}(B)\big\langle\xi\hat{\Psi}_{m:n}^{\theta} \big\rangle }
	{\big\langle\xi\hat{R}_{m:n}^{\theta} \big\rangle}
	\\
	=&
	H_{\theta,\boldsymbol Y}^{m:n}(\xi,\hat{\alpha}_{m}^{\theta}(\xi) )(B)
	-
	\hat{F}_{m-1:n}^{\theta}(B)
	\big\langle H_{\theta,\boldsymbol Y}^{m:n}(\xi,\hat{\alpha}_{m}^{\theta}(\xi) ) \big\rangle
	\\
	=&
	G_{\theta,\boldsymbol Y}^{m:n}\big(\xi,\hat{\alpha}_{m}^{\theta}(\xi) \big)(B)
	+
	\big(
	F_{\theta,\boldsymbol Y}^{m:n}(\xi)(B) 
	-
	\hat{F}_{m-1:n}^{\theta}(B)
	\big)
	\big\langle H_{\theta,\boldsymbol Y}^{m:n}(\xi,\hat{\alpha}_{m}^{\theta}(\xi) ) \big\rangle. 
\end{align*}
Hence, (6.16), (6.17) hold.

Let $\hat{\beta}_{n}^{\theta}(dx)$ be the (random) element of 
${\cal M}_{s}^{d}({\cal X})$ defined by 
\begin{align}\label{l5.3.7''}
	\hat{\beta}_{n}^{\theta}(B)
	=
	\hat{\alpha}_{n}^{\theta}\big(\hat{\xi}_{n}^{\theta} \big)(B)
\end{align}
($\hat{\alpha}_{n}^{\theta}(\xi)(dx)$ is defined in (6.4)). 
Due to (2.1), (2.8), (2.9), (5.3), (6.3), we have 
\begin{align*}
	W_{k,i}^{\theta}
	=&
	\frac{
	\sum_{j=1}^{N} 
	r_{\theta,\boldsymbol Y}^{k}\big(\hat{X}_{k,i}^{\theta} |\hat{X}_{k-1,j}^{\theta} \big) 
	W_{k-1,j}^{\theta} 
	+
	\sum_{j=1}^{N} 
	\nabla_{\theta} 
	r_{\theta,\boldsymbol Y}^{k}\big(\hat{X}_{k,i}^{\theta} |\hat{X}_{k-1,j}^{\theta} \big) }
	{\sum_{j=1}^{N} 
	r_{\theta,\boldsymbol Y}^{k}\big(\hat{X}_{k,i}^{\theta} |\hat{X}_{k-1,j}^{\theta} \big) }
	\\
	=&
	\frac{
	\sum_{j=1}^{N} 
	r_{\theta,\boldsymbol Y}^{k}\big(\hat{X}_{k,i}^{\theta} |\hat{X}_{k-1,j}^{\theta} \big) 
	V_{k-1,j}^{\theta} 
	+
	\sum_{j=1}^{N} 
	\nabla_{\theta} 
	r_{\theta,\boldsymbol Y}^{k}\big(\hat{X}_{k,i}^{\theta} |\hat{X}_{k-1,j}^{\theta} \big) }
	{\sum_{j=1}^{N} 
	r_{\theta,\boldsymbol Y}^{k}\big(\hat{X}_{k,i}^{\theta} |\hat{X}_{k-1,j}^{\theta} \big) }
	+
	\frac{1}{N} \!\sum_{j=1}^{N} W_{k-1,j}^{\theta}
	\\
	=&
	\frac{
	\int r_{\theta,\boldsymbol Y}^{k}\big(\hat{X}_{k,i}^{\theta} |x \big) 
	\hat{\zeta}_{k-1}^{\theta}(dx) 
	+
	\int 
	\nabla_{\theta} 
	r_{\theta,\boldsymbol Y}^{k}\big(\hat{X}_{k,i}^{\theta} |x \big) 
	\hat{\xi}_{k-1}^{\theta}(dx) }
	{\int r_{\theta,\boldsymbol Y}^{k}\big(\hat{X}_{k,i}^{\theta} |x \big) 
	\hat{\xi}_{k-1}^{\theta}(dx) }
	+
	\frac{1}{N} \!\sum_{j=1}^{N} W_{k-1,j}^{\theta}
	\\
	=&
	\hat{v}_{k}^{\theta}\big(\hat{X}_{k,i}^{\theta} \big)
	+
	\frac{1}{N} \!\sum_{j=1}^{N} W_{k-1,j}^{\theta}
\end{align*}
for $1\leq i\leq N$, $k\geq 1$. 
Hence, (2.8), (5.3), (6.4), (\ref{l5.3.7''}) imply  
\begin{align}\label{l5.3.7'}
	V_{k,i}^{\theta}
	=
	\hat{v}_{k}^{\theta}\big(\hat{X}_{k,i}^{\theta} \big)
	-
	\frac{1}{N} \sum_{j=1}^{N} 
	\hat{v}_{k}^{\theta}\big(\hat{X}_{k,j}^{\theta} \big)
	=&
	\hat{v}_{k}^{\theta}\big(\hat{X}_{k,i}^{\theta} \big)
	-
	\int \hat{v}_{k}^{\theta}(x)
	\hat{\xi}_{k}^{\theta}(dx) 
	\\ \nonumber
	=&
	\hat{v}_{k}^{\theta}\big(\hat{X}_{k,i}^{\theta} \big)
	-
	\big\langle\hat{\beta}_{k}^{\theta}\big\rangle.
\end{align}
Since $W_{0,i}^{\theta}=w_{\theta}(\hat{X}_{0,i}^{\theta} )$, 
(2.8), (5.3), (6.1), (6.3), (6.4), (\ref{l5.3.7''}) yield  
\begin{align}\label{l5.3.9'}
	V_{0,i}^{\theta}
	=
	u_{\theta}\big( \hat{X}_{0,i}^{\theta} \big) 
	-
	\frac{1}{N} \sum_{j=1}^{N} 
	u_{\theta}\big( \hat{X}_{0,j}^{\theta} \big) 
	=&
	\hat{v}_{0}^{\theta}\big(\hat{X}_{0,i}^{\theta} \big)
	-
	\frac{1}{N} \sum_{j=1}^{N} 
	\hat{v}_{0}^{\theta}\big(\hat{X}_{0,j}^{\theta} \big)
	\\ \nonumber
	=&
	\hat{v}_{0}^{\theta}\big(\hat{X}_{0,i}^{\theta} \big)
	-
	\int \hat{v}_{0}^{\theta}(x) 
	\hat{\xi}_{0}^{\theta}(dx)
	\\ \nonumber
	=&
	\hat{v}_{0}^{\theta}\big(\hat{X}_{0,i}^{\theta} \big)
	-
	\big\langle\hat{\beta}_{0}^{\theta} \big\rangle. 
\end{align}
Using (2.8), (5.3), (6.4), 
(\ref{l5.3.7''}) -- (\ref{l5.3.9'}), 
we conclude 
\begin{align}\label{l5.3.11'}
	\hat{\zeta}_{n}^{\theta}(B)
	=
	\frac{1}{N} \sum_{i=1}^{N} 
	\big(
	\hat{v}_{n}^{\theta}\big(\hat{X}_{n,i}^{\theta} \big)
	-
	\big\langle\hat{\beta}_{n}^{\theta} \big\rangle
	\big)
	\delta_{\hat{X}_{n,i}^{\theta} }(B)
	=&
	\!\int_{B} \hat{v}_{n}^{\theta}(x)
	\hat{\xi}_{n}^{\theta}(dx) 
	-
	\hat{\xi}_{n}^{\theta}(B) 
	\big\langle\hat{\beta}_{n}^{\theta} \big\rangle 
	\\ \nonumber
	=&
	\hat{\beta}_{n}^{\theta}(B) 
	-
	\hat{\xi}_{n}^{\theta}(B) 
	\big\langle\hat{\beta}_{n}^{\theta} \big\rangle. 
\end{align}
Thus, (\ref{l5.3.301}), (\ref{l5.3.7''}) imply  
\begin{align*}
	\big(\hat{\xi}_{m}^{\theta}\hat{\Psi}_{m:n}^{\theta} \big)(B)
	=&
	\big(\hat{\alpha}_{m}^{\theta}(\hat{\xi}_{m}^{\theta} ) \hat{R}_{m:n}^{\theta} \big)(B)  
	+
	\big(\hat{\xi}_{m}^{\theta}\hat{S}_{m:n}^{\theta} \big)(B) 
	\\
	=&
	\big(\hat{\beta}_{m}^{\theta}\hat{R}_{m:n}^{\theta} \big)(B)  
	+
	\big(\hat{\xi}_{m}^{\theta}\hat{S}_{m:n}^{\theta} \big)(B) 
	\\
	=&
	\big(\hat{\zeta}_{m}^{\theta}\hat{R}_{m:n}^{\theta} \big)(B)  
	+
	\big(\hat{\xi}_{m}^{\theta}\hat{S}_{m:n}^{\theta} \big)(B) 
	+
	\big(\hat{\xi}_{m}^{\theta}\hat{R}_{m:n}^{\theta} \big)(B) 
	\big\langle\hat{\beta}_{m}^{\theta} \big\rangle. 
\end{align*}
Combining this with (4.4), (4.5), 
(6.5), (6.7), (6.8), (6.11), we get
\begin{align}\label{l5.3.501}
	\hat{C}_{m:n}^{\theta}(B)
	=&
	\frac{\big(\hat{\zeta}_{m}^{\theta}\hat{R}_{m:n}^{\theta} \big)(B)  
	+
	\big(\hat{\xi}_{m}^{\theta}\hat{S}_{m:n}^{\theta} \big)(B) }
	{\big\langle\hat{\xi}_{m}^{\theta}\hat{R}_{m:n}^{\theta} \big\rangle }
	+
	\big\langle\hat{\beta}_{m}^{\theta} \big\rangle 
	\frac{\big(\hat{\xi}_{m}^{\theta}\hat{R}_{m:n}^{\theta} \big)(B)}
	{\big\langle\hat{\xi}_{m}^{\theta}\hat{R}_{m:n}^{\theta} \big\rangle }
	\\ \nonumber
	=&
	\hat{H}_{m:n}^{\theta}(B)
	+
	\hat{F}_{m:n}^{\theta}(B) \big\langle\hat{\beta}_{m}^{\theta} \big\rangle. 
\end{align}
As $\hat{F}_{m:n}^{\theta}({\cal X} )=1$, 
(4.4) -- (4.6), 
(6.5) -- (6.8), 
(6.12), (6.13), (\ref{l5.3.501}) 
yield 
\begin{align*}
	\hat{A}_{m:n}^{\theta}(B) 
	+
	\hat{B}_{m:n}^{\theta}(B) 
	=
	\hat{C}_{m:n}^{\theta}(B) 
	-
	\hat{F}_{m:n}^{\theta}(B) 
	\big\langle\hat{C}_{m:n}^{\theta}\big\rangle
	=&
	\hat{H}_{m:n}^{\theta}(B) 
	-
	\hat{F}_{m:n}^{\theta}(B)
	\big\langle\hat{H}_{m:n}^{\theta} \big\rangle
	\\
	=&
	\hat{G}_{m:n}^{\theta}(B). 
\end{align*}
Hence, (6.15) holds. 

Owing to (6.1) -- (6.3), we have 
$\hat{\zeta}_{-1}^{\theta}(dx)=\hat{v}_{0}^{\theta}(x) \hat{\xi}_{-1}^{\theta}(dx)$. 
Since $\hat{F}_{-1:0}^{\theta}(dx)=\hat{\xi}_{-1}^{\theta}(dx)$, 
(4.5), (6.7), (6.8), (6.9) imply  
\begin{align*}
	\hat{H}_{-1:n}^{\theta}(B)
	=
	H_{\theta,\boldsymbol Y}^{0:n}
	\big(\hat{\xi}_{-1}^{\theta}, \hat{\zeta}_{-1}^{\theta} \big)(B)
	=&
	\frac{\big(\hat{\zeta}_{-1}^{\theta}\hat{R}_{0:n}^{\theta} \big)(B) 
	+
	\big(\hat{\xi}_{-1}^{\theta}\hat{S}_{0:n}^{\theta} \big)(B) }
	{\big\langle\hat{\xi}_{-1}^{\theta}\hat{R}_{0:n}^{\theta} \big\rangle }
	\\ \nonumber
	=&
	\frac{\int \big(
	\hat{R}_{0:n}^{\theta}(x,B) \hat{v}_{0}^{\theta}(x) + \hat{S}_{0:n}^{\theta}(x,B) 
	\big) 
	\hat{\xi}_{-1}^{\theta}(dx) }
	{\big\langle\hat{\xi}_{-1}^{\theta}\hat{R}_{0:n}^{\theta} \big\rangle }
	\\ \nonumber
	=&
	\frac{\big(\hat{F}_{-1:0}^{\theta}\hat{\Psi}_{0:n}^{\theta} \big)(B) }
	{\big\langle\hat{F}_{-1:0}^{\theta}\hat{R}_{0:n}^{\theta} \big\rangle }. 
\end{align*}
Consequently, (4.4), (4.6), 
(6.5), (6.6), (6.8), (6.10) imply 
\begin{align*}
	\hat{G}_{-1:n}^{\theta}(B)
	=&
	H_{\theta,\boldsymbol Y}^{0:n}(\hat{\xi}_{-1}^{\theta}, \hat{\zeta}_{-1}^{\theta} )(B)
	-
	F_{\theta,\boldsymbol Y}^{0:n}(\hat{\xi}_{-1}^{\theta} )(B) 
	\big\langle H_{\theta,\boldsymbol Y}^{0:n}(\hat{\xi}_{-1}^{\theta}, \hat{\zeta}_{-1}^{\theta} ) \big\rangle
	\\ \nonumber
	=&
	\frac{\int 
	\big( 
	\hat{\Psi}_{0:n}^{\theta}(x,B) 
	-
	\hat{F}_{-1:n}^{\theta}(B)
	\hat{\Psi}_{0:n}^{\theta}(x,{\cal X}) 
	\big) 
	\hat{F}_{-1:0}^{\theta}(dx) }
	{\big\langle \hat{F}_{-1:n}^{\theta}\hat{R}_{0:n}^{\theta} \big\rangle }
	\\ \nonumber
	=&
	\frac{\big(\hat{F}_{-1:0}^{\theta}\hat{\Phi}_{0:n}^{\theta} \big)(B) }
	{\big\langle\hat{F}_{-1:0}^{\theta}\hat{R}_{0:n}^{\theta} \big\rangle }. 
\end{align*}
Hence, (6.18) holds for $m=0$. 
Therefore, to complete the proof of (ii), 
it is sufficient to show (6.18) for $n\geq m>0$. 
In what follows in this part of the proof, 
we assume $n\geq m>0$.

Owing to (2.2), (4.1), (4.2), (4.3), (6.8), 
we have 
\begin{align}\label{l5.3.751}
	\hat{R}_{n-1:n}(x,B)
	=
	\int_{B}r_{\theta,\boldsymbol Y}^{n}(x'|x)\mu(dx'), 
	\;\;\;\;\; 
	\hat{S}_{n-1:n}(x,B)
	=
	\int_{B}\nabla_{\theta}r_{\theta,\boldsymbol Y}^{n}(x'|x)\mu(dx'). 
\end{align}
Consequently, (4.4), (6.4), (6.5), (6.8) imply 
\begin{align}\label{l5.3.753}
	\hat{\alpha}_{n}^{\theta}\big(\hat{F}_{n-1:n}^{\theta} \big)(B)
	=&
	\frac{\int_{B} \hat{v}_{n}^{\theta}(x) \big(\hat{\xi}_{n-1}^{\theta}\hat{R}_{n-1:n}^{\theta}\big)(dx) }
	{\big\langle \hat{\xi}_{n-1}^{\theta}\hat{R}_{n-1:n}^{\theta}\big\rangle }
	\\ \nonumber
	=&
	\frac{\int_{B} \hat{v}_{n}^{\theta}(x) 
	\left(\int r_{\theta,\boldsymbol Y}^{n}(x|x')\hat{\xi}_{n-1}^{\theta}(dx')\right)\mu(dx) }
	{\big\langle \hat{\xi}_{n-1}^{\theta}\hat{R}_{n-1:n}^{\theta}\big\rangle }. 
\end{align}
Moreover, (6.3) yields  
\begin{align}\label{l5.3.755}
	\hat{v}_{\theta,\boldsymbol Y}^{n}(x) 
	\left( \int r_{\theta,\boldsymbol Y}^{n}(x|x') \hat{\xi}_{n-1}^{\theta}(dx') \right)
	= &
	\int r_{\theta,\boldsymbol Y}^{n}(x|x') \hat{\zeta}_{n-1}^{\theta}(dx') 
	\\ \nonumber
	&
	+ 
	\int \nabla_{\theta} r_{\theta,\boldsymbol Y}^{n}(x|x') \hat{\xi}_{n-1}^{\theta}(dx').
\end{align}
Combining (\ref{l5.3.751}) -- (\ref{l5.3.755}), we get 
\begin{align}\label{l5.3.757}
	\hat{\alpha}_{n}^{\theta}\big(\hat{F}_{n-1:n}^{\theta} \big)(B)
	=&
	\frac{\int\left(\int_{B} r_{\theta,\boldsymbol Y}^{n}(x|x')\mu(dx) \right) 
	\hat{\zeta}_{n-1}^{\theta}(dx') }
	{\big\langle \hat{\xi}_{n-1}^{\theta}\hat{R}_{n-1:n}^{\theta}\big\rangle }
	\\ \nonumber
	+&
	\frac{\int\left(\int_{B} \nabla_{\theta}r_{\theta,\boldsymbol Y}^{n}(x|x')\mu(dx) \right) 
	\hat{\xi}_{n-1}^{\theta}(dx') }
	{\big\langle \hat{\xi}_{n-1}^{\theta}\hat{R}_{n-1:n}^{\theta}\big\rangle }
	\\ \nonumber
	=&
	\frac{\big(\hat{\zeta}_{n-1}^{\theta}\hat{R}_{n-1:n}^{\theta} \big)(B)
	+ \big(\hat{\xi}_{n-1}^{\theta}\hat{S}_{n-1:n}^{\theta} \big)(B)}
	{\big\langle \hat{\xi}_{n-1}^{\theta}\hat{R}_{n-1:n}^{\theta}\big\rangle}. 
\end{align}

Differentiating (\ref{l5.3.303}) in $\theta$ and using (4.1), we get 
\begin{align*}
	s_{\theta,\boldsymbol Y}^{m-1:n}(x_{m-1:n} )
	=
	r_{\theta,\boldsymbol Y}^{m:n}(x_{m:n} ) s_{\theta,\boldsymbol Y}^{m-1:m}(x_{m-1:m} ) 
	+
	s_{\theta,\boldsymbol Y}^{m:n}(x_{m:n} ) r_{\theta,\boldsymbol Y}^{m-1:m}(x_{m-1:m} ). 
\end{align*}
Consequently, (4.3), (6.8) yield 
\begin{align}\label{l5.3.23}
	\big(\hat{\xi}_{m-1}^{\theta}\hat{S}_{m-1:n}^{\theta} \big)(B) 
	=&
	\begin{aligned}[t]
	&
	\int\left(\int_{{\cal X}^{n-m-1}\times B}r_{\theta,\boldsymbol Y}^{m:n}(x_{m:n} ) 
	\mu^{n-m}(dx_{m+1:n} ) \right)
	\\
	&\cdot s_{\theta,\boldsymbol Y}^{m-1:m}(x_{m-1:m} ) 
	(\hat{\xi}_{m-1}^{\theta}\times\mu)(dx_{m-1:m} )
	\end{aligned}
	\\ \nonumber 
	&+
	\begin{aligned}[t]
	&
	\int\left(\int_{{\cal X}^{n-m-1}\times B}s_{\theta,\boldsymbol Y}^{m:n}(x_{m:n} ) 
	\mu^{n-m}(dx_{m+1:n} ) \right)
	\\
	&\cdot r_{\theta,\boldsymbol Y}^{m-1:m}(x_{m-1:m} ) 
	(\hat{\xi}_{m-1}^{\theta}\times\mu)(dx_{m-1:m} )
	\end{aligned}
	\\ \nonumber
	=&
	\int \hat{R}_{m:n}^{\theta}(x,B) 
	\big(\hat{\xi}_{m-1}^{\theta}\hat{S}_{m-1:m}^{\theta} \big)(dx)
	\\ \nonumber
	&
	+
	\int \hat{S}_{m:n}^{\theta}(x,B) 
	\big(\hat{\xi}_{m-1}^{\theta}\hat{R}_{m-1:m}^{\theta} \big)(dx). 
\end{align}
Similarly, (4.2), (6.8) imply  
\begin{align}\label{l5.3.25}
	\big(\hat{\zeta}_{m-1}^{\theta}\hat{R}_{m-1:n}^{\theta} \big)(B) 
	=&
	\begin{aligned}[t]
	&
	\int\left(\int_{{\cal X}^{n-m-1}\times B}r_{\theta,\boldsymbol Y}^{m:n}(x_{m:n} ) 
	\mu^{n-m}(dx_{m+1:n} ) \right)
	\\
	&\cdot r_{\theta,\boldsymbol Y}^{m-1:m}(x_{m-1:m} ) 
	(\hat{\xi}_{m-1}^{\theta}\times\mu)(dx_{m-1:m} )
	\end{aligned}
	\\ \nonumber
	=&
	\int \hat{R}_{m:n}^{\theta}(x,B) 
	\big(\hat{\zeta}_{m-1}^{\theta}\hat{R}_{m-1:m}^{\theta} \big)(dx). 
\end{align}
Moreover, using (4.4), (6.5), (6.8), we conclude   
\begin{align}\label{l5.3.29}
	\big(\hat{F}_{m-1:m}^{\theta}\hat{S}_{m:n}^{\theta} \big)(B)   
	=&
	\frac{\int\hat{S}_{m:n}^{\theta}(x,B) 
	\big(\hat{\xi}_{m-1}^{\theta}\hat{R}_{m-1:m}^{\theta} \big)(dx) }
	{\big\langle\hat{\xi}_{m-1}^{\theta}\hat{R}_{m-1:m}^{\theta} \big\rangle }. 
\end{align}
Similarly, relying on (\ref{l5.3.757}), we deduce 
\begin{align}\label{l5.3.27'}
	\big(\hat{\alpha}_{m}^{\theta}(\hat{F}_{m-1:m}^{\theta} ) \hat{R}_{m:n}^{\theta} \big)(B)
	=&
	\frac{\int\hat{R}_{m:n}^{\theta}(x,B) 
	\big(\hat{\zeta}_{m-1}^{\theta}\hat{R}_{m-1:m}^{\theta} \big)(dx) }
	{\big\langle\hat{\xi}_{m-1}^{\theta}\hat{R}_{m-1:m}^{\theta} \big\rangle }
	\\ \nonumber
	&+
	\frac{\int\hat{R}_{m:n}^{\theta}(x,B) 
	\big(\hat{\xi}_{m-1}^{\theta}\hat{S}_{m-1:m}^{\theta} \big)(dx) }
	{\big\langle\hat{\xi}_{m-1}^{\theta}\hat{R}_{m-1:m}^{\theta} \big\rangle }. 
\end{align}

Owing to (6.4), (6.9), we have 
\begin{align*}
	\big(\hat{F}_{m-1:m}^{\theta}\hat{\Psi}_{m:n}^{\theta} \big)(B) 
	=&
	\int \hat{R}_{m:n}^{\theta}(x,B) \hat{v}_{m}^{\theta}(x) 
	\hat{F}_{m-1:m}^{\theta}(dx)
	+
	\int \hat{S}_{m:n}^{\theta}(x,B) \hat{F}_{m-1:m}^{\theta}(dx)
	\\ \nonumber
	=&
	\big(\hat{\alpha}_{m}^{\theta}(\hat{F}_{m-1:m}^{\theta} ) \hat{R}_{m:n}^{\theta} \big)(B)
	+
	\big(\hat{F}_{m-1:m}^{\theta}\hat{S}_{m:n}^{\theta} \big)(B). 
\end{align*}
Combining this with (\ref{l5.3.29}), (\ref{l5.3.27'}), we get 
\begin{align*}
	\big(\hat{F}_{m-1:m}^{\theta}\hat{\Psi}_{m:n}^{\theta} \big)(B) 
	=&
	\frac{\big(\hat{\zeta}_{m-1}^{\theta}\hat{R}_{m-1:n}^{\theta} \big)(B) 
	+ \big(\hat{\xi}_{m-1}^{\theta}\hat{S}_{m-1:n}^{\theta} \big)(B) }
	{\big\langle\hat{\xi}_{m-1}^{\theta}\hat{R}_{m-1:m}^{\theta}\big\rangle }. 
\end{align*}
Consequently, (4.5), (6.7), (6.8), (\ref{l5.3.1'}) imply   
\begin{align*}
	\hat{H}_{m-1:n}^{\theta}(B)
	=&
	\frac{\big(\hat{\zeta}_{m-1}^{\theta}\hat{R}_{m-1:n}^{\theta} \big)(B) 
	+ \big(\hat{\xi}_{m-1}^{\theta}\hat{S}_{m-1:n}^{\theta} \big)(B)  }
	{\big\langle\hat{\xi}_{m-1}^{\theta}\hat{R}_{m-1:n}^{\theta}\big\rangle }
	=
	\frac{\big(\hat{F}_{m-1:m}^{\theta}\hat{\Psi}_{m:n}^{\theta}\big)(B) }
	{\big\langle\hat{F}_{m-1:m}^{\theta}\hat{R}_{m:n}^{\theta}\big\rangle }. 
\end{align*}
Then, (4.5), (4.6), (6.6), (6.7), (6.10) yield 
\begin{align*}
	\hat{G}_{m-1:n}^{\theta}(B)
	\!=&
	\hat{H}_{m-1:n}^{\theta}(B)
	-
	\hat{F}_{m-1:n}^{\theta}(B)
	\big\langle\hat{H}_{m-1:n}^{\theta}\big\rangle
	\\ \nonumber
	=&
	\frac{\int 
	\big( 
	\hat{\Psi}_{m:n}^{\theta}(x,B) 
	\!- 
	\hat{F}_{m-1:m}^{\theta}(B) 
	\hat{\Psi}_{m:n}^{\theta}(x,{\cal X}) 
	\big) 
	\hat{F}_{m-1:m}^{\theta}(dx) }
	{\int \hat{R}_{m:n}^{\theta}(x,{\cal X})  
	\hat{F}_{m-1:m}^{\theta}(dx) }
	\\ \nonumber
	=&
	\frac{\big(\hat{F}_{m-1:m}^{\theta}\hat{\Phi}_{m:n}^{\theta}\big)(B) }
	{\big\langle\hat{F}_{m-1:m}^{\theta}\hat{R}_{m:n}^{\theta}\big\rangle }.
\end{align*}
Hence, (6.18) holds.

\end{proof}

%% file: smc_bias_shared.tex
\usepackage{amsmath,amsfonts,amssymb}
\usepackage{paralist}

\newtheorem{assumption}{Assumption}[section]

\newtheorem*{vremark}{Remark}

\title{Bias of Particle Approximations to Optimal Filter Derivative} 

\author{Vladislav Z. B. Tadi\'{c}
\thanks{School of Mathematics, University of Bristol,
Bristol, United Kingdom 
(v.b.tadic@bristol.ac.uk). }
\and 
Arnaud Doucet
\thanks{Department $\!$of Statistics, 
$\!$University $\!$of Oxford, $\!$Oxford, $\!$United  $\!$Kingdom
(doucet@stats.ox.ac.uk). } }